\documentclass[oneside,11pt,reqno]{amsart}
\usepackage[utf8]{inputenc}
\usepackage{stmaryrd}
\usepackage{bbm}
\usepackage{mathrsfs}
\usepackage{enumerate}
\usepackage{latexsym,amsxtra}
\usepackage[dvips]{graphicx}
\usepackage{dsfont}
\usepackage{slashed}
\usepackage[all]{xy}
\usepackage{amscd,graphics}
\usepackage{amsmath,amsfonts,amsthm,amssymb}
\usepackage{latexsym,amsmath}
\usepackage{graphicx,psfrag}
\usepackage{mathabx}
\usepackage{hyperref}
\usepackage{fancyhdr}

\newtheorem{thm}{Theorem}[section]
\newtheorem{introthm}{Theorem}
\newtheorem{introcor}[introthm]{Corollary}

\newtheorem{assum}[thm]{Assumption}
\newtheorem{notation}[thm]{Notation}
\newtheorem{lem}[thm]{Lemma}

\newtheorem{pro}[thm]{Proposition}
\newtheorem{cons}[thm]{Construction}
\theoremstyle{definition}
\newtheorem{defi}[thm]{Definition}
\newtheorem{ex}[thm]{Example}
\theoremstyle{remark}
\newtheorem{rmk}{Remark}
\newtheorem{ques}[thm]{Question}

\newcommand{\Q}{\mathcal{Q}} 

\newcommand{\spinc}{\textnormal{spin}^{c}}

\newcommand{\ii}{\textnormal{i}}
\newcommand{\jj}{\textnormal{j}}
\newcommand{\kk}{\textnormal{k}}

\newcommand{\SU}{\textnormal{SU}}

\newcommand{\SW}{\textnormal{SW}}
\newcommand{\gr}{\textnormal{gr}}
\newcommand{\FSW}{\textnormal{FSW}}
\newcommand{\Ind}{\textnormal{Ind}}
\newcommand{\ind}{\textnormal{ind}}
\newcommand{\dimension}{\textnormal{dim}}
\newcommand{\HMfrom}{\widehat{HM}}
\newcommand{\HMarrow}{\overrightarrow{HM}}
\newcommand{\HMbar}{\overline{HM}}
\newcommand{\HMto}{\widecheck{HM}}
\newcommand{\Diff}{\textnormal{Diff}}
\newcommand{\Symp}{\textnormal{Symp}}
\newcommand{\QHS}{\textnormal{QHS}^{3}}

\title{The Family Seiberg-Witten Invariant and non-symplectic loops of diffeomorphisms}
\author{Jianfeng Lin}
\date{}
\address{Yau Mathematical Sciences Center\\ Tsinghua University\\ Beijing\\ China.}

\email{linjian5477@mail.tsinghua.edu.cn}

\begin{document}

\maketitle

\begin{abstract}
By extending a result of Kronheimer-Mrowka to the family setting, we prove a gluing formula for the family Seiberg-Witten invariant. This formula allows one to compute the invariant for a smooth family of 4-manifolds by cutting it open along a product family of 3-manifolds and studying the induced maps on monopole Floer (co)homology. When the cutting 3-manifold is an L-space, this formula implies a relation between the family Seiberg-Witten invariant, the Seiberg-Witten invariant of the fiber and the index of the family Dirac operator. We use this relation to calculate the Seiberg-Witten invariant of families of 4-manifolds that arise when resolving an ADE singularity using a hyperk\"ahler family of complex structures near the singularity. Several applications are obtained. First, we establish a large family of simply-connected 4-manifolds $M$ (e.g. all elliptic surfaces) such that $\pi_{1}(\textrm{Diff}(M))$ has a $\mathbb{Z}^{\infty}$-summand . For such $M$, the product $S^{2}\times M$ smoothly fibers over $S^{2}$ with fiber $M$ in infinitely many distinct ways.
Second, we show that on any closed symplectic 4-manifold that contains a smoothly embedded sphere of self-intersection $-1$ or $-2$, there is a loop of diffeomorphisms that is not homotopic to a loop of symplectormorphisms. This generalizes a previous result by Smirnov and confirms a conjecture by McDuff in dimension 4. It also provides many new examples of 4-manifolds whose space of symplectic forms has a nontrivial fundamental group or first homology group.
\end{abstract}

\section{Introduction}\label{section: introduction}
The Seiberg-Witten invariant is a powerful tool in understanding smooth structures on 4-manifolds. Its family version, as proposed by Donaldson \cite{Donaldson96} and later pursued by many authors including Ruberman \cite{RubermanObstruction,Ruberman2001},  Liu \cite{Liu2000}, Li-Liu \cite{LiLiuFamily}, Nakamura \cite{Nakamura2003} and Baraglia-Konno \cite{BaragliaKonno}, has been used to prove many essential results on smooth families of 4-manifolds  \cite{Baraglia2019,LiLiuFamily,Liu2000,RubermanObstruction,Ruberman2001}). As most gauge theoretic invariants, it is difficult to directly compute this invariant by solving the Seiberg-Witten equations. Hence it is natural to develop a gluing formula that allows one to compute this invariant via a cut-and-paste technique. For classical Seiberg-Witten invariants, such formula has been extensively studied \cite{Park2002,MorganMrowkaSzabo,Taubes2001,Froyshov2010}. In particular, using monopole Floer homology, Kronheimer-Mrowka \cite{KronheimerMrowkaMonopole} proved a gluing formula for decompositions along a general $3$-manifold. For the family Seiberg-Witten invariant, Baraglia-Konno \cite{BaragliaKonno} proved a gluing formula for connected sums, generalizing an earlier result by Ruberman \cite{Ruberman2001}. These gluing formulas are powerful in computing the (family) Seiberg-Witten invariants and they have numerous applications.


The purpose of this paper is to extend Kronheimer-Mrowka's gluing formula to the family setting and to obtain new topological applications. This version of gluing formula allows a general cutting 3-manifold and we apply it in the computation of the Seiberg-Witten invariants of two smooth families that naturally raise from resolution of singularities: the ADE family and the blown-up family. As applications, we obtain new results about  the topology of the diffeomorphism groups, the symplectormorphism groups, and the space of symplectic forms. 

It takes some works to give statements of the gluing theorems (Theorem \ref{thm: FSW gluing} and Theorem \ref{thm: family switching}). For readers' convinience, we discuss  the topological applications first.
\subsection{4-manifolds with infinite type diffeomorphism groups} The first application concerns the diffeomorphism group of a 4-manifold. There has been active researches on the topology of this topological group in recent years \cite{BaragliaKonno,RubermanObstruction,rubermanpolynomial,KronheimerMrowkaDehnTwist,watanabe2018some,baraglia2019constraints,budney2019knotted}. For example, Ruberman discovered examples of $M$ such that $\ker (\pi_{0}(\Diff(M))\rightarrow \pi_{0}(\operatorname{Homeo}(M)))$ is nontrivial \cite{RubermanObstruction}. (And he later found a further example for which the kernel contains a $\mathbb{Z}^{\infty}$-summand \cite{rubermanpolynomial}).  Budney-Gabai \cite{budney2019knotted} proved that $\pi_{0}(\Diff(S^{1}\times D^{3},\partial))$ (i.e., the mapping class group of $S^{1}\times D^{3}$ relative to the boundary) is not finitely generated. Regarding higher homotopy groups, Watanabe proved that $\pi_{i}(\Diff(D^{4},\partial))\otimes \mathbb{Q}\neq 0$ for $i=1,4,8\cdots$. For any $i\geq 1$, and Auckly-Ruberman \cite{Aucklyexoticfamily} found many other examples of closed 4-manifolds so that $\pi_{i}(\Diff(M))\rightarrow \pi_{i}(\operatorname{Homeo}(M))$ is not injective. Regarding the important example of a $K3$ surface, Baraglia-Konno \cite{baraglia2019note} proved that $\pi_{1}(\Diff(K3))\rightarrow \pi_{1}(\operatorname{Homeo}(K3))$ is not surjective. And Baraglia \cite{baraglia2021non} proved that $\pi_{1}(\Diff(K3))$ has a $\mathbb{Z}^{\infty}$ summand.  
In the current paper, we establish a large family of 4-manifolds whose diffeomorphism group has a $\mathbb{Z}^{\infty}$ summand, showing that this phenomena is very common in dimension $4$.

Our construction makes use of the generalized Dehn twist  along a smoothly embedded sphere of self-intersection $-1$ or $-2$ (i.e. a $-1$-spheres or a $-2$-sphere). Let $S$ be such a sphere and let $\nu(S)$ be its closed tubular neighborhood. Then $\nu(S)$ is diffeomorphic to the disk bundle for a complex line bundle over $S$. Let $\gamma:S^{1}\rightarrow \operatorname{Diff}(\nu(S))$ be the circle action given by fiberwise complex multiplication. 
When $S\cdot S=-1$, the restriction of $\gamma$ on $\partial\nu(S)=S^{3}$ admits a factorization 
$$
S^{1}\hookrightarrow SU(2)\rightarrow \operatorname{Diff}(\partial\nu(S))
$$
Thus we can use a null-homotopy of $S^{1}$ in $SU(2)$ to modify $\gamma$ near $\partial \nu(S)$ and get a loop $\gamma':S^{1}\rightarrow \operatorname{Diff}(\nu(S),\partial)$. Then we extend by the identity map on $M\setminus \nu(S)$ and get a loop 
$
\gamma_{S}: S^{1}\rightarrow \operatorname{Diff}(M)
$ supported near $S$. The case $S\cdot S=-2$ is similar. We have $\partial \nu(S)=\mathbb{RP}^{3}$ and the restriction of the circle action  $\gamma$ on $\partial \nu(S)$ is contained in a larger group action by $SO(3)$. This allows us to modify $\gamma^{2}$ near $\partial \nu(S)$ and get a loop $\gamma':S^{1}\rightarrow \operatorname{Diff}(\nu(S),\partial)$. Then we extend by the identity map on $M\setminus \nu(S)$ and get the loop $\gamma_{S}: S^{1}\rightarrow \operatorname{Diff}(M)$. Unlike the ordinary Dehn twist, which is a single diffeomorphism, $\gamma_{S}$ is a loop of diffeomorphisms. Let's call $\gamma_{S}$ the generalized Dehn twists along $S$.

\begin{introthm}\label{thm: infinitely generated diff}
Let $M$ be a smooth closed oriented 4-manifold that contains an infinite collection of smoothly embedded $-2$-spheres $\{S_{i}\}_{i\in \mathbb{N}}$, all representing different homology classes. Suppose  $b_{1}(M)=0$ and assume either one of the following conditions hold:
\begin{enumerate}
    \item $b^{+}_{2}(M)>1$ and   $\SW(M,\mathfrak{s})\neq 0$ for some $\spinc$-structure $\mathfrak{s}$ with $d(\mathfrak{s},M)=0$. Here $$d(\mathfrak{s},M)= \frac{c^{2}_{1}(\mathfrak{s})-3\sigma(M)-2\chi(M)}{4}$$ denotes the expected dimension of the Seiberg-Witten moduli space. Or
    \item $b^{+}_{2}(M)=1$ and there exists a $\spinc$ structure $\mathfrak{s}$ such that $d(\mathfrak{s},M)=0$ and that 
    $
    \langle c_{1}(\mathfrak{s}),[S_{i}]\rangle =0 
    $ for any $i$. (We don't impose any constraint on the Seiberg-Witten invariant.)
\end{enumerate}
Then $\pi_{1}(\operatorname{Diff}(M))$ contains a $\mathbb{Z}^{\infty}$ summand, generated by generalized Dehn twists along a subcollection $\{S_{n_{i}}\}_{i\in \mathbb{N}}$.
\end{introthm}

Recall that a space is of finite type if it is weakly homotopy equivalent to a $CW$-complex with finitely many cells in each dimension. Since the fundamental group of a finite type space is always finitely generated, Theorem \ref{thm: infinitely generated diff} implies the following corollary.

\begin{introcor}
Let $M$ be a manifold that satisfies the conditions of Theorem \ref{thm: infinitely generated diff}. Then $\Diff_{0}(X)$ (the unit component of $\Diff(X)$) is \textbf{not} of finite type. 
\end{introcor}

\begin{rmk}
It is interesting to compare Baraglia's result and Theorem \ref{thm: infinitely generated diff} with a recent result by Bustmante-Krannich-Kupers \cite{bustamante2021finiteness}, which states that the $\pi_{i}(\operatorname{Diff}(M))$ is finitely generated for any $i\geq 1$ when $M$ is smooth manifold of dimension $2n\geq 6$ with a finite fundamental group. In an earlier paper \cite{kupers2019some}, Kupers proved the same result for any dimension $\neq 4,5,7$, under the additional assumption that $M$ is $2$-connected. As mentioned in \cite{bustamante2021finiteness}, it is expected that this $2$-connectedness condition can also  be removed in odd dimensions $\geq 7$. Therefore, Baraglia's result and Theorem \ref{thm: infinitely generated diff} reflect some very special phenomena in dimension 4 which doesn't happen in higher dimensions.      
\end{rmk}
\begin{rmk} It's natural to ask whether this $\mathbb{Z}^{\infty}$ summand constructed in Theorem \ref{thm: infinitely generated diff} contains some exotic loops of diffeomorphisms (i.e. nonzero elements in  $\ker(\pi_{1}(\Diff(M))\rightarrow \pi_{1}(\operatorname{Homeo}(M)))$). For those generators (i.e. generalized Dehn twists), one can show that they are nontrivial in $\pi_{1}(\operatorname{Homeo}(M))$ by showing the total space of the corresponding fiber bundle over $S^{2}$ has a different cohomology ring as $S^{2}\times M$. (See Remark \ref{rmk: topologically nontrivial}. This argument was pointed out to the author by Kronheimer). However, the argument breaks down when one tries to show these generators are linearly independent in $\pi_{1}(\operatorname{Homeo}(M))$. Actually we will prove in Corollary \ref{cor: block} that some linear combinations are indeed contratible in the space of auotomorphisms on $M$. So it's unclear whether this $\mathbb{Z}^{\infty}$ summand contains some exotic elements. We leave this as a curious open question.
\end{rmk}

There are many ways to construct interesting examples that satisfy the conditions of Theorem \ref{thm: infinitely generated diff}. We give some of them in the following corollary.

\begin{introcor}\label{cor: construction of examples}  Let $M$ be a smooth 4-manifold that belongs to the following list.  Then $M$ satisfies the conditions of Theorem \ref{thm: infinitely generated diff} and hence $\pi_{1}(\Diff(M))$ has a $\mathbb{Z}^{\infty}$-summand.
\begin{enumerate}
    \item $M$ is an elliptic surface with $b_{1}(M)=0$. 
    \item $M$ is a hypersurface or a complete intersection in $\mathbb{CP}^{n}$ with $b^{+}_{2}(M)>1$.
    \item $M$ is a fiber sum $M'\#_{T^{2}}E$. Here $E$ is an elliptic surface with $b_{1}(E)=0$ and $M'$ is a simply-connected $4$-manifold with $b^{+}(M)>1$ and $\SW(M',\mathfrak{s})\neq 0$ for some $\mathfrak{s}$ with $d(\mathfrak{s},M)=0$. The fiber sum is applied along a regular fiber of $E$ and a smoothly embedded torus $F\subset M'$ with self-intersection $0$ and $[F]\neq 0\in H_{2}(M';\mathbb{Z})$.
    
    \item $M$ is a smooth $4$-manifold with $b_{1}(M)=0$, $b^{+}_{2}(M)>1$ and  $\SW(M,\mathfrak{s})\neq 0$ for some $\mathfrak{s}$ with $d(\mathfrak{s},M)=0$. And $M$ contains two smoothly embedded spheres $S_0,S_1$ that satisfy
$$
S_0\cdot S_0=-2,\quad S_1\cdot S_1=-2,\quad |S_0\cdot S_1|\geq 3.
$$
\item $M$ is a blow up of one of the above manifolds.
\end{enumerate} 
\end{introcor}
By the finiteness theorem proved in \cite{kupers2019some}, one can further deduce the following corollary.
\begin{introcor}\label{cor: block} Let $M$ be a simply-connected 4-manifold that satisfies the conditions of Theorem \ref{thm: infinitely generated diff}. Then the following conclusions hold.
\begin{enumerate}
    \item Let $\operatorname{haut^{id}}(M)$ be the space of continuous maps $f:M\rightarrow M$ which are homotopic to the identity. Consider the map 
    $$
    \pi_{1}(\Diff(M))\rightarrow \pi_{1}(\operatorname{haut^{id}}(M)) 
    $$
    induced by the inclusion. Then the kernel of this map contains a $\mathbb{Z}^{\infty}$ subgroup.
    \item Let $\widetilde{\Diff}(M)$ be the simplicial group of block diffeomorphisms on $M$. Consider the map 
    $$\pi_{1}(\Diff(M))\rightarrow \pi_{1}(\widetilde{\Diff}(M)) 
    $$
    induced by the inclusion. Then the kernel of this map contains a $\mathbb{Z}^{\infty}$ subgroup. Concretely, that means any element in this subgroup can be represented by a loop $\gamma:S^{1}\rightarrow \Diff(M)$ such that the map
    $$
    \widetilde{\gamma}:S^{1}\times M\rightarrow S^{1}\times M\text{ defined by } (t,x)\mapsto (t,\gamma(t)\cdot x)
    $$
    can be extended to a diffeomorphism $\widehat{\gamma}$ on $D^{2}\times M$.
    \item $S^{2}\times M$ fibers over $S^{2}$ in infinitely many distinct ways. I.e., there exist infinitely many nonisomorphic smooth fiber bundles whose fibers, bases and total spaces are diffeomorphic to $M$, $S^{2}$ and  $S^{2}\times M$ respectively. \end{enumerate}

\end{introcor}

\subsection{Nonsymplectic loops of diffeomorphisms} The next application concerns the topology of of the symplectormorphism group. This has also been a classical question in symplectic geometry. In a general dimension, the technique of pseudo-holomorphic curves is powerful. In dimension four, because of Taubes's fundamental result  ``SW=Gr'' \cite{Taubes1996,Taubes2000}, the Seiberg-Witten theory provides extra tools to study this problem. Given a closed symplectic 4-manifold $(M,\omega)$, we let $\Symp(M,\omega)$  be the sympltectormorphism group. We have a natural inclusion map 
\begin{equation}\label{eq: inclusion}
i:     \Symp(M,\omega)\rightarrow \operatorname{Diff}(M,[\omega]).
\end{equation}
Here $\operatorname{Diff}(M,[\omega])$ denotes the group diffeomorphisms that preserves the symplectic class $[\omega]\in H_{2}(M;\mathbb{R})$.  One may ask the following question:
\begin{ques}\label{ques: symplectic=diffeo}
Given $(M,\omega)$, is the map (\ref{eq: inclusion}) a homotopy equivalence?
\end{ques}
By Moser's theorem \cite{Moser}, the analogous question in dimension two always has an affirmative answer. But in higher dimensions, it is expected that $ \Symp(M,\omega)$ is significantly ``smaller'' than $\operatorname{Diff}(M,[\omega])$ and the answer to Question \ref{ques: symplectic=diffeo}) should usually (if not always) be negative. See \cite[Chapter 10]{McDuffSalamon} for relevant discussions. In dimension four, for some small rational or ruled surface (e.g., $\mathbb{CP}^{2}, \mathbb{CP}^{2}\# \overline{\mathbb{CP}^{2}}$ and $S^{2}\times S^{2}$), the topology of  $\Symp(M,\omega)$ is well-understood 
and a negative answer to Question \ref{ques: symplectic=diffeo} is known. For a general $M$, finding the homotopy type of $\Symp(M,\omega)$ seems out of reach. But one can still approach Question \ref{ques: symplectic=diffeo} by studying the induced map on homotopy groups. For any $n\geq 0$, we pick the identity element as the base point and let $$
i_{n,*}:\pi_{n}(\Symp(M,\omega))\rightarrow \pi_{n}(\Diff(M,[\omega]))    
$$
be the map induced by $i$. Note that when $n\geq 1$, we have $\pi_{n}(\Diff(M,[\omega]))=\pi_{n}(\Diff(X))$ because $\Diff(M,[\omega])$ contains the unit component of $\operatorname{Diff}(M)$.  

Among $\pi_{n}(\Symp(M,\omega))$ for various $n$, the most interesting ones are $\pi_0(\Symp(M,\omega))$ and $\pi_1(\Symp(M,\omega))$. The former is  the so-called symplectic mapping class group, and the later is related to Hofer's geometry \cite{Polterovich} and quantum cohomology \cite{Seidel97}. So one may start by asking about the surjetivity/injectivity of $i_{n,*}$ when $n=0,1$. It's an open question whether $i_{0,*}$ is always surjective. For the injectivity of $i_{0,*}$, Seidel proved a remarkable result in 
his thesis \cite{Seidel97}. Let $L$ be an embedded Lagrangian sphere in $M$ and let $\lambda_{L}$ be a Dehn-Seidel twist along $L$ (see \cite{Seidel08}). Under mild assumptions on the homology of $M$, Seidel proved that $\Lambda^{2}_{L}$, which is smoothly isotopic to the identity, is not symplectically isotopic to the identity. This gives a large family of examples where $i_{0,*}$ is not injective. There has been many subsequent works that prove various generalization and variation of Seidel's result. In particular, Kronheimer \cite{KronheimerSymplectic} established nontrivial higher dimensional families of sympltectic forms, using the family Seiberg-Witten invariants. By a refinement of Kronheimer's argument, Smirnov \cite{smirnov2020flops} established a large family of K\"ahler surfaces for which the map 
$
i_{1,*}$ is not surjective. In other words, on such $M$, there exists a loop of diffeomorphisms that is not smoothly isotopic to any loop of symplectormorphisms. (Previously such phenomena was only known for rational or ruled surfaces.) Both Kronheimer and Smirnov's proofs rely on computations of the family Seiberg-Witten invariants via algebraic geometric methods (counting rational holmorphic curves in a K\"ahler family). Using the gluing theorem (Theorem \ref{thm: FSW gluing}), we can compute this invariant for many smooth families. And we show that any nontrivial linear combinations of generalized Dehn twists along homological distinct spheres are non-symplectic. 


To give the precise statement, we let $K$ be the canonical class of a symplectic 4-manifold $(M,\omega)$ and we consider the following sets 
\begin{equation}\label{eq: the set of roots}
\begin{split}
R^{1}_{M}:= \{\alpha\in H_{2}(M;\mathbb{Z})&\mid  \alpha\cdot\alpha=-1,\ \langle K, \alpha\rangle=-1,\\ &\alpha \text{ is represented by a smoothly embedded sphere}\}\\
R^{2}_{M}:= \{\alpha\in H_{2}(M;\mathbb{Z})&\mid  \alpha\cdot\alpha=-2,\ \langle K, \alpha\rangle=0,\\ &\alpha \text{ is represented by a smoothly embedded sphere}\}/\pm 1\\
&R_{M}:=R^{1}_{M}\cup R^{2}_{M}
\end{split}
\end{equation}


\begin{introthm}\label{thm: linearly independence} Let $(M,\omega)$ be a closed symplectic 4-manifold. We pick an embedded sphere $S_{\alpha}$ for  each $\alpha\in R_{M}$. Then the generalized Dehn twists $\{\gamma_{S_{\alpha}}\}$ generate a free abelian subgroup 
$$ 
\mathbb{Z}^{R_{M}}\subset \pi_{1}(\rm{Diff}(M)),
$$
which intersects trivially with the image of $i_{1,*}$. 
\end{introthm}


\begin{introcor}
\label{thm: non-symplectic loop}
Let $(M,\omega)$ be a closed symplectic 4-manifold that contains a smoothly embedded sphere $S$ of self-intersection $-1$ or $-2$. Then the map $
i_{1,*}$
is not surjective. Therefore, the map $i$ is \textbf{not} a homotopy equivalence.
\end{introcor}



\begin{rmk}
McDuff made a conjecture about the symplectomorphism group of a blow up in all even dimensions (see \cite{McDuff08}, second paragraph before Proposition 1.7). And she made progresses in the 4-dimensional case. Theorem \ref{thm: linearly independence} confirms the 4-dimensional case of this conjecture in full generality. The higher-dimensional case remains open.
\end{rmk}

We can also reinterpret Theorem \ref{thm: non-symplectic loop} as results on the space of symplectic forms. Let $\mathbb{S}(M)$ be the space of all symplectic forms on $M$ and let $\mathbb{S}_{a}(M)$ be its subspace consisting of sympletctic forms in a fixed cohomology class $a\in H^{2}(M;\mathbb{R})$. These  spaces have also been of interest. For example, it is an open question whether $\mathbb{S}_{a}(M)$ is connected for all symplectic 4-manifold $M$ (i.e., whether cohomologous symplectic forms are always isotopic in dimension 4).  One way to study the topology of $\mathbb{S}_{a}(M)$ is via the embedding 
\begin{equation}\label{eq: symplectic to almost symplectic}
\mathbb{S}_{a}(M)\hookrightarrow\mathcal{S}(M):=\{\alpha\in \Omega^{2}(M)\mid \alpha \wedge \alpha \text{ is nowhere vanishing}\}
\end{equation}
 into the space of \emph{almost symplectic forms}, whose homotopy type can be understood by the classical  obstruction theory. Note that a theorem by Gromov \cite{Gromov69} states that (\ref{eq: symplectic to almost symplectic}) is always a  homotopy equivalence when $M$ is an open manifold.

Using a long exact sequence relating the homotopy groups of $\Diff(M)$, $\Symp(M,\omega)$ and $\mathbb{S}_{[\omega]}(M)$, Theorem \ref{thm: non-symplectic loop} implies the following corollary.

\begin{introcor}\label{cor: noncontractible loop of symplectic forms} The  following results hold for any closed symplectic 4-manifold $(M,\omega)$.
\begin{enumerate}
\item There is a subgroup $\mathbb{Z}^{R_{M}}$ in $H_{1}(\mathbb{S}_{[\omega]}(M);\mathbb{Z})$.
\item There is a subgroup $\mathbb{Z}^{R^{2}_{M}}$ in $\ker(H_{1}(\mathbb{S}_{[\omega]}(M);\mathbb{Z})\rightarrow H_{1}(\mathcal{S}(M);\mathbb{Z}))$.
\item Suppose $b^{+}_{2}(M)\neq 2, 3$. Then there is a subgroup  
$\mathbb{Z}^{R^{1}_{M}}$ in $H_{1}(\mathbb{S}(M);\mathbb{Z})$.
\item Suppose $b^{+}_{2}(M)\neq 3$. Then there is a subgroup  $\mathbb{Z}^{R^{1}_{M}}$ in $\pi_{1}(\mathbb{S}(M);\mathbb{Z})$.
\end{enumerate}

\end{introcor}



\subsection{Fiberwise symplectic structures on the ADE family and the blown-up family}
Proof of Theorem \ref{thm: non-symplectic loop} relies on the study of fiberwise symplectic structures on two families of 4-manifolds: the blown-up family and the ADE family. Both appears  naturally when one studies surface singularities. 

We start by explaining our terminologies. Let 
$$
M\hookrightarrow \widetilde{M}\xrightarrow{p}\Q
$$
be a smooth fiber bundle whose fiber $M$ is a smooth 4-manifold and whose base $\Q$ is another smooth, connected manifold.  We fix a base point $b_{0}\in \Q$ and fix a diffeomorphism $M\cong p^{-1}(b_{0})$.  We denote this bundle by $\widetilde{M}/\Q$ and treat it as a smooth family of 4-manifolds parameterized by $\Q$. For simplicity, we always assume that $\widetilde{M}/\Q$ is a homology product, i.e., $\pi_{1}(\Q,b_{0})$ acts trivially on $H^{*}(M;\mathbb{Z})$.

By a fiberwise symplectic structure on $\widetilde{M}/\Q$, we mean a smooth family $\{\omega_{b}\}_{b\in \Q}$ of symplectic forms on the fibers $\{p^{-1}(b)\}_{b\in \Q}$. A fiberwise symplectic structure induces a characteristic map
\begin{equation}\label{eq: characteristic map}
\rho: \Q\rightarrow \{a\in H^{2}(M;\mathbb{R})\mid a\cdot a>0\}\simeq S^{b^{+}_{2}(M)}
\end{equation}
defined by $\rho(b)=[\omega_{b}]$.

Now we introduce the blown-up family. Let $M'$ be a smooth 4-manifold. We pick a point $p\in M'$ and identify a  neighborhood of $U$ of $p$ with the unit disk $D(\mathbb{H})$ in $\mathbb{H}$ (the space of quoternions). Consider the space 
$$
\Q:=\{b \in \mathbb{H}\mid b^{2}=-1\}=\{x\ii+y\jj+z\kk\mid x^{2}+y^{2}+z^{2}=1\}\cong S^{2}
$$
Via the right multiplication, each point $b\in \Q$ gives a complex structure $J_{b}$ on $U$. We take the trivial fiber bundle $\Q\times M'$ and blow up the fiber $\{b\}\times M'$ at $(b,p)$ using $J_{b}$. This gives a smooth family
$$
M=M'\# \overline{\mathbb{CP}^{2}}\hookrightarrow \widetilde{M}\rightarrow \Q,
$$
called the blown-up family. We use $S$ to denote the exceptional divisor in $M$ and use $\operatorname{PD}[S]$ to denote its Poincar\'e dual. By computing the family Seiberg-Witten invariant of the blown-up family, we prove the following result.
\begin{introthm}\label{thm: blown-up family has fiberwise symplectic structure}
Let $\widetilde{M}/\Q$ be a blown-up family that admits a fiberwise symplectic structure $\{\omega_{b}\}_{b\in \Q}$ with canonical class $K\in H^{2}(M;\mathbb{Z})$. Suppose $\langle K,[S]\rangle =-1$. 
Then both of the following conclusions hold:
\begin{enumerate}[(i)]
    \item $b^{+}_{2}(M)=3$;
    \item The characteristic map (\ref{eq: characteristic map}) has mapping degree $\pm 1$.
\end{enumerate}
\end{introthm}
By \cite[Theorem 1]{Li99}, the condition $\langle K,[S]\rangle=-1$ is automatically satisfied unless $M$ is rational or ruled. So we get the following corollary.
\begin{introcor}
Let $\widetilde{M}/\Q$ be a blown-up family. Suppose the fiber $M$ satisfies $b^{+}_{2}(M)\neq 1, 3$. Then $\widetilde{M}/\Q$ does \textbf{not} admit a fiberwise symplectic structure.
\end{introcor}
Next, we introduce the ADE family. Recall that an ADE singularity is a discrete surface singularity that locally models on the quotient $\mathbb{H}/\Gamma$, where $\Gamma$ is a finite subgroup of $Sp(1)$. Up to conjugacy, such groups one-one correspond to Dynkin diagrams of type A, D or E.  

Given any complex structure on $\mathbb{H}$ that is compatible with the action of $\Gamma$, the singular variety $\mathbb{H}/\Gamma$ has a minimal resolution. Therefore, we can take the product family
$
\Q\times D(\mathbb{H})/\Gamma
$
and resolve the fiber $\{b\}\times D(\mathbb{H})/\Gamma$ using the complex structure $J_{b}$. By doing this, we get a smooth (but not algebraic) family 
\begin{equation}\label{eq: ADE with boundary}
M_{1}\rightarrow \widetilde{M}_{1}\rightarrow \Q    
\end{equation}
whose fiber $M_{1}$ is a complex 4-manifold with boundary $Y=S^{3}/\Gamma$. Let $M$ be any closed 4-manifold such that $M_{1}$ can be embedded smoothly. (e.g., $M$ is the resolution of an algebraic surface with a single ADE singularity.) Then $M$ admits a smooth decomposition $M_{1}\cup _{Y}M_{2}$. By gluing the family (\ref{eq: ADE with boundary}) with the product family $\Q\times M_{2}$, we obtain a smooth family
$$
M\hookrightarrow \widetilde{M}=\widetilde{M}_{1}\cup (\Q\times M_{2})\rightarrow \Q,
$$
called the ADE family. (Note that for an $A_{1}$-singularity, this construction works for any $M$ that contains a smoothly embedded $-2$-sphere $S$. In this case, $M_{1}$ is just a tubular neighborhood of $S$.) 

\begin{introthm}\label{thm: ADE family doesn't carry symplectic structure} Let $\widetilde{M}/\Q$ be an ADE family whose fiber is minimal. Suppose there is a fiberwise symplectic structure $\{\omega_{b}\}_{b\in \Q}$ on $\widetilde{M}/\Q$. Then at least one of following two conclusion holds:
\begin{enumerate}
\item For homology class $a\in H_{2}(M_{0};\mathbb{Z})$ with $a\cdot a=-2$ (such classes one-one correspond to roots of the corresponding Dynkin diagram), the function 
    $
    \Q\rightarrow \mathbb{R}$ defined by $b\mapsto \langle [\omega_{b}],a\rangle 
    $
    takes both positive and negative values, or 
    \item $b^{+}_{2}(M)=3$, $K$ is torsion, and characteristic map (\ref{eq: characteristic map})
    has mapping degree $\pm 1$.   
\end{enumerate}
\end{introthm}

\begin{introcor}
Any ADE family does \textbf{not} admit a fiberwise minimal symplectic structure in a constant cohomology class.
\end{introcor}

\subsection{Statement of the gluing theorem} Now we state our gluing theorem. We describe our setting as follows: Let $$M\hookrightarrow \widetilde{M}\xrightarrow{p} \Q$$
be a smooth family of closed 4-manifolds. Again, we assume that  $\widetilde{M}/\Q$ is a homology product. 
We further assume that $\widetilde{M}/\Q$ admits a smooth decomposition
\begin{equation}\label{eq: decomposition1}
\widetilde{M}=\widetilde{M}_{1}\cup_{\Q\times Y} (\Q\times M_{2}),    
\end{equation}
where $\widetilde{M}_{0}/\Q$ is a (possibly nontrivial) family whose fiber $M_{1}$ is a 4-manifolds with boundary $Y$, and $M_{2}$ is a 4-manifold with boundary $-Y$. 
By removing a small disk from $M_{2}$ and each fiber of $\widetilde{M}_{1}$, we obtain 
a single cobordism 
$$
W_{12}=M_{2}\setminus \mathring{D}^{4}: Y\rightarrow S^{3}
$$
and a family of cobordisms 
$$
\widetilde{W}_{01}= \widetilde{M}_{1}\setminus (\Q\times \mathring{D}^{4}): S^{3}\rightarrow Y.
$$
We pick a family $\spinc$ structure $\mathfrak{s}_{\widetilde{W}_{01}}$ on $\widetilde{W}_{01}/\Q$ and a $\spinc$ structure $\mathfrak{s}_{W_{12}}$ on $W_{12}$ such that $\mathfrak{s}_{\widetilde{W}_{01}}|_{\Q\times Y}$ is the pull-back of $\mathfrak{s}_{W_{12}}|_{Y}$. We let $\mathbb{S}$ be the set of family $\spinc$-structures $\widetilde{\mathfrak{s}}$ on $\widetilde{M}/\Q$ that satisfies 
$$
\widetilde{\mathfrak{s}}|_{\widetilde{W}_{01}}=\mathfrak{s}_{\widetilde{W}_{01}},\quad \widetilde{\mathfrak{s}}|_{W_{12}}=\mathfrak{s}_{W_{12}}\quad \text{and}\quad  d(\widetilde{\mathfrak{s}},\widetilde{M})=0.
$$
Here $$d(\widetilde{\mathfrak{s}},\widetilde{M}):=\frac{c_{1}(\widetilde{\mathfrak{s}}|_{M})^{2}[M]-2\chi(M)-3\sigma(M)}{4}+\operatorname{dim}(\Q)$$ denotes the expected dimension of the parametrized moduli space. Note that when $Y$ is a rational homology 3-sphere ($\QHS$), $\mathbb{S}$ contains at most one element.

\begin{introthm}\label{thm: FSW gluing} We have the following results regarding the family Seiberg-Witten invariant. 
\begin{enumerate}
    \item Suppose $b^{+}_{2}(M_{2})>1$. Then we have 
    $$
\sum_{\widetilde{\mathfrak{s}}\in \mathbb{S}}\FSW(\widetilde{M},\widetilde{\mathfrak{s}})= \langle\HMfrom_{*}(\widetilde{W}_{01},\mathfrak{s}_{\widetilde{W}_{01}})(\widehat{1})\ ,\  \HMarrow^{*}(W_{12},\mathfrak{s}_{W_{12}})(\widecheck{1})\rangle $$
\item Suppose $b^{+}_{2}(M)>b^{+}_{2}(M_{2})=1$. Then we have 
    $$
\sum_{\widetilde{\mathfrak{s}}\in \mathbb{S}}\FSW(\widetilde{M},\widetilde{\mathfrak{s}})= \langle\HMfrom_{*}(\widetilde{W}_{01},\mathfrak{s}_{\widetilde{W}_{01}})(\widehat{1})\ ,\  \HMarrow^{*}(W_{12},\mathfrak{s}_{W_{12}})_{\pm}(\widecheck{1})\rangle .$$
\item Suppose $b^{+}_{2}(M)=b^{+}_{2}(M_{2})=1$. Then we have 
    $$
\sum_{\widetilde{\mathfrak{s}}\in \mathbb{S}}\FSW_{\pm}(\widetilde{M},\widetilde{\mathfrak{s}})= \langle\HMfrom_{*}(\widetilde{W}_{01},\mathfrak{s}_{\widetilde{W}_{01}})(\widehat{1})\ ,\  \HMarrow^{*}(W_{12},\mathfrak{s}_{W_{12}})_{\pm}(\widecheck{1})\rangle .$$
\end{enumerate}
Here $\FSW(\widetilde{M},\widetilde{\mathfrak{s}})$ and  $\FSW_{\pm}(\widetilde{M},\widetilde{\mathfrak{s}})$ denote the family Seiberg-Witten invariants for canonical chambers (see Section \ref{section: family SW}); $\widecheck{1}$ and $\widehat{1}$ denote canonical generators of the monopole Floer (co)homology of $S^{3}$; $\HMarrow^{*}(-)$ and $\HMfrom_{*}(-)$ denote two different versions of cobordism-induced maps; and $\langle -,-\rangle$ denotes the canonical pairing between the monopole Floer homology of $Y$ and the monopole Floer cohomology of $Y$. 
\end{introthm} 
Theorem \ref{thm: FSW gluing} is a generalization of Kronheimer-Mrowka's result \cite{KronheimerMrowkaMonopole} to the family setting and the proof is adapted from their proof. When the cutting 3-manifold $Y$ is an L-space, one can recover $\HMfrom_{*}(\widetilde{W}_{01},\mathfrak{s}_{\widetilde{W}_{01}})(\widehat{1})$ by studying the wall-crossing of the family Dirac opeartors and recover $\HMarrow^{*}(W_{12},\mathfrak{s}_{W_{12}})(\widecheck{1})$ from the Seiberg-Witten invariant of the fiber $M$. In this case,  Theorem \ref{thm: FSW gluing} can be used to deduce a relation between the Seiberg-Witten invariants of the family $\widetilde{M}/\Q$ and the fiber $M$. 
\begin{introthm}\label{thm: family switching}
Let $\widetilde{M}/\Q$ be a family of closed 4-manifolds with a connected, even dimensional base $\Q$. Suppose $\pi_{1}(\Q)$ acts trivially on $H_{*}(M;\mathbb{Z})$ and suppose $\widetilde{M}$ has a decomposition of the form (\ref{eq: decomposition1}) with $b_{1}(M_1)=0$ and $Y$ being an $L$-space.  Let $\widetilde{\mathfrak{s}}$ be a family $spin^{c}$ structure on $\widetilde{M}$ and let $\mathfrak{s}$ be a $spin^{c}$ structure on the fiber $M$ which satisfy $\widetilde{\mathfrak{s}}|_{M_{2}}=\mathfrak{s}|_{M_{2}}$
and 
$
d(\mathfrak{s})\geq d(\widetilde{\mathfrak{s}},\widetilde{M})=0. 
$ Then we have the following results:

\begin{enumerate}
    \item Suppose $b^{+}_{2}(M_{1})=0$ and $b^{+}_{2}(M_{2})>1$. Then 
    \begin{equation*}
    \FSW(\widetilde{M},\widetilde{\mathfrak{s}})=\langle c_{\frac{\operatorname{dim}(\Q)}{2}}(-(\operatorname{Ind}(\slashed{D}^{+}(\widetilde{M},\widetilde{\mathfrak{s}}))\otimes L^{-1})),[\Q]\rangle \cdot \SW(M,\mathfrak{s}).
\end{equation*}
\item Suppose $b^{+}_{2}(M_{1})=0$ and $b^{+}_{2}(M_{2})=1$. Then \begin{equation*}
    \FSW_{\pm}(\widetilde{M},\widetilde{\mathfrak{s}})=\langle c_{\frac{\operatorname{dim}(\Q)}{2}}(-(\operatorname{Ind}(\slashed{D}^{+}(\widetilde{M},\widetilde{\mathfrak{s}}))\otimes L^{-1})),[\Q]\rangle \cdot \SW_{\pm}(M,\mathfrak{s}).
\end{equation*}
\item Suppose $b^{+}_{2}(M_{1})>0$ and $b^{+}_{2}(M)>1$. Then $\FSW(\widetilde{M},\widetilde{\mathfrak{s}})=0.$ 
\item Suppose $b^{+}_{2}(M_{1})>0$ and $b^{+}_{2}(M)=1$. Then  $\FSW_{\pm}(\widetilde{M},\widetilde{\mathfrak{s}})=0$.
\end{enumerate}
Here $L$ is the unique complex line bundle over $\Q$ such that $(\widetilde{\mathfrak{s}}|_{\Q\times M_{2}})\otimes p^{*}(L^{-1})$ 
is a pull back of $\widetilde{\mathfrak{s}}|_{M_{2}}$. And $\operatorname{Ind}(\slashed{D}^{+}(\widetilde{M},\widetilde{\mathfrak{s}}))$ denotes the index bundle of the family Dirac operator.
\end{introthm}
We note that when $M_{1}$ is a tubular neighborhood of an embedded sphere with negative self-intersection number (which is the case in many of our applications), Theorem \ref{thm: family switching} can be recovered from Liu's family blow-up formula and family switching formula \cite{Liu2000,liu2003,liu2003switching}.

\subsection{Organization of the paper} In Section 2, we review the definition of the family Seiberg-Witten invariants and recall some of its important properties. In Section 3, we give a very quick review on the monopole Floer homology, following \cite{KronheimerMrowkaMonopole}. In Section 4, maps induced by family of cobordisms are defined. The gluing formula (Theorem \ref{thm: FSW gluing}) is proved. In Section 6 and 7, we stdudy the ADE family and the blown-up family respectively. The family Seiberg-Witten invariants for these families are computed in Section 8. Theorem \ref{thm: blown-up family has fiberwise symplectic structure} and Theorem \ref{thm: ADE family doesn't carry symplectic structure} are also proved in Section 8. Finally in Section 9, we prove Theorem \ref{thm: infinitely generated diff}, Theorem \ref{thm: linearly independence} and their corollaries. Readers who are mainly interested in the topological applications could assume Theorem \ref{thm: FSW gluing} and skip Section 2 through Section 5.
\\
\\
\textbf{Acknowledgment
} We would like to thank Hokuto Konno, Sander Kupers, Tianjun Li, Daniel Ruberman, Weiwei Wu for many enlightening discussions. We are very grateful to Gleb Smirnov for explaining his work in \cite{smirnov2020flops} and for many useful comments. And we especially like to thank Peter Kronheimer for suggesting to study the ADE families and for explaining the argument that the $A_1$-family is topologically nontrivial.   

\section{Recollection of the family Seiberg-Witten theory}\label{section: family SW}
\subsection{The family Seiberg-Witten invariants} In this section, we briefly recall the definition and basic properties of the family Seiberg-Witten invariant. For a more detailed exposition, we refer to Li-Liu's paper \cite{LiLiuFamily}, which carefully sets up the family Seiberg-Witten theory and proves a general wall-crossing formula.

Let $M$ be a connected, closed, smooth oriented 4-manifold. Consider a smooth fiber bundle $p:\widetilde{M}\rightarrow \Q$ whose fibers are diffeomorphic to $M$ and whose base is another connected, closed, oriented smooth manifold $\Q$. We use $M_{b}$ to denote the fiber over $b$. We fix a base point $b_{0}\in Q$ and set $M=M_{b_{0}}$. We treat this fiber bundle as a family of 4-manifolds over $\Q$ and denote it by $\widetilde{M}/\Q$. By taking monodromies, one obtains an action of $\pi_{1}(\Q,b_{0})$ on $H_{*}(M;\mathbb{Z})$

\begin{defi}\label{defi: homology product(closed)}
 We say $\widetilde{M}/\Q$ is a \emph{homology product} if the monodromy action is trivial.
\end{defi}

For a homology product, one has a canonical identification between homology of a general fiber $M_{b}$ with the homology of $M$ . From now on, we will assume all families we work with are homology products. (This is automatically satisfied if $\Q$ is simply-connected.) Let $$T(\widetilde{M}/\Q)=\ker(p_{*}:T\widetilde{M}\rightarrow TQ)$$
be the tangent bundle in the fiber direction. We equip $T(\widetilde{M}/\Q)$ with a smooth Riemannian metric $g_{\widetilde{M}/\Q}=\{g_{b}\}_{b\in \Q}$. A family $\spinc$ structure $\widetilde{\mathfrak{s}}$ is a lift of the corresponding frame bundle to an  principal bundle $P$ with fiber 
$$
\textnormal{Spin}^{c}(4)=(\SU(2)\times \SU(2)\times S^{1})/\pm 1
.$$ Associated to $P$, we have a rank-4 Hermitian bundle $S$ over $\widetilde{M}$ (the ``spinor bundle''), which is canonically decomposed as the direct sum $S^{+}\oplus S^{-}$. Let $\Lambda^{*}T^{*}(\widetilde{M}/\Q)$ be the bundle (over $\widetilde{M}$) of vertical differential forms. The Clifford multiplication gives an action on $S$ by vertical differential forms, which restricts to $$\rho:\Lambda^{\text{even}}(T^{*}(\widetilde{M}/\Q))\otimes \mathbb{C}\rightarrow \operatorname{End}(S^+).$$ 
We denote the restriction of $P,\ g_{\widetilde{M}/\Q},\  S^{\pm}$ to the fiber $M_{b}$ by $P_{b}$, $g_{b}$ and $S^{\pm}_{b}$ respectively. We use $\mathfrak{s}$ to denote the restriction of $\widetilde{\mathfrak{s}}$ to $M$. 

Now we consider the configuration space 
\begin{equation*}
\begin{split}
\mathcal{C}:=\{(b,A_{b},\phi_{b})\mid 
b \text{ is a point in }\Q.\ \text{$A_{b}$ is a $\spinc$-connection on $P_{b}$},&\\\ \text{and}\ \phi_{b} \text{ is a smooth section of }S^{+}_{b}.&\}    
\end{split}
\end{equation*}
This is a fiber bundle over $\Q$ whose fiber is a Fr\'echet space. The (unperturbed) Seiberg-Witten equations for a point $(b,A_{b},\phi_{b})\in \mathcal{C}$ can be written as
\begin{equation}\label{eq: unperturbed SW}
\begin{split}
    \rho(F^{+}_{A_{b}^{t}})+(\phi_{b}\phi_{b}^{*})_{0}&=0\\
    \slashed{D}^{+}_{A_{b}}\phi_{b}&=0.
\end{split}
\end{equation}
Here $A_{b}^{t}$ denotes the induced connection on the determinant bundle $\operatorname{det}(S^{+}_{b})$ and $F^{+}_{A_{b}^{t}}\in \Omega^{+}_{2}(M_{b};i\mathbb{R})$ denotes the self-dual part of its curvature. $\slashed{D}^{+}_{A_{b}}$ denotes the Dirac operator. The endomorphism $(\phi_{b}\phi_{b}^{*})_{0}\in \operatorname{End}(S^+_{b})$ is defined by the formula 
$$
\psi\mapsto \langle\psi,\phi_{b}\rangle\cdot \phi_{b}-\frac{|\phi_{b}|^2}{2}\psi.
$$
A solution of (\ref{eq: unperturbed SW}) is called reducible if $\phi_{b}=0$. To achieve transversality and to avoid reducible solutions, one needs to introduce a perturbation on the Seiberg-Witten equations as follows: Take a smooth family $\widetilde{\omega}=\{\omega_{b}\}_{b\in \Q}$ of self-dual 2-forms on the fibers $\{M_{b}\}_{b\in \Q}$. The $\widetilde{\omega}$-perturbed Seiberg-Witten equations are written as 
\begin{equation}\label{eq: perturbed SW}
\begin{split}
    \rho(F^{+}_{A_{b}^{t}})+(\phi_{b}\phi_{b}^{*})_{0}&=\rho(\ii\omega_{b})\\
    \slashed{D}^{+}_{A_{b}}\phi_{b}&=0.
\end{split}
\end{equation}
There is a condition on $(g_{\widetilde{M}/\Q},\widetilde{\omega})$ which can ensure that  (\ref{eq: perturbed SW}) has no reducible solutions. To state this condition, we use 
$I(M_{b})$ to denote $H^{2}(M_{b};\mathbb{R})$. Then $I(M_{b})$ can be identified with the space of harmonic 2-forms on $M_{b}$ with respect to the metric $g_{b}$. There is a metric-dependent decomposition
\begin{equation}\label{eq: decompostion of IM}
I(M_{b})=I^{+}(M_{b},g_{b})\oplus I^{-}(M_{b},g_{b})
\end{equation}
where $I^{+}(M_{b},g_{b})$ and $I^{-}(M_{b},g_{b})$ denote the space of self-dual and anti-self-dual harmonic 2-forms
respectively. This decomposition induces a projection $$P_{I^{+}}:I(M_{b})\rightarrow I^{+}(M_{b},g_{b}).$$ 
We have another projection map 
$$
P_{I}: \Omega^{2}(M_{b})\rightarrow I(M_{b})
$$
defined by the following condition
$$
\int_{M_{b}} P_{I}(\omega)\wedge \kappa =\int_{M_{b}}\omega\wedge \kappa \text{ for any harmonic 2-form } \kappa. 
$$
Note that $P_{I}(\omega)\in I^{+}(M_{b},g_{b})$ if $\omega$ is self-dual.
 \begin{defi}\label{defi: admissible pair closed 4-manifolds}
The pair $(g_{\widetilde{M}/\Q},\widetilde{\omega})$ is called \emph{admissible} if for any $b\in \Q$, we have
$$
P_{I}(\omega_{b})+2\pi P_{I^{+}} c_{1}(\mathfrak{s})\neq 0.
$$
\end{defi}

Suppose $(g_{\widetilde{M}/\Q},\widetilde{\omega})$ is admissible. Then the Chern-Weil formula ensures that the curvature equation
$
F^{+}_{A^{t}_{b}}=i\omega_{b}
$ has no solution for any $b$. Hence (\ref{eq: perturbed SW}) has no reducible solutions.  Consider the Seiberg-Witten moduli space $$\mathcal{M}(g_{\widetilde{M}/\Q},\widetilde{\omega}):=\{(b,A_{b},\phi_{b})\in \mathcal{C} \text{ that solves (\ref{eq: perturbed SW})}\}/\text{gauge transformations}.$$

We say $(g_{\widetilde{M}/\Q},\widetilde{\omega})$ is  regular if it  is admissible and if the moduli space $\mathcal{M}(g_{\widetilde{M}/\Q},\widetilde{\omega})$ is regular at each point. (That means the linearized Seiberg-Witten equations combined with the gauge fixing equation is form a surjective Fredholm operator. See \cite[Page 6]{LiLiuFamily} for a precise definition). By the Sard-Smale theorem, this holds for a generic pair $(g_{\widetilde{M}/\Q},\widetilde{\omega})$. In this case, the moduli space is a compact, smooth manifold of dimension 
$$
d(\widetilde{M},\widetilde{\mathfrak{s}}):=\frac{c^{2}_{1}(S^{+}_{b})[M]+3\sigma(M)+2\chi(M)}{4}+\dimension\Q.
$$
Since the bundle $\widetilde{M}/Q$ is a homology product, the moduli space $\mathcal{M}(g_{\widetilde{M}/\Q},\widetilde{\omega})$ is orientable and an orientation is specified by an orientation of $\Q$ together with a homological orientation on $M$ (i.e. an orientation of $H^{1}(M;\mathbb{R})\oplus H^{2}_{+}(M;\mathbb{R})$).

Although the family Seiberg-Witten invariant can be defined in general. We will focus on the case $d(\widetilde{M},\widetilde{\mathfrak{s}})=0$ in this paper. Then we use $\#\mathcal{M}(g_{\widetilde{M}/\Q},\widetilde{\omega})$ to denote the number of points in $\mathcal{M}(g_{\widetilde{M}/\Q},\widetilde{\omega})$, counted with sign. 

Just like the classical Seiberg-Witten theory, the family Seiberg-Witten theory has a wall-crossing phenomena. Namely, the space of all admissible pairs $(g_{\widetilde{M}/\Q},\widetilde{\omega})$ has a certain chamber structure and the number  $\#\mathcal{M}(g_{\widetilde{M}/\Q},\widetilde{\omega})$ will change if the pair goes from one chamber to another. These chambers correspond to elements in the cohomotopy group of $\Q$. To give a precise description, we consider the positive cone
$$
V^{+}(M):=\{\mathfrak{a}\in H^{2}(M;\mathbb{R}) \mid  (\mathfrak{a}\cup \mathfrak{a})[M] > 0\}.
$$
We use $[\Q,V^{+}(M)]$ to denote the set of homotopy classes of maps from the base $\Q$ to $V^{+}(M)$. Given any admissible pair $(g_{\widetilde{M}/\Q},\widetilde{\omega})$, we have a map from $\Q$ to $V^{+}(M)$ defined by
\begin{equation}\label{eq: chamber map}
b\mapsto  P_{I}(\omega_{b})+2\pi P_{I^{+}} c_{1}(\mathfrak{s}) \in I^{+}(M_{b},g_{b})-\{0\}\subset V^{+}(M).
\end{equation}

This map represents an element $\xi(g_{\widetilde{M}/\Q},\widetilde{\omega})\in [\Q,V^{+}(M)]$.
\begin{lem}
Any element in $\xi\in [\Q,V^{+}(M)]$ can be realized by some regular pair $(g_{\widetilde{M}/\Q},\widetilde{\omega})$. Furthermore, the number $\#\mathcal{M}(g_{\widetilde{M}/\Q},\widetilde{\omega})$ only depends on $\xi(g_{\widetilde{M}/\Q},\widetilde{\omega})$.
\end{lem}
\begin{proof}
This is essentially proved in \cite{LiLiuFamily}. The authors constructed a bundle $\mathcal{P}$ over $\Q$ whose fiber at $b$ consists of all pairs $(W,l)$, where $W$ is a maximal subspace of $H^{2}(M_{b};\mathbb{R})$ on which the intersection form is positive definite, and $l$ is a vector in $H^{2}(M_{b};\mathbb{R})$ that is not perpendicular to $W$. In our situation, since $\widetilde{M}$ is a homology product, the bundle $\mathcal{P}$ is trivial. And there is a bundle map from $\mathcal{P}$ to the trivial bundle $\Q\times V^{+}(M)$ defined by 
$$
(W,l)\mapsto \text{orthogonal projection of $l$ to $W$}
.$$
Since this bundle map induces a homotopy between the fibers, our lemma is essential the same as \cite[Proposition 3.6]{LiLiuFamily}, which states that the chamber structure is classified by homotopy classes of sections of $\mathcal{P}$.
\end{proof}
Proof of the following lemma is elementary so we omit.

\begin{lem}\label{lem: topology of positive cone}
Let $V$ be a real vector space with a symmetric non-degenerate bilinear form $\langle-,-\rangle$, with $k$ positive eigenvalues and $l$ negative eigenvalues. Assume $k\geq 1$. Then we have the the following conclusions:
\begin{enumerate}[(i)]
    \item The space $V^{+}:=\{\alpha\in V\setminus \{0\}\mid \langle \alpha,\alpha\rangle>0 \}$ is homotopy equivalent to $S^{k-1}$.
    \item For any $x\in V\setminus \{0\}$ such that $\langle x,x\rangle \geq 0$, the space 
    $$
    V^{+}_{x}:=\{\alpha\in V\setminus \{0\}\mid \langle \alpha,\alpha\rangle>0, \langle \alpha,x\rangle>0 \}.
    $$
    $V^{+}_{x}$ is always contractible. 
\end{enumerate}
\end{lem}

By Lemma \ref{lem: topology of positive cone}, we have $V^{+}(M)\simeq S^{b^{+}_{2}(M)-1}$. 

\begin{defi}
For any $\widetilde{\mathfrak{s}}$ with $d(\widetilde{\mathfrak{s}},\widetilde{M})=0$ and any $\xi\in [\Q,V^{+}(M)]$, we define the family Seiberg-Witten invariant $$\FSW_{\xi}(\widetilde{M},\widetilde{\mathfrak{s}}):=\#\mathcal{M}(g_{\widetilde{M}/\Q},\widetilde{\omega})$$ for any regular pair $(g_{\widetilde{M}/\Q},\widetilde{\omega})$ with $\xi(g_{\widetilde{M}/\Q},\widetilde{\omega})=\xi$.
\end{defi}

A chamber is called a canonical chamber if it is represented by a constant map. When $b^{+}_{2}(M)\geq 2$, the space path $V^{+}(M)$ is connected, so there is only one canonical chamber, denoted by $\xi_{c}$. When $b^{+}_{2}(M)=1$, the space $V^+$ has two components. By fixing an orientation on $H^{2}_{+}(M;\mathbb{R})$, we can talk about the positive component and the negative component. The corresponding chambers are called the ``positive canonical chamber'' and the ``negative canonical chamber'', denoted by $\xi^{+}_{c}$ and $\xi_{c}^{-}$ respectively. We will also write the family Seiberg-Witten invariant for $\xi_{c},\ \xi^{+}_{c},\ \xi^{-}_{c}$ as $\FSW(-),\ \FSW_{+}(-)$ and $\FSW_{-}(-)$ respectively.

\begin{ex}\label{ex: winding chamber}
Let $\Q=S^{2}$. Then $[\Q,V^{+}(M)]$ contains a single element $\xi_{c}$ when $b^{+}_{2}(M)\neq 1, 3$. When $b^{+}_{2}(M)=1$, we have $[\Q,V^{+}(M)]=\{\xi^{+}_{c}, \xi^{-}_{c}\}$. When $b^{+}_{2}(M)=3$, there is a bijection $[\Q,V^{+}(M)]\cong \mathbb{Z}$ provided by the mapping degree. Let $\xi^{j}$ be chamber corresponding to $j\in \mathbb{Z}$. Then the general wall-crossing formula \cite[Theorem 4.10]{LiLiuFamily} specializes to  
\begin{equation}\label{eq: wall-crossing}
    \FSW_{\xi^{j}}(\widetilde{M},\widetilde{\mathfrak{s}})=\FSW_{\xi_{c}}(\widetilde{M},\widetilde{\mathfrak{s}})+j\cdot  C(M,\mathfrak{s}).
\end{equation}
Here $C(M,\mathfrak{s})$ is the so-called ``wall-crossing number'' and it only depends on the fiber $M$. (See \cite{LiLiuFamily}) for the precise formula.) We call $\xi^{1}$ the \emph{winding chamber}.
\end{ex}

\begin{lem}\label{lem: twisting by a line bundle doesn't change FSW}
Suppose $d(\widetilde{\mathfrak{s}},\widetilde{M})=0$. Then one has 
\begin{equation}\label{twisting equality}
\FSW_{\xi}(\widetilde{M}/\Q,\widetilde{\mathfrak{s}})=\FSW_{\xi}(\widetilde{M}/\Q,\widetilde{\mathfrak{s}}\otimes p^{*}(L))
\end{equation}
for any complex line bundle $L$ over $\Q$ and any  chamber $\xi\in [\Q,V^{+}(M)]$.
\end{lem}
\begin{proof}
We pick a pair $(g_{\widetilde{M}/\Q},\widetilde{\omega})$ which is regular with respect to both $\widetilde{\mathfrak{s}}$ and $\widetilde{\mathfrak{s}}\otimes p^{*}(L)$. Then the moduli spaces 
$\mathcal{M}(\widetilde{\mathfrak{s}},g_{\widetilde{M}/\Q},\widetilde{\omega})
$ and $\mathcal{M}(\widetilde{\mathfrak{s}}\otimes p^{*}(L),g_{\widetilde{M}/\Q},\widetilde{\omega})$ are both finite sets. Therefore, there are only finitely many $b_{1},\cdots ,b_{n}\in \Q$ that appears as the first component of some point in either moduli spaces. For each $b_i$, we pick a unit vector $v_{i}$ in the corresponding fiber $L_{b_{i}}$ of $L$. The isomorphism $\phi\mapsto \phi\otimes v_{i}$ between the spinor bundles $
S_{b_i}^{\pm}\rightarrow S_{b_i}^{\pm}\otimes L_{b_{i}}$ induces an isomorphism 
$$
\widetilde{\mathfrak{s}}|_{M_{b_{i}}}\cong (\widetilde{\mathfrak{s}}\otimes p^{*}(L))|_{M_{b_{i}}}\quad \text{for }i=1,\cdots n.
$$
These isomorphisms together induce a sign-preserving bijection between $\mathcal{M}(\widetilde{\mathfrak{s}},g_{\widetilde{M}/\Q},\widetilde{\omega})$ and 
$\mathcal{M}(\widetilde{\mathfrak{s}}\otimes p^{*}(L),g_{\widetilde{M}/\Q},\widetilde{\omega})$. This finishes the proof.
\end{proof}
\begin{rmk}
In \cite{LiLiuFamily}, the authors studied a more general version of the family Seiberg-Witten invariant when $d(\mathfrak{s}_{
\widetilde{M}},\widetilde{M})>0$, by paring the fundamental class $[\mathcal{M}(g_{\widetilde{M}/\Q},\widetilde{\omega})]$ with various cohomology classes of the parametrized configuration space. To adapt (\ref{twisting equality}) to this more general setting, one needs to introduce certain correction terms. We will not pursue this generality in this paper.
\end{rmk}
The next lemma addresses the behavior of the family Seiberg-Witten invariant under a pull-back of the base. 

\begin{lem}\label{lem: FSW pull-back} Let $\widetilde{M}/\Q$ be a homology product with fiber $M$. Consider a smooth map $f:\Q'\rightarrow \Q$ from some manifold $\Q'$ of the same dimension with $\Q$. Then the pull-back family $\widetilde{M}'/\Q'$ is also a homology product. Moreover, given any family $\spinc$ structure $\widetilde{\mathfrak{s}}$ and any chamber $\xi\in [\Q,V^{+}(M)]$, we let $\mathfrak{s}_{\widetilde{M}'}$ and  $\xi'\in [\Q',V^{+}(M)]$ be their pull-backs via $f$. Assume $d(\widetilde{\mathfrak{s}},\widetilde{M})=0$. Then we have 
$$
\FSW_{\xi'}(\widetilde{M}'/\Q',\mathfrak{s}_{\widetilde{M}'})=\operatorname{deg}(f)\cdot 
\FSW_{\xi}(\widetilde{M}/\Q,\widetilde{\mathfrak{s}}).
$$
\end{lem}
\begin{proof} 
We pick an generic, admissible pair $(g_{\widetilde{M}/\Q},\widetilde{\omega})$ in the chamber $\xi$. Let $(g_{\widetilde{M}'/\Q'},\widetilde{\omega}')$ be the pull-back of  $(g_{\widetilde{M}/\Q},\widetilde{\omega})$ via $f$. We use $(g_{\widetilde{M}/\Q},\widetilde{\omega})$ and $(g_{\widetilde{M}'/\Q'},\widetilde{\omega}')$ to define the family Seiberg-Witten equations on $\widetilde{M}/\Q$ and $\widetilde{M}'/\Q'$ respectively.  By transversality, there are only finitely many points $b_{1},\cdots,b_{n}\in \Q$ such that the Seiberg-Witten equations have a solution on the fiber $M_{b_{i}}$. By suitably homotope $f$, we may assume that $b_{i}$ is a regular values of $f$. For a point $b'\in f^{-1}(b_{i})$, we let $\# \mathcal{M}_{b'}\in \mathbb{Z}$ be the the contribution of solutions on $M_{b'}$ to $\FSW_{\xi'}(\widetilde{M}'/\Q',\mathfrak{s}_{\widetilde{M}'})$ and let  $\# \mathcal{M}_{b_{i}}$ be the contribution of solutions on $M_{b}$ to $\FSW_{\xi}(\widetilde{M}/\Q,\widetilde{\mathfrak{s}})$. By definition, we have 
$$
\# \mathcal{M}_{b'}=\operatorname{deg}(f,b') \cdot \# \mathcal{M}_{b_{i}} 
$$
where $\operatorname{deg}(f,b')$ denotes the local degree of $f$ at $b'$. Summing over these equalities over all points in $\bigcup f^{-1}(b_{i})$, we finish the proof.
\end{proof}

It turns out that one can define group homomorphisms from $\pi_{k}(\Diff(M))$ to $\mathbb{Z}$ using the family Seiberg-Witten invariants (see \cite{RubermanObstruction,Baraglia2019,BaragliaKonno}). 
We briefly recall this construction here. In Section \ref{section: linearly independence}, we will use them to detect the $\mathbb{Z}^{\infty}$ summand in Theorem \ref{thm: infinitely generated diff}.

For simplicity, we focus on the case  $b_{1}(M)=0$. Let $\gamma$ be an element in $\pi_{k}(\Diff(M))$ for some $k\geq 0$. When $k=0$, we assume that $\gamma$ induces the identity map on $H_{*}(M;\mathbb{Z})$. Let $\widetilde{\gamma}: S^{k}\times M\rightarrow S^{k}\times M$ be a diffeomorphism that represents $\gamma$. Then we use $\widetilde{\gamma}$ as the clutching function and form the smooth fiber bundle 
$$
\widetilde{M}:=(D^{k+1}\times M)\cup_{\widetilde{\gamma}} (D^{k+1}\times M) 
$$
over $\Q=S^{k+1}$. The following lemma can be proved by a straightforward application of the Serre spectral sequence. 

\begin{lem}\label{lem: extension of spinc b1=0}
Any $\spinc$ structure $\mathfrak{s}$ on the fiber $M$  can be extended to a family $\spinc$ structure $\widetilde{\mathfrak{s}}$ on $\widetilde{M}/\Q$. Such extension is unique when $k\neq 1$. When $k=1$, any two extensions are related by a twisting with a complex line bundle pulled back from $\Q$.  
\end{lem}

Let $\mathfrak{s}$ be a $\spinc$-structures on $M$ with $d(\mathfrak{s},M)=-k-1$. We assume $b^{+}(M)>1$ for now. Then the Seiberg-Witten invariant for $\gamma$ is defined as
$$
\SW(\gamma,\mathfrak{s}):=\FSW_{\xi_{c}}(\widetilde{M},\widetilde{\mathfrak{s}}).
$$
Here $\widetilde{\mathfrak{s}}$ is any extension of $\mathfrak{s}$. By Lemma \ref{lem: extension of spinc b1=0} and Lemma \ref{lem: twisting by a line bundle doesn't change FSW}, $\SW(\gamma,\mathfrak{s})$ is well-defined. 

The following proposition is proved in \cite[Theorem 2.6 and Theorem 2.9]{baraglia2021non}.
\begin{pro}\label{pro: finiteness b+>1} The map $\SW(-,\mathfrak{s}):\pi_{k}(\Diff(M))\rightarrow \mathbb{Z}$ is a group homomorphism. Furthermore, for any $\gamma$, the number $\SW(\gamma,\mathfrak{s})$ is nonzero for only finitely many $\mathfrak{s}$ with $d(\mathfrak{s},M)=-k-1$.
\end{pro}

Next we turn to the case $b^{+}(M)=1$. We fix a choice of $H\in H^{2}(M;\mathbb{R})$ with $H\cdot H>0$ and use it to fix an orientation on $H^{+}(M;\mathbb{R})$. Then we define the ``small-perturbation Seiberg-Witten invariant'' $$
\SW_{H}(\gamma,\mathfrak{s}):= \begin{cases} \FSW_{\xi^{+}_{c}}(\widetilde{M},\widetilde{\mathfrak{s}})\quad &\text{if }c_{1}(\mathfrak{s})\cdot H\geq 0\\ 
\FSW_{\xi^{-}_{c}}(\widetilde{M},\widetilde{\mathfrak{s}}) \quad &\text{if }c_{1}(\mathfrak{s})\cdot H< 0
\end{cases}. 
$$

 The following lemma helps us to better understand the dependence of  $\SW_{H}(\gamma,\mathfrak{s})$ on $H$. 

\begin{lem}\label{lem: H and H' different sign} Suppose $H,H'$ belongs to the same component of $V^{+}(M)$. Then for any fixed $n$, there exists only finitely many $\spinc$ structures $\mathfrak{s}$ with $d(\mathfrak{s},M)=n$ such that
$
c_{1}(\mathfrak{s})\cdot H\geq 0$ whereas $c_{1}(\mathfrak{s})\cdot H<0
$.
\end{lem}
\begin{proof}
By rescaling $H$ and $H'$, we may assume that $H\cdot H=H\cdot H'=1$. Let $m=b^{-}_{2}(M)$. Then we can choose a basis $\{e_{0},e_{1},\cdots ,e_{m}\}$ for $H^{2}(M;\mathbb{R})$ such that
$$
e_{i}\cdot e_{j}:= \begin{cases} 1\quad &\text{if }i=j=0\\ -1 \quad &\text{if }i=j>0\\
0 &\text{if }i\neq j
\end{cases}. 
$$
and that $H=e_{0}$, $H'=e_{0}-ae_{1}$ for some $a\in [0,1)$. Assume that $$c_{1}(\mathfrak{s})=\sum\limits^{m}_{i=0}x_{i}e_{i}\text{ for }x_{i}\in \mathbb{R}.$$ Then our assumptions can be written as
$$
x^2_{0}-x^{2}_{1}-\cdots -x^{2}_{m}=4n-m+9,\quad x_{0}\geq 0,\quad x_0-a x_{1}<0.
$$
These (in)equalities together imply a uniform bound on $|x_{i}|$. So $c_{1}(\mathfrak{s})$ can take only finitely many values. Since there can be only finitely many $\mathfrak{s}$ with identical $c_{1}(\mathfrak{s})$, the proof is finished.
\end{proof}

The following proposition is an analogue of Proposition \ref{pro: finiteness b+>1} in the case $b^{+}_{2}(M)=1$.
\begin{pro}\label{pro: finiteness b+=1} The map $\SW_{H}(-,\mathfrak{s}):\pi_{k}(\Diff(M))\rightarrow \mathbb{Z}$ is a group homomorphism. Furthermore, for any $\gamma$, the number $\SW_{H}(\gamma,\mathfrak{s})$ is nonzero for only finitely many $\mathfrak{s}$ with $d(\mathfrak{s},M)=-k-1$.
\end{pro}
\begin{proof}
Proof of the first assertion is exactly the same as Proposition \ref{pro: finiteness b+>1}. So we focus on the second assertion. We fix a family metric $g_{\widetilde{M}/\Q}$. By Lemma \ref{lem: H and H' different sign}, changing $H$ in the same component of $V^{+}(M)$ does not affect our conclusion. Therefore, we may assume that $H$ can be represented by self-dual harmonic 2-form on $(M_{b_{0}},g_{b_{0}})$. We set the perturbation 2-form $\widetilde{\omega}$ to be identically $0$. 
A compactness argument for Seiberg-Witten moduli space (see \cite[Theorem 5.2.4]{morgan1996seiberg}) implies the following claim: There exist only finitely many $\mathfrak{s}$ with $d(\mathfrak{s},M)=-k-1$ such that the moduli space 
$
\mathcal{M}(g_{\widetilde{M}/\Q},0,\widetilde{\mathfrak{s}})
$
is nonempty for some extension $\widetilde{\mathfrak{s}}$ of $\mathfrak{s}$.
In particular, this implies that for all but finitely many $\mathfrak{s}$ with $d(\mathfrak{s},M)=-k-1$, the  pair $(g_{\widetilde{M}/\Q},0)$ is admissible with respect to some extension $\widetilde{\mathfrak{s}}$ of $\mathfrak{s}$ (in particular, we have $c_{1}(\mathfrak{s})\cdot H\neq 0$) and family Seiberg-Witten invariant   
$
\FSW_{\xi}(\widetilde{M}/\Q,\widetilde{
\mathfrak{s}})
$ equals $0$ for the corresponding chamber $\xi$. Since $\xi=\xi^{+}_{c}$ when $c_{1}(\mathfrak{s})\cdot H>0$ and $\xi=\xi^{-}_{c}$ when $c_{1}(\mathfrak{s})\cdot H<0$, we see that  
$$
\SW_{H}(\gamma,\mathfrak{s})=\FSW_{\xi}(\widetilde{M}/\Q,\widetilde{
\mathfrak{s}})=0.
$$ This finishes the proof.
\end{proof}

\subsection{Seiberg-Witten invariants of symplectic 4-manifolds}
In this section, we recall Taubes' remakable theorems about the Seiberg-Witten invariants of symplectic 4-manifolds. We will also state some of its corollaries that we need.

\begin{thm}[Taubes \cite{Taubes94,Taubes1996}]\label{thm: taubes's vanishing} Let $(M,\omega)$ be a symplectic 4-manifold. Suppose $b^{+}_{2}(X)>1$. Then the following result holds.
\begin{enumerate}
    \item $\SW(M,\mathfrak{s}_{J})=\pm 1$. Here  $\mathfrak{s}_{J}$ is the canonical $\spinc$ structure.
    \item Suppose $\SW(M,\mathfrak{s}_{J}+e)\neq 0$ for some $e\in H^{2}(M;\mathbb{Z})$. Then $0\leq e \cdot  [\omega]$. The equality is achieved only when $e=0$. Furthermore, for any $\omega$-compatible almost complex structure $J$, the Poincar\'e dual of $e$ can be represented by a $J$-holomorphic curve.  \end{enumerate}
When $b^{+}_{2}(M)=1$ and we orient $H^{+}(M;\mathbb{R})$ using $[\omega]$, the results hold for $\SW_{-}(M)$.
\end{thm}

\begin{rmk}
To avoid discussing wall-crossing, the statements of the main theorems in \cite{Taubes1996,Taubes94} assume $b^{+}_{2}(M)>1$. But the whole proof actually works for  $\SW_{-}(M)$ when $b^{+}_{2}(M)=1$ (see \cite[Proof of Theorem 0.3]{Taubes1996}).
\end{rmk}

Part of Theorem \ref{thm: taubes's vanishing} can be adapted to the family setting in a straightforward way.

\begin{pro}\label{pro: family vanishing}
Let $\widetilde{M}/\Q$ be a smooth family that carries a fiberwise symplectic structure $\{\omega_{b}\}_{b\in \Q}$. (I.e. symplectic structures on the fibers that varies smoothly with $b$). Let $\xi_{\omega}$ be the ``symplectic chamber'' defined by the map  $\Q\rightarrow V^{+}(M)$ that sends $b$ to $-[\omega_{b}]\in V^{+}(M)$. Let $\widetilde{\mathfrak{s}}$ be family a $\spinc$ structure such that $d(\widetilde{\mathfrak{s}},\widetilde{M})=0$ and $\widetilde{\mathfrak{s}}|_{M}=\mathfrak{s}_{J}+e$ for some nonzero $e\in H^{2}(M;\mathbb{Z})$. Suppose $\FSW_{\xi_{\omega}}(\widetilde{M};\widetilde{\mathfrak{s}})\neq 0$. Then there exists some $b\in \Q$ such that $e\cdot [\omega_{b}]>0$.  
\end{pro}
\begin{proof}
We pick a family metric $g_{\mathcal{M}/\Q}$ such that $g_{b}$ is compatible with $\omega_{b}$ for any $b$. Then we consider the the perturbation
\begin{equation}\label{eq: taubes' perturbation}
\omega'_{b}=-\\i\cdot F^{+}_{A_{\omega_{b}}}-r\omega_{b},  \end{equation}
where $A_{\omega_{b}}$ is a canonical $\spinc$ connection determined by $\omega_{b}$, and $r$ is a large positive number. Suppose $e\cdot [\omega]\leq 0$ for any $b\in \Q$. By the argument in \cite{Taubes94}, when $r$ is large enough, the Seiberg-Witten equations on any fiber have no solutions. So the corresponding family Seiberg-Witten invariant is $0$. Note that when $r\gg 0$, the term $-r\omega_{b}$ dominates the other terms in (\ref{eq: taubes' perturbation}). So the chamber corresponding to $(g_{\mathcal{M}/\Q},\{\omega'_{b}\}_{b\in\Q})$ is the symplectic chamber $\xi_{\omega}$.
\end{proof}

Another important theorem of Taubes is the following result, which states that symplectic $4$-manifolds are all of ``simple-type''.
\begin{thm}[Taubes \cite{Taubes1996}]\label{thm: simple-type} Let $(M,\omega)$ be a symplectic $4$-manifold with $b^{+}_{2}>1$ and let $\mathfrak{s}$ be a $\spinc$-structure such that $d(\mathfrak{s},M)>0$. Then $\SW(M,\mathfrak{
s})=0$.
\end{thm}

Next, we state some results deduced from Theorem \ref{thm: taubes's vanishing}. 

\begin{thm}[Taubes \cite{Taubes1996}]\label{thm: Taubes2} Let $(M,\omega)$ be a symplectic 4-manifold. Suppose $e\in H_{2}(M;\mathbb{Z})$ is represented by a smoothly embedded 2-sphere with self-intersection $-1$. Then for any $\omega$-compatible almost complex structure $J$, there is an embedded $J$-holomorphic curve $S$ in $M$ with self-intersection $-1$. Furthermore, if we assume $b^{+}_{2}(M)>1$, then we may pick $S$ such that it represents $\pm e$. 
\end{thm}

In the case $b^{+}_{2}(M)=1$, Li-Liu proved the following related result. 
\begin{thm}[Li-Liu \cite{LiLiu95}]\label{thm: Li-Liu} Let $(M,\omega)$ be a symplectic 4-manifold with $b^{+}_{2}(M)=1$. Let $K$ be the canonical class and let $e\in H_{2}(M;\mathbb{Z})$ be a class with $e\cdot e=K\cdot e=-1$. Suppose $\operatorname{PD}(e)$ can be represented by a smoothly embedded 2-sphere with self-intersection $-1$. Then for any $\omega$-compatible almost complex structure $J$, there is a $J$-holomorphic curve $C$ that represents $\operatorname{PD}(e)$.   
\end{thm}

As for smoothly embedded sphere with self-intersection $-2$, we have the following result: 

\begin{lem}\label{pro: minimal implies adjunction}
Let $(M,\omega)$ be a symplectic 4-manifold that contains a smoothly embedded sphere $S$ with self-intersection $-2$. Assume $M$ is minimal. Then the pairing $\langle K,[S]\rangle$ must be $0$.\end{lem}
\begin{proof}
The proof is adapted from Theorem \ref{thm: Taubes2}. We use the \emph{reflection diffeomorphism} 
$
r: M\rightarrow M
$ constructed as follows. Let $N$ be a closed tubular neighborhood of $S$. Then we identify $N$ with the unit disk bundle of the sphere $S^{2}$. A reflection on $S^{2}$ induces an orientation preserving diffeomorphism $r':N\rightarrow N$. The restriction $r'|_{\partial N}$ is an orientation preserving isometry on $\partial N=\mathbb{RP}^{3}$. The group of such isometry is $SO(3)\times SO(3)$. So $r'|_{\partial N}$  is smoothly isotopic to the identity map. By modifying $r'$ near $\partial N$, we get a diffeomorphism $r''$ on $N$ that is identity near $\partial N$. By extending $\tau''$ with the identity map, we get the diffeomorphism $r:M\rightarrow M$ that is supported near $S$. The action $r^{*}$ on $H^{2}(M;\mathbb{Z})$ is just the reflection along $\operatorname{PD}[S]$. Therefore, we have 
$$
r^{*}(\mathfrak{s}_{J})=\mathfrak{s}_{J}-\langle K,[S]\rangle \cdot \operatorname{PD}([S])
$$
Now we assume $b^{+}_{2}(M)>1$. (The argument when $b^{+}_{2}(M)=1$ is the same if we use $\SW_{-}$ instead of $\SW$.) Suppose $K(S)\neq 0$. Then 
$$
\SW(M,r^{*}(\mathfrak{s}_{J}))=\SW(M,\mathfrak{s}_{J})=\pm 1.
$$
By Theorem \ref{thm: taubes's vanishing} (2), the class $-\langle K,[S]\rangle\cdot [S]$ is represented by a $J$-holomorphic curve $C$. However, this is impossible: The self-intersection of $C$ equals $-\langle K,[S]\rangle^{2}<0$. But by the regularity theorem of $J$-holomorphic curves (\cite[Proposition 7.1]{Taubes1996}), for a generic $J$, any $J$-holomorphic curve with negative self-intersection must contain a genus-$0$ component with self-intersection $-1$. This contradicts with our assumption that $M$ is minimal.
\end{proof}

\section{Recollection of the monopole Floer homology}\label{sec: definition of HM}

\subsection{Definition of monopole Floer homology}
In this section, we briefly recall the construction of the monopole Floer homology of 3-manifolds, following \cite{KronheimerMrowkaMonopole}. See also \cite{linmonopolelectures} for a self-contained survey. 

Let $Y$ be a compact, connected, closed and oriented 3-manifold equipped with a Riemannian metric $g_{Y}$. Given a $\spinc$ structure $\mathfrak{s}$ on $Y$, we can associate the spinor bundle $S_{Y}$, which a rank-2 Hermitian bundle over $Y$. We fix an integer $k\geq 3$ and let $\mathcal{C}(Y)$ be the $L^{2}_{k-1/2}$-completion of the affine space: 
$$
\{(B,\phi)\mid\text{$B$ is a smooth  $\spinc$-connection on $S_{Y}$ and $\phi$ is a smooth section of $S_{Y}$}\}
\}.
$$
We call $\mathcal{C}(Y)$ \emph{the configuration space downstairs}. There is a related space $\mathcal{C}^{\sigma}(Y)$ called \emph{the configuration space upstairs}, which is the  $L^{2}_{k-1/2}$-completion of the space: 
\begin{equation}\label{eq: 3d blown-up configuration space}
\begin{split} \{(B,s,\psi)\mid
\text{$B$ is a smooth  $\spinc$-connection on $S_{Y}$}, s\in \mathbb{R}^{\geq 0}\\ \text{ and $\phi$ is a smooth section of $S_{Y}$ with $|\phi|_{L^{2}}=1$.}  \}\end{split}
\end{equation}
There is a natural projection map $p: \mathcal{C}^{\sigma}(Y)\rightarrow \mathcal{C}(Y)$ defined by $p(B,s,\phi)=(B,s\phi)$. We call a point $(B,\phi)\in \mathcal{C}(Y)$ (resp. a point $(B,s,\phi)\in \mathcal{C}^{\sigma}(Y)$) reducible if $\phi=0$ (resp. $s=0$) and call it irreducible otherwise. We use $\mathcal{C}^{*}(Y)$ to denote the subspace of $\mathcal{C}(Y)$ consisting of irreducible points. 

The gauge group $\mathcal{G}$ is defined as the $L^{2}_{k+1/2}$-completion of the space of smooth maps from $Y$ to $S^{1}$. Then $\mathcal{G}$ acts on $\mathcal{C}(Y)$ as 
$$
g(B,\psi):=(B-g^{-1}dg,g\phi)
$$
and acts on  $\mathcal{C}^{\sigma}(Y)$. as
$$
g(B,s,\phi):=(B-g^{-1}dg,s,g\phi).
$$
We denote the quotient of $\mathcal{C}(Y),\ \mathcal{C}^{*}(Y)$ and $\mathcal{C}^{\sigma}(Y)$ by $\mathcal{G}$ respectively by $\mathcal{B}(Y),\ \mathcal{B}^{*}(Y)$ and $\mathcal{B}^{\sigma}(Y)$. Since the action of $\mathcal{G}$ on $\mathcal{C}^{*}(Y)$ and $\mathcal{C}^{\sigma}(Y)$ are free, the space $\mathcal{B}^{*}(Y)$ is a Hilbert manifold and $\mathcal{B}^{\sigma}(Y)$ is a Hilbert manifold with boundary. The map $p$ induces a homeomorphism between $\mathcal{B}^{*}(Y)$ and the interior of $\mathcal{B}^{\sigma}(Y)$. 

There is a Chern-Simons-Dirac functional $\mathcal{L}:C(Y)\rightarrow \mathbb{R}$ (see \cite[Definition 4.1.1]{KronheimerMrowkaMonopole} for the precise formula). Roughly speaking, the monopole Floer homology is modeled on an infinite dimensional Morse homology for $\mathcal{L}$. To achieve transversality, one needs to add a generic, gauge invariant perturbation function $f: \mathcal{C}(Y)\rightarrow \mathbb{R}$. Let $\mathfrak{q}$ be the formal gradient of $f$. Then $\mathfrak{q}$ has two components $(\mathfrak{q}^{1},\mathfrak{q}^{2})$, where $\mathfrak{q}^{1}$ is the connection component and $\mathfrak{q}^{2}$ is the spinor component. 

The gradient of the perturbed functional $$\mathcal{L}_{\mathfrak{q}}:=\mathcal{L}+f:\mathcal{C}(Y)\rightarrow \mathbb{R}$$ is invariant under $\mathcal{G}$ and hence induces a vector field $v^{\sigma}_{\mathfrak{q}}$ on $\mathcal{B}^{\sigma}(Y)$. (More precisely,  $v_{\mathfrak{q}}^{\sigma}$ is a section of the $L^{2}_{k-3/2}$-completion of the tangent bundle $T\mathcal{B}^{\sigma}(Y)$.)

The vector field $v_{\mathfrak{q}}^{\sigma}$ has finitely many critical points. 
Furthermore, there are two types of reducible critical points: the boundary stable ones and the boundary unstable ones. To recall their difference, we note that a reducible point $(B,0)\in C(Y)$ is a critical point of $\mathcal{L}_{\mathfrak{q}}$ if and only if $B$ satisfies the perturbed curvature equation
\begin{equation}\label{eq: perturbed curvature equation}
*F_{B^{t}}=\mathfrak{q}^{1}(B,0),
\end{equation}
where $B^{t}$ denotes the induced connection on determinant bundle $\bigwedge ^{2}S_{Y}$. We consider the perturbed Dirac operator 
\begin{equation}\label{eq: perturbed Dirac operator 3d}
\slashed{D}_{B,\mathfrak{q}}:L^{2}_{k}(S_{Y})\rightarrow L^2_{k-1}(S_{Y}),
\end{equation}
defined by 
$$
\slashed{D}_{B,\mathfrak{q}}(\psi)=\slashed{D}_{B}\phi+\mathcal{D}_{(B,0)}\mathfrak{q}^{1}(0,\psi),
$$
where the last term denotes the derivative of $\mathfrak{q}^{1}$ at $(B,0)$ in the direction $(0,\psi)$. This is a self-adjoint elliptic operator. Part of the transversality assumption requires that the eigenvalues of $\slashed{D}_{B,\mathfrak{q}}$ are all simple and nonzero. We arrange these eigenvalues as 
$$
\cdots \lambda_{-2}<\lambda_{-1}<0<\lambda_{0}<\lambda_{1}<\cdots .
$$
Then each $\lambda_{i}$ corresponds to a reducible critical point $[(B,0,\psi_{i})]$ upstairs, where $\psi_{i}$ is a unit-length eigenvector for $\lambda_{i}$. We call the critical point $[(B,0,\psi_{i})]$ boundary stable (resp. boundary unstable) if $\lambda_{i}>0$ (resp. $\lambda_{i}<0$).

Let $\mathfrak{C}(Y)$ be the set of critical points in $\mathcal{B}^{\sigma}(Y)$. Then we have a decomposition
$$
\mathfrak{C}(Y)=\mathfrak{C}^{o}(Y)\cup \mathfrak{C}^{s}(Y)\cup \mathfrak{C}^s(Y)
$$
according to the three types: irreducible, boundary stable and boundary unstable. We will suppress $Y$ from our notations when it is clear from the context. 

Given any pairs of critical points $[\mathfrak{a}],[\mathfrak{b}]\in \mathfrak{C}$, we consider the moduli space $\breve{\mathcal{M}}([\mathfrak{a}],[\mathfrak{b}])$ of \emph{unparametrized}, nonconstant flow lines of $v^{\sigma}_{\mathfrak{q}}$ going from $[\mathfrak{a}]$ to $[\mathfrak{b}]$. This moduli space is an oriented manifold. Note that different components of $\breve{\mathcal{M}}([\mathfrak{a}],[\mathfrak{b}])$ may have different dimensions. We use $\# \breve{\mathcal{M}}([\mathfrak{a}],[\mathfrak{b}])$ to denote the signed count of points in the 0-dimensional components. Similarly, for each pair of \emph{reducible} critical points $[\mathfrak{a}],[\mathfrak{b}]$, we use $\breve{\mathcal{M}}^{\text{red}}([\mathfrak{a}],[\mathfrak{b}])$ to denote the moduli space of unparametrized, nonconstant \emph{reducible} flow lines going from $[\mathfrak{a}]$ to $[\mathfrak{b}]$ and use  $\# \breve{\mathcal{M}}^{\text{red}}([\mathfrak{a}],[\mathfrak{b}])$ to denote the signed count of points in the 0-dimensional components. The compactness theorem for the Seiberg-Witten trajectories (\cite[Proposition 16.1.4]{KronheimerMrowkaMonopole}) ensures that $\# \breve{\mathcal{M}}([\mathfrak{a}],[\mathfrak{b}])$ and $\# \breve{\mathcal{M}}^{\text{red}}([\mathfrak{a}],[\mathfrak{b}])$ are always finite numbers.

Let $C^{o}$ (respectively $C^{s}$ and $C^{u}$) be the vector space over $\mathbb{Q}$ whose basis vectors $e_{[\mathfrak{a}]}$ are indexed by elements $[\mathfrak{a}]\in \mathfrak{C}^{o}$ (respectively $\mathfrak{C}^{o}$ and $\mathfrak{C}^{u}$). We define a linear map 
$
\partial^{o}_{s}:C^{o}\rightarrow C^{s}
$
by the formula 
\begin{equation}\label{eq: partial}
\partial^{o}_{s} e_{[\mathfrak{a}]}:=\sum\limits_{[\mathfrak{b}]\in \mathfrak{C}^{s}} \# \breve{\mathcal{M}}([\mathfrak{a}],[\mathfrak{b}]) \cdot e_{[\mathfrak{b}]}.
\end{equation}
The linear maps $\partial^{o}_{o}:C^{o}\rightarrow C^{o},\  \partial^{u}_{s}:C^{u}\rightarrow C^{s}$ and $\partial^{u}_{o}:C^{u}\rightarrow C^{o}$ are defined similarly. Moreover, by counting only reducible flow lines, one can define another four linear maps $$\widebar{\partial}^{s}_{s}:C^{s}\rightarrow C^{s},\  \widebar{\partial}^{s}_{u}:C^{s}\rightarrow C^{u},\ \widebar{\partial}^{u}_{s}:C^{s}\rightarrow C^{s},\quad \widebar{\partial}^{u}_{u}:C^{u}\rightarrow C^{u}$$
by the formula 
\begin{equation}\label{eq: partialbar}
\widebar{\partial}^{*}_{*}[\mathfrak{a}]:=\sum\limits_{[\mathfrak{b}]\in \mathfrak{C}^{*}} \# \breve{\mathcal{M}}^{\text{red}}([\mathfrak{a}],[\mathfrak{b}]) \cdot e_{[\mathfrak{b}]}.
\end{equation}
Consider the vector spaces 
$$
\widebar{C}=C^{s}\oplus C^{u},\quad \widecheck{C}=C^{o}\oplus C^{s},\quad \widehat{C}=C^{o}\oplus C^{u}.
$$
The three flavors of monopole Floer homology $\HMbar_{*}(Y,\mathfrak{s}), \widecheck{HM}_{*}(Y,\mathfrak{s}),\ \widehat{HM}_{*}(Y,\mathfrak{s})$ are defined as homology of the chain complexes $(\widebar{C},\widebar{\partial}),\  (\widehat{C},\widehat{\partial}),\ (\widecheck{C},\widecheck{\partial})$ with differentials 
\begin{equation}\label{eq: differential}
\widebar{\partial}=\left(\begin{array} {cc}
 \widebar{\partial}^{s}_{s}  & \widebar{\partial}^{u}_{s}  \\
 \widebar{\partial}^{s}_{u}  & \widebar{\partial}^{u}_{u}
\end{array}\right), \quad\widehat{\partial}=\left(\begin{array} {cc}
 \partial^{o}_{o}  & \partial^{u}_{o}  \\
 -\widehat{\partial}^{s}_{u}\partial^{o}_{s}  & -\widehat{\partial}^{u}_{u}-\widehat{\partial}^{s}_{u}\partial^{u}_{s}
\end{array}\right),\quad \widecheck{\partial}=\left(\begin{array} {cc}
 \partial^{o}_{o}  & -\partial^{u}_{o}\widehat{\partial}^{s}_{u}  \\
 \partial^{o}_{s}  & \widehat{\partial}^{s}_{s}-\partial^{u}_{s}\widehat{\partial}^{s}_{u}
\end{array}\right).
\end{equation}

The corresponding monopole Floer cohomology
$\HMbar^{*}(Y,\mathfrak{s}), \widehat{HM}^{*}(Y,\mathfrak{s}),\ \widecheck{HM}^{*}(Y,\mathfrak{s})$ are defined as homology of the dual complexes $\text{Hom}(\widebar{C},\mathbb{Q}), \text{Hom}(\widecheck{C},\mathbb{Q})$ and $\text{Hom}(\widehat{C},\mathbb{Q}).$ It's shown in \cite[Corollary 23.1.6]{KronheimerMrowkaMonopole} that all these flavors of monopole Floer (co)homology are independent with the choice of auxiliary data (e.g. perturbation and Riemannian matric) and hence are invariants of the pair $(Y,\mathfrak{s})$. Many results in \cite{KronheimerMrowkaMonopole} holds for all three flavors. In this situation, we will simply use the notation $HM^{\circ}_{*}(Y,\mathfrak{s})$ for homology and $HM_{\circ}^{*}(Y,\mathfrak{s})$ for cohomology. We define 
$$
HM^{\circ}_*(Y):=\bigoplus\limits_{\mathfrak{s}\in \spinc(Y)} HM^{\circ}_*(Y,\mathfrak{s}),\quad HM_{\circ}^*(Y):=\bigoplus\limits_{\mathfrak{s}\in \spinc(Y)} HM_{\circ}^*(Y,\mathfrak{s}).
$$
We use 
\begin{equation}\label{eq: pairing}
    \langle-,-\rangle: HM^{\circ}_*(Y,\mathfrak{s})\otimes HM_{\circ}^*(Y,\mathfrak{s})\rightarrow \mathbb{Q}
\end{equation}
to denote the canonical pairing induced by the obvious pairing on the chain level. 
\subsection{Basic properties of the monopole Floer homology}
In this section, we recall some basic properties of monopole Floer homology which will be useful in later sections. 

First, there are $U$-actions on $HM^{\circ}_{*}(Y,\mathfrak{s})$ and  $HM_{\circ}^{*}(Y,\mathfrak{s})$, making them into modules over the polynomial ring $\mathbb{Q}[U]$. The actions are defined as follows: Consider the canonical $S^{1}$-bundle $P$ over $\mathcal{B}^{\sigma}(Y)$, obtained as the quotient of $\mathcal{C}^{\sigma}(Y)$ by the based gauge group (the group of gauge transformations that equal $1$ at the base point). Let $\mathcal{M}([\mathfrak{a}],[\mathfrak{b}])$  (resp. $\mathcal{M}^{\text{red}}([\mathfrak{a}],[\mathfrak{b}])$) be the moduli spaces of \emph{parametrized} nonconstant flow lines (resp. reducible flow lines). Evaluating these flow lines at $t=0$ defines the maps  $$r:\mathcal{M}([\mathfrak{a}],[\mathfrak{b}])\rightarrow \mathcal{B}^{\sigma}(Y)\text{ and } r':\mathcal{M}^{\text{red}}([\mathfrak{a}],[\mathfrak{b}])\rightarrow \mathcal{B}^{\sigma}(Y).$$
We take a generic section of the line bundle associated to $r^{*}(P)$ (respectively $r'^{*}(P)$) and use it to cut the 2-dimensional components of the moduli space $\mathcal{M}([\mathfrak{a}],[\mathfrak{b}])$ (respectively $\mathcal{M}^{\text{red}}([\mathfrak{a}],[\mathfrak{b}])$) into a 0-dimensional manifold.  By counting points in these cut moduli spaces, one can define linear maps $l^{*}_{*},\ \overline{l}^{*}_{*}$ similar to $\partial^{*}_{*}$ and $\overline{\partial}^{*}_{*}$ (see (\ref{eq: partial}) and (\ref{eq: partialbar})). Then we define the chain maps 
$$
\widebar{l}:\widebar{C}\rightarrow \widebar{C},\quad \widecheck{l}:\widecheck{C}\rightarrow \widecheck{C},\quad \widehat{l}:\widehat{C}\rightarrow \widehat{C}
$$
by replacing $\partial^{*}_{*}$ and $\widehat{\partial}^{*}_{*}$ in (\ref{eq: differential}) with  $l^{*}_{*}$ and $\widehat{l}^{*}_{*}$ respectively. The induced map on (co)homology defines the $U$-action.

Second, the monopole (co)homology has two gradings: the rational grading and the mod-2 grading. Here we focus on the rational grading, which is only defined when  $c_{1}(\mathfrak{s})$ is torsion. There is a grading function
$$
\gr^{\mathbb{Q}}:\mathfrak{C}\rightarrow \mathbb{Q}.
$$
(see \cite[Definition 28.3.1]{KronheimerMrowkaMonopole} for definition) which induce a rational grading on $\widecheck{HM}_{*}(Y,\mathfrak{s})$ and $\widehat{HM}_{*}(Y,\mathfrak{s})$. The rational grading on $\HMbar_{*}(Y,\mathfrak{s})$ is induced by a lightly modified grading function $\widebar{\gr}^{\mathbb{Q}}: \mathfrak{C}^{u}\cup \mathfrak{C}^{s}\rightarrow \mathbb{Q}$ defined as
$$
\widebar{\gr}^{\mathbb{Q}}(\mathfrak{a})=\begin{cases} \gr^{\mathbb{Q}}(\mathfrak{a}) &\mbox{if } \mathfrak{a}\in \mathfrak{C}^{s} \\ 
\gr^{\mathbb{Q}}(\mathfrak{a})-1 &\mbox{if } \mathfrak{a}\in \mathfrak{C}^{u} \end{cases}.
$$
When $c_{1}(\mathfrak{s})$ is torsion, we will use the notation $HM^{\circ}_{k}(Y,\mathfrak{s})$ to denote the component of $HM^{\circ}_{*}(Y,\mathfrak{s})$ with respect to this rational grading $k$. The $U$-action decrease the rational on $HM^{\circ}_{*}(Y,\mathfrak{s})$ by 2.

Next, there is an exact triangle relating the three flavours of monopole Floer homology:
\begin{equation}\label{exact triangle}
    \cdots \xrightarrow{i_{*,Y}} \widecheck{HM}_{*}(Y,\mathfrak{s})\xrightarrow{j_{*,Y}}\widehat{HM}_{*}(Y,\mathfrak{s})\xrightarrow{p_{*,Y}}\HMbar_{*}(Y,\mathfrak{s})\xrightarrow{i_{*,Y}}\widecheck{HM}_{*}(Y,\mathfrak{s})\xrightarrow{j_{*,Y}} \cdots.
\end{equation}
(See \cite[Proposition 22.2.1]{KronheimerMrowkaMonopole} for a construction.)
All the maps in this exact triangle preserve the $U$-action. When $c_{1}(\mathfrak{s})$ is torsion, the maps $i_{*,Y},\ j_{*,Y}$ preserve the rational grading and the map $p_{*,Y}$ decrease the rational grading by $1$. 

Lastly, we discuss how the theory can be simplified when $Y$ is a rational homology 3-sphere. In this case, the curvature equation $F_{B^{t}}=0$ has a unique solution $B_{0}$ up to gauge transformations. One can achieve transversality by introducing a small perturbation $\mathfrak{q}$ with the property that the perturbed curvature equation (\ref{eq: perturbed curvature equation}) still have $B_{0}$ as its unique solution. (See \cite[Section 2.2]{LinSplitting} for an explicit construction of such perturbations.) Therefore, the reducible critical points upstairs are exactly $\{[(B_{0},0,\psi_{i})]\}_{i\in \mathbb{Z}}$, where $\psi_{i}$ is a unit-length eigenvalue vector of the perturbed Dirac operator (\ref{eq: perturbed Dirac operator 3d}).

Moreover, the moduli space of reducible trajectories can be explicitly described. In particular, $\breve{\mathcal{M}}^{\text{red}}([(B_{0},0,\psi_{i})],[(B_{0},0,\psi_{j})])$ never has a zero-dimensional component so the maps $\widebar{\partial}^{*}_{*}$ are vanishing. The only case that $\mathcal{M}^{\text{red}}([(B_{0},0,\psi_{i})],[(B_{0},0,\psi_{j})])$ has a 2-dimensional component is when $j=i+1$. We have
$$
\mathcal{M}^{\text{red}}([(B_{0},0,\psi_{i})],[(B_{0},0,\psi_{i+1})])=\mathbb{CP}^{1}
$$
and $r'^{*}(P)$ is the tautological bundle over it. Therefore, 
 the $U$-action on $\HMbar_{*}(Y,\mathfrak{s})$ sends $e_{[(B_{0},0,\psi_{i})]}$ to $e_{[(B_{0},0,\psi_{i-1})]}$. Hence we have an isomorphism: $$\HMbar_{*}(Y,\mathfrak{s})\cong \mathbb{Q}[U,U^{-1}].$$  The image of $p_{*,Y}$ is isomorphic to $\mathbb{Q}[U]$. The Fr\o yshov invariant $h(Y,\mathfrak{s})\in \mathbb{Q}$ is defined by the formula
 $$
 -2h(Y,\mathfrak{s})-2=\min\{\widebar{\gr}^{\mathbb{Q}}(\mathfrak{a})\mid \mathfrak{a}\in \text{Im}(p_{*})\}
 $$
 Namely, $-2h(Y,\mathfrak{s})-2$ equals the rational grading of the generator of the submodule $\text{Im}(p_{*})\subset \HMbar_{*}(Y,\mathfrak{s})$. 
 
 Recall that a rational homology 3-sphere $Y$ is called an \emph{L-space} (over $\mathbb{Q}$) if $p_{*,Y}$ is injective for any $\spinc$-structure $\mathfrak{s}$. In this case, $\widehat{HM}_{*}(Y,\mathfrak{s})$ is a free $\mathbb{Q}[U]$-module of rank 1 and the generator has the rational grading $-2h(Y,\mathfrak{s})-1$. 
 
 We end this section by introducing some notations.
 \begin{notation}\label{notation: generators}
 For any $k\in -2h(Y,\mathfrak{s})+2\mathbb{Z}$, we use $\widebar{e}_{k}(Y)$ to denote the unique element in $\HMbar_{k}(Y,\mathfrak{s})$ represented by a reducible critical point  $(B_{0},0,\psi_{*})$. When $k\geq -2h(Y,\mathfrak{s})$, the element $i_{*,Y}(\widebar{e}_{k})$ is nonzero in $\widecheck{HM}_{k}(Y,\mathfrak{s})$ and we denote it as $[\widecheck{e}_{k}(Y)]$. Furthermore, suppose $Y$ is an L-space. Then for any 
 $
 k\in -2h(Y,\mathfrak{s})-1+2\mathbb{Z}^{\geq 0},
 $ there is a unique element in $\widehat{HM}_{k}(Y,\mathfrak{s})$ that is sent to $\widebar{e}_{k-1}(Y)$ under $p_{*,Y}$. We use $\widehat{e}_{k}(Y)$ to denote this element. Note that unlike $\widecheck{e}_{k}(Y)$ and $\widebar{e}_{k}(Y)$, the class $\widehat{e}_{k}(Y)$ is not well-defined for a general rational homology 3-sphere because $p_{*,Y}$ is not necessarily injective. 
 \end{notation}
 
\section{Maps on induced by family cobordisms}
In this section, we discuss maps between monopole Floer homology induced by families of cobordisms. These maps were defined by Kronheimer-Mrowka \cite{KronheimerMrowkaMonopole} for a single cobordism. Although most of our constructions are adapted from there, we still provide some details for the sake of concreteness and completeness.
\subsection{Definition of maps induced by family cobordisms}
We start by clarifying the meaning of family  cobordisms in our setting.

\begin{defi}
Let $Y_{0},Y_{1}$ be two nonempty, connected oriented 3-manifolds and let $\Q$ be a closed, connected, oriented smooth manifold. A family of cobordisms $\widetilde{W}/\Q$ from $Y_{0}$ to $Y_{1}$ consists of the following data:
\begin{itemize}
    \item A smooth fiber bundle $p: \widetilde{W}\rightarrow \Q$ whose total space $\widetilde{W}$ is an oriented manifold and whose fiber $W$ is a smooth, connected 4-manifold with boundary. 
    \item An orientation preserving diffeomorphism 
    \begin{equation}\label{eq: boundary parameterization}
   f:\partial \widetilde{W}\xrightarrow{\cong}     (-Y_{0}\bigsqcup Y_{1})\times \Q
    \end{equation}
compatible with $p$, called the \emph{boundary parameterization} of the family. Here $-Y_{0}$ denotes the orientation reversal of $Y_{0}$. (In particular, we only consider bundles that are trivial when restricted to the boundary.)
\end{itemize}
\end{defi}
For family cobordisms, we will use the same set of notations as what we used for families of closed manifolds (see Section \ref{section: family SW}). In particular, we use $W_{b}$ to denote the fiber over $b\in \Q$. And  we say the family cobordism is a homology product if the monodromy action of $\pi_{1}(\Q)$ on $H_{*}(W;\mathbb{Z})$ is trivial. 
We equip the family with a smooth family of metrics $\{g_{W_{b}}\}$. We require that when restricted to a collar neighborhood $N_{i}$ of the boundary component $Y_{i}$, the metric $g_{W_{b}}$ is isomorphic to the product $[0,1)\times Y_{i}$  for a fixed metric $g_{Y_{i}}$ (independent with $b$). 
\begin{defi}
Let $T(\widetilde{W}/\Q)$ be the vertial tangent bundle and let $\text{Fr}$ be its frame bundle. A family $\spinc$ structure $\mathfrak{s}_{\widetilde{W}}$ on $\widetilde{W}/\Q$ is a lift of $\text{Fr}$ to a principal $\textnormal{Spin}^{c}(4)$-bundle $P_{\mathfrak{s}_{\widetilde{W}}}$. We say $\mathfrak{s}_{\widetilde{W}}$ and $\mathfrak{s}'_{\widetilde{W}}$ are isomorphic if there is a bundle isomorphism $P_{\mathfrak{s}_{\widetilde{W}}}\cong P_{\mathfrak{s}'_{\widetilde{W}}}$ that covers the identity map on $\text{Fr}$. 
\end{defi}
\begin{defi}\label{defi: spin-c cobobordism}
We say $\mathfrak{s}_{\widetilde{W}}$ is a product on the boundary if under the boundary parameterization (\ref{eq: boundary parameterization}), the restriction of  $\mathfrak{s}_{\widetilde{W}}$ to $\partial \widetilde{W}$ is isomorphic to the pull-back of a $\spinc$-structure $\mathfrak{s}_{Y_{i}}$ on $Y_{i}$ for $i=0,1$. In this case, we call $(\widetilde{W}/\Q,\mathfrak{s}_{\widetilde{W}})$ a family $\spinc$-cobordism from $(Y_{0},\mathfrak{s}_{Y_{0}})$ to $(Y_{1},\mathfrak{s}_{Y_{1}})$. 
\end{defi}

\begin{rmk}
In Definition \ref{defi: spin-c cobobordism}, we do not specify an identification between the spinor bundles $S^{\pm }_{\widetilde{W}}|_{\partial \widetilde{W}}$ and $\Q\times S^{\pm }_{Y_{i}}$. We will choose any such identification when we later define the induced map on monopole Floer homology. This choice does not affect the induced map.
\end{rmk}

By attaching cylindrical ends $(-\infty,0]\times Y_{0}$ and $[0,\infty)\times Y_{1}$ fiberwisely (via the boundary parameterization $f$), we obtain a family of noncompact cylindrical end manifolds, denoted by $\widetilde{W}^{*}/\Q$. The induced map  $$HM^{\circ}_{*}(\widetilde{W},\mathfrak{s}_{\widetilde{W}}):HM^{\circ}_{*}(Y_{0},\mathfrak{s}_{Y_{0}})
\rightarrow HM^{\circ}_{*}(Y_{1},,\mathfrak{s}_{Y_{1}})$$ is defined by studying the family Seiberg-Witten equations on $\widetilde{W}^{*}/\Q$. 

Since solutions of the Seiberg-Witten equations on $W_{b}^{*}$ usually have infinite $L^{2}$-norms, the blown-up configuration space is defined in a different way from the 3-dimensional case (\ref{eq: 3d blown-up configuration space}). We define 
$\mathcal{C}^{\sigma}(\widetilde{W}^{*},\mathfrak{s}_{\widetilde{W}})$ be the space of quadruples $(b,A,\Phi,\phi)$, where $b$ is a point in $\Q$; $A$ is a $L^{2}_{k,\text{loc}}$-connection on the spinor bundle $S^{+}_{W^*_{b}}$; $\phi$ is an element in the topological quotient 
$$
\{\text{nonzero }L^{2}_{k,\text{loc}}\text{-sections of }S^{+}_{W^{*}_{b}}\}/\mathbb{R}_{>0}
$$
and $\Phi\in \mathbb{R}_{\geq 0}\cdot \phi$. 
We treat $\mathcal{C}^{\sigma}(\widetilde{W}^{*},\mathfrak{s}_{\widetilde{W}})$ as a bundle over $\Q$. The fiber over $b$ is the blown-up configuration space for $W_{b}$ and it is acted upon by the gauge group 
$$
\mathcal{G}(W_b):=L^{2}_{k+3/2,\text{loc}}(\operatorname{Map}(W_{b},S^{1})).
$$
We use $\mathcal{B}^{\sigma}(\widetilde{W}^{*},\mathfrak{s}_{\widetilde{W}})$ to denote the quotient space of $\mathcal{C}^{\sigma}(\widetilde{W}^{*},\mathfrak{s}_{\widetilde{W}})$ by this fiberwise action.

The blown-up Seiberg-Witten equations on the fiber $W^{*}_{b}$ are written as 
\begin{equation}\label{eq: blown-up SW}
\begin{split}
    \rho(F^{+}_{A^{t}})+(\Phi\Phi^{*})_{0}=0\\
    \slashed{D}^{+}_{A}\phi=0.
\end{split}
\end{equation}

Again, one has to introduce suitable perturbations to achieve transversality. 
\begin{cons}\label{construction: perturbations on cobordism}
We perturb (\ref{eq: blown-up SW}) in three steps:
\begin{enumerate}[(i)]
    \item We choose a smooth family $\tilde{\omega}=\{\omega_{b}\}_{b\in \Q}$ of self-dual two forms on the fibers $W_{b}$. We require that $\omega_{b}$ is supported in  $W_{b}\setminus (N_{0}\cup N_{1})$. We use these two forms to perturb the curvature equation into
    \begin{equation}\label{eq: perturbed curvature equation on cobordism}
        \rho(F^{+}_{A^{t}})+(\Phi\Phi^{*})_{0}=\rho(i\omega_{b}).
    \end{equation}
    \item To define the perturbation on the the cylindrical ends $(-\infty,0]\times Y_{0}$ and $[0,+\infty)\times Y_1$, we choose 3-dimensional perturbations $\mathfrak{q}_{i}$ on $Y_{i}$ ($i=0,1$) so that  moduli spaces of Seiberg-Witten flow lines on $Y_{i}$ are regular. We apply the corresponding 4-dimensional perturbations $\hat{\mathfrak{q}}_{0}$ and $\hat{\mathfrak{q}}_{1}$ on the  cylindrical ends. (Recall that on cylindrical ends, the 4-dimensional Seiberg-Witten equations can be rewritten as the equations for the Seiberg-Witten flow on $Y_{i}$. Therefore, $\mathfrak{q}_{i}$ induces a perturbation $\hat{\mathfrak{q}}_{i}$ for the 4-dimensional Seiberg-Witten equations.)
\item To define the perturbation on the collar neighborhood $N_{i}$ of $\partial W_{b}$, we choose a smooth family of 3-dimensional perturbations $\{\mathfrak{p}'_{i,b}\}_{b\in \Q}$ on $Y_{i}$. Then we take the following perturbation on $N_{i}=[0,1]\times Y_{i}$:
\begin{equation}\label{eq: perturbation near boundary}
\hat{\mathfrak{p}}_{i,b}:=\beta(t)\hat{\mathfrak{q}}_{i}+\beta_{0}(t)\hat{\mathfrak{p}}'_{i,b}.
\end{equation} Here $\beta(t)$ is a cutoff function that equals $1$ near $t=0$ and equals $0$ near $t=1$. And $\beta_{0}$ is a bump function supported in $(0,1)$. This is the third part of the perturbation. \end{enumerate}
 \end{cons}
\begin{rmk}\label{rmk: perturbation 2-forms} We briefly explain the purpose of these perturbations. The perturbations $\{\omega_{b}\}$ in step (i) are introduced to avoid reducible solutions when the fiber $W_{b}$ has positive $b^{+}_{2}$ (see Lemma \ref{lem: no reducible}). Absence of reducible solutions is essential in the definition of $\HMarrow_{*}(\widetilde{W},\mathfrak{s}_{\widetilde{W}})$ but not required in the definition of $\widehat{HM}_{*}(\widetilde{W},\mathfrak{s}_{\widetilde{W}})$, $ \widecheck{HM}_{*}(\widetilde{W},\mathfrak{s}_{\widetilde{W}})$ and $\HMbar_{*}(\widetilde{W},\mathfrak{s}_{\widetilde{W}})$. So when defining the later three maps, one can simply set $\omega_{b}=0$. The perturbation $\hat{\mathfrak{q}}_{i}$ in step (ii)  is introduced to achieve the transversality on $\mathbb{R}\times Y_{i}$, which is necessary to define the monopole Floer homology groups.  The perturbations $\{\hat{\mathfrak{p}}_{i,b}\}$ in step (iii) are introduced to ensure the transversality of the family Seiberg-Witten equations on $\widetilde{W}^{*}$. Note that the only perturbations we impose on the cylindrical ends are $\hat{\mathfrak{q}}_{i}$, which is fiber-independent.
\end{rmk}

To simplify our notations, we will use $\hat{\mathfrak{p}}_{b}$ denote the sum of the perturbations in (ii) and (iii). We use $\hat{\mathfrak{p}}^{0}_{b}$ and $\hat{\mathfrak{p}}^{1}_{b}$ to denote its connection component and the spinor component respectively. The  triply perturbed blown-up Seiberg-Witten equations on $W^{+}_{b}$ can be written as
\begin{equation}\label{eq: triply perturbed SW}
\mathfrak{F}^{\sigma}_{\hat{\mathfrak{p}}_{b},\omega_{b}}(A,\phi,\Phi)=0.
\end{equation}
We refer to \cite[Page 156]{KronheimerMrowkaMonopole} for an explicit expression of $\mathfrak{F}^{\sigma}_{\hat{q}_{b},\omega_{b}}$. We mention that when $\Phi=0$(i.e. on the reducibles), the equation (\ref{eq: triply perturbed SW}) can be written as 
\begin{equation}
 \begin{split}
&\rho(F^{+}_{A^{t}})=-2\hat{\mathfrak{p}}_{b}^{0}(A,0)\\
&\slashed{D}^{+}_{A,\hat{\mathfrak{p}}_{b}}(\phi)=0.
 \end{split} 
\end{equation}
Here $\slashed{D}^{+}_{A,\hat{\mathfrak{p}}_{b}}$ is the perturbed 4-dimensional Dirac operator, defined as 
\begin{equation}\label{eq: perturbed 4d Dirac equation}
\slashed{D}^{+}_{A,\hat{\mathfrak{p}}_{b}}(\phi):=\slashed{D}^{+}_{A}\phi+\mathcal{D}_{(A_{b},0)}\hat{\mathfrak{p}}^{1}_{b}(0,\phi).
\end{equation}

Given critical points $[\mathfrak{a}]\in \mathfrak{C}(Y_{0})$ and $[\mathfrak{b}]\in \mathfrak{C}(Y_{1})$, we use $\mathcal{B}^{\sigma}([\mathfrak{a}],\widetilde{W}^{*},\mathfrak{s}_{\widetilde{W}},[\mathfrak{b}])$ to denote the subspace of $\mathcal{B}^{\sigma}(W^{*},\mathfrak{s}_{\widetilde{W}})$ consisting of points that satisfy the asymptotic condition
\begin{equation}\label{eq: cylindrical boundary condition}
\begin{split}
\lim_{t\rightarrow -\infty} [(A|_{\{t\}\times Y_{0}},\frac{\phi|_{\{t\}\times Y_{0}}}{|\phi|_{\{t\}\times Y_{0}}|_{L^{2}}}, |\Phi|_{\{t\}\times Y_{0}}|_{L^{2}})]=[\mathfrak{a}]\\
\lim_{t\rightarrow +\infty} [(A|_{\{t\}\times Y_{1}},\frac{\phi|_{\{t\}\times Y_{1}}}{|\phi|_{\{t\}\times Y_{1}}|_{L^{2}}}, |\Phi|_{\{t\}\times Y_{1}}|_{L^{2}})]=[\mathfrak{b}]\end{split}.   
\end{equation}  
(In particular, this condition requires that $\phi|_{\{t\}\times Y_{i}}\neq 0$ when $|t|\gg 0$.) 

We consider the moduli space $$\mathcal{M}([\mathfrak{a}],\widetilde{W}^{*},\mathfrak{s}_{\widetilde{W}},[\mathfrak{b}]):=\{[(b,A,\phi,\Phi)]\in \mathcal{B}^{\sigma}([\mathfrak{a}],\widetilde{W}^{*},\mathfrak{s}_{\widetilde{W}},[\mathfrak{b}]) \text{ that satisfies }(\ref{eq: triply perturbed SW})\}. $$
We let $\mathcal{M}^{\text{red}}([\mathfrak{a}],\widetilde{W}^{*},\mathfrak{s}_{\widetilde{W}},[\mathfrak{b}])$ be the subspace of $\mathcal{M}([\mathfrak{a}],\widetilde{W}^{*},\mathfrak{s}_{\widetilde{W}},[\mathfrak{b}])$ that consists of all reducible solutions. 

According to the path component of $\mathcal{B}^{\sigma}([\mathfrak{a}],\widetilde{W}^{*},\mathfrak{s}_{\widetilde{W}},[\mathfrak{b}])$, we have the  canonical decompositions 
$$
\mathcal{M}([\mathfrak{a}],\widetilde{W}^{*},\mathfrak{s}_{\widetilde{W}},[\mathfrak{b}])=\bigsqcup_{z\in \pi_{0}(\mathcal{B}^{\sigma}([\mathfrak{a}],\widetilde{W}^{*},\mathfrak{s}_{\widetilde{W}},[\mathfrak{b}]))}\mathcal{M}_{z}([\mathfrak{a}],\widetilde{W}^{*},[\mathfrak{b}]),
$$
$$
\mathcal{M}^{\text{red}}([\mathfrak{a}],\widetilde{W}^{*},\mathfrak{s}_{\widetilde{W}},[\mathfrak{b}])=\bigsqcup_{z\in \pi_{0}(\mathcal{B}^{\sigma}([\mathfrak{a}],\widetilde{W}^{*},\mathfrak{s}_{\widetilde{W}},[\mathfrak{b}]))}\mathcal{M}^{\text{red}}_{z}([\mathfrak{a}],\widetilde{W}^{*},[\mathfrak{b}]).
$$
\begin{rmk}  For each $b\in \Q$, the inclusion map $i_{b}: W_{b}\hookrightarrow \widetilde{W}$ induces a map $$i^{*}_{b}: \pi_{0}(\mathcal{B}^{\sigma}([\mathfrak{a}],\widetilde{W}^{*},\mathfrak{s}_{\widetilde{W}},[\mathfrak{b}]))\rightarrow \pi_{0}(\mathcal{B}^{\sigma}([\mathfrak{a}],W_{b}^{*},\mathfrak{s}_{W},[\mathfrak{b}])).$$ This map is always a bijection. To see this, we take $[\gamma],[\gamma']\in \mathcal{B}^{\sigma}([\mathfrak{a}],W_{b}^{*},\mathfrak{s}_{W},[\mathfrak{b}])$. On the cylindrical end $(-\infty,0]\times Y_{0}$, $\gamma$ (resp. $\gamma'$) converges to a critical point $\mathfrak{a}_{1}$ (resp. $\mathfrak{a}_{2}$) in $\mathcal{C}^{\sigma}(Y_0,\mathfrak{s}_{Y_{0}})$. Let $\gamma_{\mathfrak{a}}$ be a path from $\mathfrak{a}_{1}$ to $\mathfrak{a}_{2}$ and let $\gamma_{\mathfrak{b}}$ be a path in $\mathcal{C}^{\sigma}(Y_1,\mathfrak{s}_{Y_{1}})$ defined similarly. Then $\gamma_{\mathfrak{a}}\cdot \gamma'\cdot \overline{\gamma}_{\mathfrak{b}}$ have the same end points as $\gamma$. The linear homotopy between $\gamma_{\mathfrak{a}}\cdot \gamma'\cdot \overline{\gamma}_{\mathfrak{b}}$ and $\gamma$ implies that $[\gamma_{\mathfrak{a}}\cdot \gamma'\cdot \overline{\gamma}_{\mathfrak{b}}]$ and $[\gamma']$ are in the same component of $\mathcal{B}^{\sigma}([\mathfrak{a}],W_{b}^{*},\mathfrak{s}_{W},[\mathfrak{b}])$. Note that $\gamma_{\mathfrak{a}}$ and $\gamma_{\mathfrak{b}}$ represent elements in $\pi_{1}(\mathcal{B}^{\sigma}(Y_{0},\mathfrak{s}_{Y_{0}}),[\mathfrak{a}])$ and $\pi_{1}(\mathcal{B}^{\sigma}(Y_{1},\mathfrak{s}_{Y_{1}}),[\mathfrak{b}])$ respectively. From this construction, we get a transitive action of $$\pi_{1}(\mathcal{B}^{\sigma}(Y_{0},\mathfrak{s}_{Y_{0}}),[\mathfrak{a}])\times \pi_{1}(\mathcal{B}^{\sigma}(Y_{1},\mathfrak{s}_{Y_{1}}),[\mathfrak{b}])$$ on $\pi_{0}(\mathcal{B}^{\sigma}([\mathfrak{a}],W_{b}^{*},\mathfrak{s}_{W},[\mathfrak{b}]))$. Now we treat  $\mathcal{B}^{\sigma}([\mathfrak{a}],\widetilde{W}^{*},\mathfrak{s}_{\widetilde{W}},[\mathfrak{b}])$ as a fiber bundle over $\Q$ with fibers $\{\mathcal{B}^{\sigma}([\mathfrak{a}],W_{b}^{*},\mathfrak{s}_{W},[\mathfrak{b}])\}$. Since $\widetilde{W}$ is a product near the boundary, $\pi_{1}(\Q)$ acts trivially on  $\pi_{1}(\mathcal{B}^{\sigma}(Y_{0},\mathfrak{s}_{Y_{0}}),[\mathfrak{a}])\times \pi_{1}(\mathcal{B}^{\sigma}(Y_{1},\mathfrak{s}_{Y_{1}}),[\mathfrak{b}])$ and hence also trivially on the set $\pi_{0}(\mathcal{B}^{\sigma}([\mathfrak{a}],W_{b}^{*},\mathfrak{s}_{W},[\mathfrak{b}]))$. Therefore, the map $i^{*}_{b}$ is an bijection. In Section \ref{subsection: composition}, this observation will allow us to compose homotopy classes $z_{1}, z_{2}$ of configurations for two family cobordisms $\widetilde{W}_{01}$ and $\widetilde{W}_{12}$ and obtain a homotopy class $z_{1}\circ z_{2}$ for the composed cobordism $\widetilde{W}_{01}\cup \widetilde{W}_{12}$, without worrying the choice of a fiber. 
\end{rmk}

With a choice of generic perturbations, these moduli spaces are smooth manifolds (possibly with boundaries). Moreover, they are canonically oriented after we choose an orientation of the space 
$$
H^{1}(Y_{0};\mathbb{R})\oplus H^{1}(W_{b};\mathbb{R})\oplus H^{2}_{+}(W_{b};\mathbb{R})
$$
for some $b$ (and hence all $b$). By counting points in the zero-dimensional components of $\mathcal{M}([\mathfrak{a}],\widetilde{W}^{*},\mathfrak{s}_{\widetilde{W}},[\mathfrak{b}])$ and $\mathcal{M}^{\text{red}}([\mathfrak{a}],\widetilde{W}^{*},\mathfrak{s}_{\widetilde{W}},[\mathfrak{b}])$, one can define operators
\begin{equation}\label{eq: m operator}
m ^{o}_{o}:C^{o}(Y_{0})\rightarrow C^{o}(Y_{1}),\quad m ^{o}_{s}:C^{o}(Y_{0})\rightarrow C^{s}(Y_{1}),\quad m ^{u}_{o}:C^{u}(Y_{0})\rightarrow C^{o}(Y_{1})
\end{equation}
and 
\begin{equation}\label{eq: m-bar operator}
\begin{split}
    \widebar{m}^{s}_{s}:C^{s}(Y_{0})\rightarrow C^{s}(Y_{1}),\quad \widebar{m}^{u}_{u}:C^{u}(Y_{0})\rightarrow C^{u}(Y_{1})\\
    \widebar{m}^{u}_{s}:C^{u}(Y_{0})\rightarrow C^{s}(Y_{1}),\quad \widebar{m}^{s}_{u}:C^{s}(Y_{0})\rightarrow C^{u}(Y_{1})
\end{split}    
\end{equation}
respectively. As in \cite{KronheimerMrowkaMonopole}, we pack these maps as follows:
\begin{equation}\label{eq: mbar}
\widebar{m}=\left[\begin{array} {cc}
 \widebar{m}^{s}_{s}  & \widebar{m}^{u}_{s}  \\
 \widebar{m}^{s}_{u}  & \widebar{m}^{u}_{u}
\end{array}\right]:\quad \widebar{C}(Y_{0})\rightarrow \widebar{C}(Y_{1}),
\end{equation}

\begin{equation}\label{eq: mhat}
\widehat{m}=\left[\begin{array} {cc}
 m^{o}_{o}  & m^{u}_{o}  \\
 \widebar{m}^{s}_{u}\partial^{o}_{s}(Y_{0})-\widebar{\partial}^{s}_{u}(Y_{1})m^{o}_{s}  & \widebar{m}^{u}_{u}+\widebar{m}^{s}_{u}\partial^{u}_{s}(Y_{0})-\widebar{\partial}^{s}_{u}(Y_{1})m^{u}_{s}
\end{array}\right]:\quad \widehat{C}(Y_{0})\rightarrow \widehat{C}(Y_{1}),
\end{equation}

\begin{equation}\label{eq: mcheck}
\widecheck{m}=\left[\begin{array} {cc}
 m^{o}_{o}  & -m^{u}_{o}\widebar{\partial}^{s}_{u}(Y_{0})-\partial^{u}_{o}(Y_{1})\widebar{m}^{s}_{u}  \\
 m^{o}_{s}  & \widebar{m}^{s}_{s}-m^{u}_{s}\widebar{\partial}^{s}_{u}(Y_{0})-\partial^{u}_{s}(Y_{1})\widebar{m}^{s}_{u}
\end{array}\right]:\quad \widecheck{C}(Y_{0})\rightarrow \widecheck{C}(Y_{1}).
\end{equation}

\begin{pro}\label{pro: cobordism map well-defined}
The maps $\widebar{m},\widecheck{m},\widehat{m}$ are chain maps. The induced maps on homology
\begin{equation}\label{eq: maps induced by cobordisms}
\begin{split}
\widehat{HM}_{*}(\widetilde{W},\mathfrak{s}_{\widetilde{W}}):\widehat{HM}_{*}(Y_{0},\mathfrak{s}_{Y_0})\rightarrow \widehat{HM}_{*}(Y_{1},\mathfrak{s}_{Y_1}) \\
\widecheck{HM}_{*}(\widetilde{W},\mathfrak{s}_{\widetilde{W}}):\widecheck{HM}_{*}(Y_{0},\mathfrak{s}_{Y_0})\rightarrow \widecheck{HM}_{*}(Y_{1},\mathfrak{s}_{Y_1}) \\
\HMbar_{*}(\widetilde{W},\mathfrak{s}_{\widetilde{W}}):\HMbar_{*}(Y_{0},\mathfrak{s}_{Y_0})\rightarrow \HMbar_{*}(Y_{1},\mathfrak{s}_{Y_1}) 
\end{split}
\end{equation}
are smooth invariants of $(\widetilde{W},\mathfrak{s}_{\widetilde{W}})$. (Namely, they do not depend on auxiliary choices made in the construction such as metrics and perturbations.) They preserve the $U$-actions and commutes with the maps $i_{*},j_{*},p_{*}$ in the long exact sequence (\ref{exact triangle}).
\end{pro}
\begin{proof}When $\Q$ is a single point, all these claims are proved \cite{KronheimerMrowkaMonopole}. In particular, it is proved in  \cite[Proposition 25.3.4]{KronheimerMrowkaMonopole} that  $\widebar{m},\widecheck{m},\widehat{m}$ are chain maps. By \cite[Corollary 25.3.9]{KronheimerMrowkaMonopole}, the induced maps are independent with the metric $g_{\widetilde{W}/\Q}$ and the perturbations in Step (i) and (iii) of Construction \ref{construction: perturbations on cobordism}. Therefore, these maps could only possibly depend on the metric/perturbation $(g_{Y_{i}},\mathfrak{q}_{i})$ on the cylindrical ends. This potential dependence is ruled out as a consequence of the composition law of the cobordism-induced maps (\cite[Proposition 26.1.2]{KronheimerMrowkaMonopole}). Furthermore, by \cite[Proposition 25.5.1 and Corollary 25.5.2]{KronheimerMrowkaMonopole}, these maps preserve the $U$-action and they commute with  $i_{*},j_{*},p_{*}$ up to sign. 

For a general $\Q$, all the arguments can adapted without essential change. In particular, we will give a proof of the composition law for the induced map by family of cobordisms (see Theorem \ref{thm: composition law for HM}). 
\end{proof}

Next, we discuss a fourth map $$\HMarrow_{*}(\widetilde{W},\mathfrak{s}_{\widetilde{W}})_{\xi}:\widehat{HM}_{*}(Y_{0},\mathfrak{s}_{Y_0})\rightarrow \widecheck{HM}_{*}(Y_{1},\mathfrak{s}_{Y_1}).$$ Unlike its cousins in $(\ref{eq: maps induced by cobordisms})$, this map is only defined when $b^{+}_{2}(W)>0$. Furthermore, when $\mathfrak{s}_{Y_{0}}$ and $\mathfrak{s}_{Y_{1}}$ are both torsion and $b^{+}_{2}(W)<\text{dim}\Q+1$, this map depends on the choice of a chamber $\xi$ for $\widetilde{W}/\Q$. 

We start by assuming that that both $\mathfrak{s}_{Y_0}$ and $\mathfrak{s}_{Y_1}$ are torsion. For any $b\in \Q$,
we consider the space
\begin{equation}\label{equation: IW}
I(W_{b}):= \operatorname{Image}(H^{2}(W_{b},\partial W_{b};\mathbb{R})\xrightarrow{l} H^{2}(W_{b};\mathbb{R})),  \end{equation}
where $l$ is induced by the inclusion map. Then there is a non-degenerate pairing on $I(W_{b})$ defined as
\begin{equation}\label{eq: intersection pairing}
    \langle l(\alpha),l(\beta)\rangle :=(\alpha\cup \beta)[W_{b}].
\end{equation}
(Recall that $b^{+}_{2}(W_{b})$ is defined as the maximal dimension of a subspace of $I(W)$ on which this pairing is positive definite.) We can also identify $I(W_{b})$ with the space of square-integrable harmonic 2-forms $\kappa$ on $W_{b}^{*}$. For such $\kappa$, we use $[\kappa]\in I(W_{b})$ to denote the corresponding cohomology class. For any square-integrable 2-form  $\omega$ on $W_{b}^{*}$ (not necessarily harmonic), we use $P_{I}(\omega)$ to denote the unique element in $I(W_{b})$ such that
$$
\langle P_{I}(\omega),[\kappa]\rangle=\int_{W_{b}} \omega\wedge \kappa
$$
for any square-integrable harmonic 2-form $\kappa$.

The identification between $I(W_{b})$ and the space of harmonic forms gives us an involution $$*_{g_{b}}:I(W_{b})\rightarrow I(W_{b})$$ 
As a result, we have an eigenspace decomposition
$$
I(W_{b})=I^{+}(W_{b},g_{b})\oplus I^{-}(W_b,g_{b}),
$$
where $I^{+}(W_{b},g_{b})$ and $ I^{-}(W_{b},g_{b})$ denotes the subspace respesnted by self-dual and anti-self dual harmonic 2-forms respectively. We use 
$$
P_{I^{+}}:I(W_{b})\rightarrow I^{+}(W_{b},g_{b})
$$
to denote the corresponding projection.

Recall that we chose a family $\widetilde{\omega}=\{\omega_{b}\}$ of compactly supported 2-forms on $W$ as a first perturbation of the family Seiberg-Witten equations (see Construction \ref{construction: perturbations on cobordism}). We make the following definition, which is analogous to Definition \ref{defi: admissible pair closed 4-manifolds} in the case of family cobordisms.
\begin{defi} Suppose $c_{1}(\mathfrak{s}_{Y_{i}})$ is torsion for $i=0,\ 1$. Then the pair $(g_{\widetilde{W}/\Q},\widetilde{\omega})$ is called \emph{admissible} (with respect to $\mathfrak{s}_{\widetilde{W}}$) if  the following inequality holds for any for any $b\in Q$:
$$
P_{I}(\omega_{b})+2\pi P_{I^{+}} c_{1}(\mathfrak{s}_{W}) \neq 0. 
$$
\end{defi}
Note that the family metrics $g_{\widetilde{W}/\Q}$ has been chosen so that its restriction to the tubular neighborhood $N_{i}$ is isometric to $[0,1]\times Y_{i}$. We can modify the cobordism by inserting another cylinder $[1,T]\times Y_{i}$ between $N_{i}$ and $W_{b}\setminus N_{i}$. On this new cylinder, we apply no perturbation on the Seiberg-Witten equations. So this modification will not change the admissibility condition.

\begin{lem}\label{lem: no reducible}
Consider any collections of perturbations 
$$
\{\omega_{b}\}_{b\in Q},\ \mathfrak{q}_{0},\ \mathfrak{q}_{1},\ \{\mathfrak{q}'_{0,b}\}_{b\in Q},\ \{\mathfrak{q}'_{1,b}\}_{b\in Q}.
$$
as in Construction \ref{construction: perturbations on cobordism}.
Suppose $(g_{\widetilde{W}/\Q},\widetilde{\omega})$ is admissible. Then after inserting the cylinder $[1,T]\times Y_{i}$ for $T\gg 0$, we may assume that the  perturbed Seiberg-Witten equation (\ref{eq: triply perturbed SW}) has no finite energy reducible solutions for any $b\in \Q$. 
\end{lem}
\begin{proof}
The case $\Q$ equals a point is proved in \cite[Proposition 27.2.4]{KronheimerMrowkaMonopole}. The key observation is that by inserting $[1,T]\times Y_{i}$, the perturbations introduced in step (ii) and (iii) of Construction \ref{construction: perturbations on cobordism} are pushed far away from the non-cylindrical part of $\widetilde{W}$. And in the limit case, the  $\widetilde{\omega}$-perturbed curvature equation (\ref{eq: perturbed curvature equation on cobordism}) has no reducible solution because $(g_{\widetilde{W}/\Q},\widetilde{\omega})$ is admissible. With this observation in mind, proof of the general case is the same.  
\end{proof}
Under the the perturbations and modified metric in Lemma \ref{lem: no reducible}, the maps $\widebar{m}^{*}_{*}$ are all zero since there are no reducible solutions. We consider the linear map
\begin{equation}\label{eq: mixed chain map}
\overrightarrow{m}=\left[\begin{array} {cc}
 m^{o}_{o}  & m^{u}_{o}  \\
 m^{o}_{s}  & m^{u}_{s}
\end{array}\right]:\quad \widehat{C}(Y_{0})\rightarrow \widecheck{C}(Y_{1}),
\end{equation}

\begin{lem}
$\overrightarrow{m}$ is a chain map. The induced map 
\begin{equation}\label{eq: mixed induced map}
    \overrightarrow{m}_{*}:\widehat{HM}_*(Y_{0},\mathfrak{s}_{Y_{0}})\rightarrow \widecheck{HM}_*(Y_{1},\mathfrak{s}_{Y_{1}})
\end{equation}
preserves the $U$-action. The composition $$j_{*,Y_{1}}\circ \overrightarrow{m}_{*}: \widehat{HM}_*(Y_{0},\mathfrak{s}_{Y_{0}})\rightarrow \widehat{HM}_*(Y_{1},\mathfrak{s}_{Y_{1}})$$ 
with the map $j_{1,Y_1}$ in (\ref{exact triangle}) equals the cobordism-induced map $\widehat{HM}_*(\widetilde{W},\mathfrak{s}_{\widetilde{W}})$.
\end{lem}
\begin{proof}
When $\Q$ equals a point, this is proved as \cite[Theorem 3.5.3]{KronheimerMrowkaMonopole}. The general case is essentially the same.
\end{proof}
Like the family Seiberg-Witten invariant, there is a chamber structure on the space of all admissible pairs and the induced map $\overrightarrow{m}_{*}$ depends on the choice of chamber. To explain this, we consider the positive cone
\begin{equation}\label{eq: positive cone}
V^{+}(W):=\{\mathfrak{a}\in I(W)\mid \langle \mathfrak{a}, \mathfrak{a}\rangle>0\}.
\end{equation}
where $\langle\cdot,\cdot\rangle$ is the intersection pairing defined in (\ref{eq: intersection pairing}). Because $\widetilde{W}$ is a homology product, we can canonically identify $V^{+}(W)$ with the corresponding space for each fiber $W_{b}$. 
Therefore, for any admissible pair $(g_{\widetilde{W}/\Q},\widetilde{\omega})$, we have a  map 
$$
b\mapsto P_{I}(\omega_{b})+2\pi P_{I^{+}} c_{1}(\mathfrak{s}_{W})\in V^{+}(W).
$$ 
This gives an element $\xi(g_{\widetilde{W}/\Q},\widetilde{\omega})\in [\Q,V^{+}(W)]$. We call this homotopy class the chamber of $(g_{\widetilde{W}/\Q},\widetilde{\omega})$. 
\begin{lem}\label{lem: chamber structure}
Any element $\xi\in [\Q,V^{+}(W)]$ can be realized as the chamber of an admissible pair $(g_{\widetilde{W}/\Q},\widetilde{\omega})$. Furthermore, the induced map $(\ref{eq: mixed induced map})$ only depends on the chamber. We use the notation
$$
\HMarrow_{*}(\widetilde{W},\mathfrak{s}_{\widetilde{W}})_{\xi}: \widehat{HM}(Y_{0},\mathfrak{s}_{Y_{0}})\rightarrow \widecheck{HM}(Y_{1},\mathfrak{s}_{Y_{1}})
$$
to denote this map.
\end{lem}
\begin{proof}
When $\Q$ is a single point and $b^{+}_{2}(W)\geq 2$, this is proved as Proposition 27.3.4 in \cite{KronheimerMrowkaMonopole}. To adapt the proof in our case, we note that any two admissible pairs $(g_{\widetilde{W}/\Q,0},\widetilde{\omega}_{0})$ and $(g_{\widetilde{W}/\Q,1},\widetilde{\omega}_{1})$ in the same chamber can be connected by a path $\{(g_{\widetilde{W}/\Q,t},\widetilde{\omega}_{t})\}_{t\in [0,1]}$ of admissible pairs. After a further small perturbation on this path, one can assume that the parmatrized moduli space over $[0,1]$ is regular and has no reducible solutions. The rest of the proof follows exactly the same argument as \cite{KronheimerMrowkaMonopole}. Note that just as Proposition \ref{pro: cobordism map well-defined}, proof of the independence with $(g_{Y_{i}},\mathfrak{q}_{i})$ relies on the composition law of the cobordism-induced maps (Theorem \ref{thm: composition law for HM}). 
\end{proof} 

Similar to the case of closed 4-manifolds, a chamber represented by a constant map is called a canonical chamber. Recall that $V^{+}(W)$ is homeomorphic to $(\mathbb{R}^{b^{+}_{2}(W)}-\{\overrightarrow{0}\})\times \mathbb{R}^{b^{-}(W)}$. So when $b^{+}_{2}(W)\geq 2$, there is a unique canonical chamber $\xi_{c}$. When $b^{+}_{2}(W)=1$, there are two canonical chambers. By picking an orientation on $H^{2}_{+}(W;\mathbb{R})$, we can talk about the positive canonical chamber $\xi^{+}_{c}$ and the negative canonical chamber $\xi^{-}_{c}$. To simply our notation, we also use 
$$
\HMarrow_{*}(\widetilde{W},\mathfrak{s}_{\widetilde{W}}),\quad \HMarrow_{*}(\widetilde{W},\mathfrak{s}_{\widetilde{W}})_{+},\quad \HMarrow_{*}(\widetilde{W},\mathfrak{s}_{\widetilde{W}})_{-}
$$
to denote the maps for $\xi_{c},\ \xi^{+}_{c}$ and $\xi_{c}^{-}$.

Now we turn to the case that $c_{1}(\mathfrak{s}_{Y_{k}})$ is nontorsion from $k=0$ or $1$. This case is actually easier: the unperturbed curvature equation on $Y_{k}$ has no solution. As a result, as along as one choose the 3-dimensional perturbation $\mathfrak{q}_{k}$ to be small enough, there will be no reducible critical point on $Y_{k}$, and hence no finite energy reducible solution on $W_{b}$. 
As a consequence, the map $\overrightarrow{m}$ defined in (\ref{eq: mixed chain map}) is always a chain map. With $\mathfrak{q}_{0},\mathfrak{q}_{1}$ fixed, any two choices of possible combinations of perturbations $\{\omega_{b}\}_{b\in Q},\  \{\mathfrak{q}'_{0,b}\}_{b\in Q},\ \{\mathfrak{q}'_{1,b}\}_{b\in Q}
$ can be connected by a generic path. Along this path, a reducible solution never appears and a standard cobordism argument shows that the induced map is independent with $\{\omega_{b}\}_{b\in Q},\  \{\mathfrak{q}'_{0,b}\}_{b\in Q},\ \{\mathfrak{q}'_{1,b}\}_{b\in Q}
$. The proof of the independence with $\mathfrak{q}_{0},\mathfrak{q}_{1}$ can be proceed in the same way as Lemma \ref{lem: chamber structure}. We denote the induced map $\overrightarrow{m}_{*}$ by 
$$
\HMarrow_*(\widetilde{W},\mathfrak{s}_{\widetilde{W}}): \widehat{HM}_*(Y_{0},\mathfrak{s}_{Y_{0}})\rightarrow \widecheck{HM}_*(Y_{1},\mathfrak{s}_{Y_{1}}).
$$
By our discussion above, this is an invariant of $(\widetilde{W},\mathfrak{s}_{\widetilde{W}})$. But it is actually not a new one: when $c_{1}(\mathfrak{s}_{Y_{k}})$ is nontorsion, the map $j_{*,Y_k}$ in (\ref{exact triangle}) is an isomorphism. And one has
$$
\HMarrow_*(\widetilde{W},\mathfrak{s}_{\widetilde{W}})=j^{-1}_{*,Y_1}\circ \widehat{HM}_{*}(\widetilde{W},\mathfrak{s}_{\widetilde{W}})
$$
when $c_{1}(\mathfrak{s}_{Y_{1}})$ is nontorsion and 
$$
\HMarrow(\widetilde{W},\mathfrak{s}_{\widetilde{W}})= \widecheck{HM}_{*}(\widetilde{W},\mathfrak{s}_{\widetilde{W}})\circ j^{-1}_{*,Y_0}
$$
when $c_{1}(\mathfrak{s}_{Y_{0}})$ is nontorsion. To unify our notation with the torsion case, when $c_{1}(\mathfrak{s}_{Y_{k}})$ is nontorsion for some $k$, we simply define $$\HMarrow_{*}(\widetilde{W},\mathfrak{s}_{\widetilde{W}})_{\xi}:=
\HMarrow_{*}(\widetilde{W},\mathfrak{s}_{\widetilde{W}})$$ for any $\xi\in [\Q,V^{+}(W)]$. 

By taking the map induced by $\overrightarrow{m}$ on the dual complexes, one can define the induced map on the monopole Floer cohomology
$$
\HMarrow^{*}(\widetilde{W},\mathfrak{s}_{\widetilde{W}})_{\xi}: \widehat{HM}^{*}(Y_{1},\mathfrak{s}_{Y_{1}})\rightarrow \widecheck{HM}^{*}(Y_{0},\mathfrak{s}_{Y_{0}}).$$
\subsection{Composition laws for cobordism-induced maps}\label{subsection: composition} In this section, we prove the composition laws for various maps induced by family of cobordisms (Theorem \ref{thm: composition law for HM} and Theorem \ref{thm: gluing for hm-arrow}). These results are proved in \cite{KronheimerMrowkaMonopole} when $\Q$ is a single point. So we only sketch the proofs, indicating what adaptions are needed in the family setting. This section will assume readers' familiarity with \cite{KronheimerMrowkaMonopole}.
\subsubsection{The gluing theorem for $\HMbar,\ \HMfrom$ and $\HMto$.}\label{section: composition}
We start by setting up some notations. Given $\spinc$ 3-manifolds  $(Y_{i},\mathfrak{s}_{Y_{i}})$ for $i=0,1,2$, we consider a family $\spinc$ cobordism $(\widetilde{W}_{01}/\Q,\mathfrak{s}_{\widetilde{W}_{01}})$ from $(Y_{0},\mathfrak{s}_{Y_{0}})$ to $(Y_{1},\mathfrak{s}_{Y_{1}})$, and a family $\spinc$ cobordism $(\widetilde{W}_{12}/\Q,\mathfrak{s}_{\widetilde{W}_{12}})$ from $(Y_{1},\mathfrak{s}_{Y_{1}})$ to $(Y_{2},\mathfrak{s}_{Y_{2}})$. We assume both families are homology products. By composing the cobordisms fiberwisely, we form a family cobordism $\widetilde{W}_{02}/\Q$ from $Y_0$ to $Y_2$. For $b\in \Q$, we use $W_{01,b},W_{12,b}$ and $W_{02,b}$ to denote the corresponding fibers over $b$. As before, we denote the fibers over the base point $b_{0}\in \Q$ by $W_{01},W_{12}$ and $W_{02}$. 

Note that we did not specify a way to the spinor bundles over $\Q\times Y_{1}$ so we do not have a single family $\spinc$-structure on $\widetilde{W}_{02}/\Q$. Instead, we consider the  set $\mathbb{S}$  of isomorphism classes of family $\spinc$ structures  $\mathfrak{s}$ on $\widetilde{W}_{02}$ that satisfies
$$\mathfrak{s}|_{\widetilde{W}_{01}}\cong \mathfrak{s}_{\widetilde{W}_{01}}\quad \text{and}\quad \mathfrak{s}|_{\widetilde{W}_{12}}\cong \mathfrak{s}_{\widetilde{W}_{12}}.$$ 
\begin{lem}\label{lem: S}
$\mathbb{S}$ is an affine set over image of the boundary map  
$$
\partial:H^{1}(Y_{1};\mathbb{Z})\rightarrow H^{2}(W_{02};\mathbb{Z})
$$
in the Mayer-Vietoris sequence arising from the decomposition $W_{02}=W_{01}\cup W_{12}$. In particular, $\mathbb{S}$ contains a single element if $Y_1$ is a rational homology 3-sphere.
\end{lem}
\begin{proof}
We let $P_{01}$ and $P_{12}$ be the principal $\textnormal{Spin}^{c}(4)$-bundle for the family $\spinc$ structure  $\mathfrak{s}_{\widetilde{W}_{01}}$ and  $\mathfrak{s}_{\widetilde{W}_{12}}$. We let $P_{1}$ be the principal bundle for the product family $\Q\times (Y_{1},\mathfrak{s}_{Y_{1}})$.  We pick specific identifications  $$\tau_{01}:P_{01}|_{\Q\times Y_{1}}\xrightarrow{\cong }P_{1}\text{ and }\tau_{12}:P_{12}|_{\Q\times Y_{1}}\xrightarrow{\cong }P_{1}.$$ 
Given any continuous map $f:\Q\times Y_{1}\rightarrow S^{1}$,  multiplying $f$ fiberwisely defines a bundle automorphism $\tau_{f}:P_{1}\rightarrow  P_{1}.$ We can construct a family $\spinc$ structure $\mathfrak{s}_{f}\in \mathbb{S}$ by gluing $P_{01}$ and $P_{12}$ together along $\Q\times Y_{1}$ using $\tau^{-1}_{12}\circ \tau_{f}\circ \tau_{01}$. Furthermore, it's straightforward to check that any element in $\mathbb{S}$ arises as some $
\mathfrak{s}_{f}$. 

Given $f,f'$, an isomorphism $\mathfrak{s}_{f}\xrightarrow{\cong} \mathfrak{s}_{f'}$ is given by gluing together bundle automorphisms  $$\tau_{g_{01}}:P_{01}\rightarrow P_{01}\quad\text{ and }\quad \tau_{g_{12}}:P_{12}\rightarrow P_{12}$$
which are defined as fiberwise multiplication with continuous functions $g_{01}:\widetilde{W}_{01}
\rightarrow S^{1}$ and  $g_{12}:\widetilde{W}_{12}
\rightarrow S^{1}$ respectively. They need to satisfy the compatibility condition \begin{equation}\label{eq: compatibility}
    (g_{12}|_{\Q\times Y_2})\cdot f=(g_{01}|_{\Q\times Y_2})\cdot f'.
\end{equation}
Here ``$\cdot$'' means pointwise multiplication. In other words, $\mathfrak{s}_{f}\cong \mathfrak{s}_{f'}$ if and only if there exists $g_{01},g_{12}$ that satisfies (\ref{eq: compatibility}). Note that $[A,S^{1}]\cong H^{1}(A;\mathbb{Z})$ for a general space $A$.
We obtain a bijection between $\mathbb{S}$ and cokernel of the map
$$
\tilde{i}_{12}:H^{1}(\widetilde{W}_{01};\mathbb{Z})\oplus H^{1}(\widetilde{W}_{12};\mathbb{Z})\rightarrow H^{1}(\Q\times Y_{1};\mathbb{Z})
$$
in the Mayer-Vietoris sequence for the decomposition $\widetilde{W}_{02}=\widetilde{W}_{01}\cup \widetilde{W}_{12}$. 

Note that the inclusion maps from $\Q\times Y_{1}$ to $\widetilde{W}_{01}$ and $\widetilde{W}_{12}$ preserve the fiber bundle structure. By the naturality of the Serre spectral sequence, we get isomorphisms \begin{equation*}
\begin{split}
H^{1}(\widetilde{W}_{01};\mathbb{Z})\cong H^{1}(W_{01};\mathbb{Z})\oplus H^{1}(\Q;\mathbb{Z}),\\
H^{1}(\widetilde{W}_{12};\mathbb{Z})\cong H^{1}(W_{12};\mathbb{Z})\oplus H^{1}(\Q;\mathbb{Z}),\\
H^{1}(\Q\times Y;\mathbb{Z})\cong H^{1}(Y;\mathbb{Z})\oplus H^{1}(\Q;\mathbb{Z}).
\end{split}
\end{equation*}
Therefore, the cokernel of $\tilde{i}_{12}$ is isomorphic to the cokernel of the map 
$$
i_{12}: H^{1}(W_{01};\mathbb{Z})\oplus H^{1}(W_{12};\mathbb{Z})\rightarrow H^{1}(Y_{1};\mathbb{Z})
$$
in the Mayer-Vietoris sequence of $W_{02}=W_{01}\cup W_{12}$. This is exactly the image of $\partial$.
\end{proof}

\begin{lem}\label{lem: finiteness for HM-to from and bar}
The map $HM^{\circ}_{*}(\widetilde{W}_{02},\mathfrak{s})$ is nontrivial for only finitely many $\mathfrak{s}\in \mathbb{S}$.
\end{lem}
\begin{proof}
There is nothing to prove if $\mathbb{S}$ is a finite set. Suppose $\mathbb{S}$ is an infinite set. Then by Lemma \ref{lem: S}, the map $$
\partial\otimes \mathbb{R}:H^{1}(Y_{2};\mathbb{R})\rightarrow H^{2}(W_{12};\mathbb{R})
$$
is nontrivial. By \cite[Lemma 3.9.2]{KronheimerMrowkaMonopole}, this implies that $b^{+}_{2}(W_{02})>0$. As a result, we may choose the perturbation forms $\{\omega_{b}\}$ in (\ref{construction: perturbations on cobordism}) such that the perturbed family Seiberg-Witten equations has no reducible solutions. In particular, the maps $\overline{m}^{s}_{s},\ \overline{m}^{u}_{s},\ \overline{m}^{s}_{u},\ \overline{m}^{u}_{u}$ are all vanishing. Furthermore, if 
if $[\mathfrak{a}]$ is boundary stable or $[\mathfrak{b}]$ is boundary unstable then we have 
\begin{equation}\label{eq: empty moduli space}
\mathcal{M}([\mathfrak{a}],\widetilde{W}^{*}_{02},\mathfrak{s},[\mathfrak{b}])=\emptyset
\end{equation}
because it can't contain an irreducible solution.

We first assume that $\Q$ is a single point. We fix a critical point $[\mathfrak{a}]$. Then \cite[Lemma 25.3.1]{KronheimerMrowkaMonopole} implies the following statement: There are only finitely many pairs $(\mathfrak{s},[\mathfrak{b}])$ such that both conditions are satisfied:
\begin{enumerate}[(i)]
    \item The moduli space $
\mathcal{M}([\mathfrak{a}],\widetilde{W}^{*}_{02},\mathfrak{s},[\mathfrak{b}])
$ has a (nonempty) $0$-dimensional component.
\item $[\mathfrak{b}]$ is either irreducible or boundary stable.
\end{enumerate}
In our case, (ii) is implied by (i) because of (\ref{eq: empty moduli space}). So we can conclude that for any $[\mathfrak{a}]$, there are only finitely many $\mathfrak{s}$ such that the corresponding maps $m^{*}_{*}$ do not all vanish on $[\mathfrak{a}]$. From this, we can further deduce that for any fixed element $a\in HM^{\circ}_{*}(Y_{0},\mathfrak{s}_{Y_{0}})$, there are only finitely many $\mathfrak{s}\in \mathbb{S}$ such that $HM^{\circ}_{*}(\widetilde{W}_{02},\mathfrak{s})(a)\neq 0$. This finishes the proof for $\HMfrom$ and $\HMbar$ because they are finitely generated as $\mathbb{Q}[U]$-modules. The same is true for $\HMto$ if $c_{1}(\mathfrak{s}_{Y_0})$ is nontorsion. When $c_{1}(\mathfrak{s}_{Y_0})$ is torsion, we may decompose $\HMto_{*}(Y_{0},\mathfrak{s}_{Y_{0}})$ as $\HMto_{\leq N}(Y_{0},\mathfrak{s}_{Y_{0}})\oplus \HMto_{> N}(Y_{0},\mathfrak{s}_{Y_{0}})$ according to the rational grading. The first summand is a finite dimensional vector space so our old argument works. Elements in the second summand can be represented by linear combinations of boundary stable critical points with grading greater than $N$. When $N$ is large enough, these critical points are all annihilated by the map $\overline{\partial}^{s}_{u}$. By (\ref{eq: mcheck}), we see that $\HMto_{> N}(Y_{0},\mathfrak{s}_{Y_{0}})$ is in the kernel of $\HMto(\widetilde{W}_{02},\mathfrak{s})$ for any $\mathfrak{s}\in \mathbb{S}$. This finishes the proof when $\Q$ is a single point. 

The above argument relies on the finiteness result \cite[Lemma 25.3.1]{KronheimerMrowkaMonopole}, which is a consequence of the compactness theorem \cite[Proposition 24.6.4]{KronheimerMrowkaMonopole}. For a general compact base $\Q$, this compactness theorem (and hence the finiteness result) can be adapted  without essential difficulty. The rest of the proof can be proceeded in the same way. 
\end{proof}
Now we state the first composition law.
\begin{thm}\label{thm: composition law for HM}
(1)
 Suppose $(\widetilde{W}_{01},\mathfrak{s}_{\widetilde{W}_{01}})$ is the product family associated to a single $\spinc$ cobordism $(W_{01},\mathfrak{s}_{W_{01}})$. Then we have 
$$
\sum_{\mathfrak{s}\in \mathbb{S}}\HMfrom_{*}(\widetilde{W}_{02},\mathfrak{s})= \HMfrom_{*}(W_{01},\mathfrak{s}_{W_{01}})\circ \HMfrom_{*}(\widetilde{W}_{12},\mathfrak{s}_{\widetilde{W}_{12}})
$$

(2) Suppose $(\widetilde{W}_{12},\mathfrak{s}_{\widetilde{W}_{12}})$ is the product family associated to a $\spinc$ cobordism $(W_{12},\mathfrak{s}_{W_{12}})$. Then we have 
$$
\sum_{\mathfrak{s}\in \mathbb{S}}\HMfrom_{*}(\widetilde{W}_{02},\mathfrak{s})= \HMfrom_{*}(\widetilde{W}_{01},\mathfrak{s}_{\widetilde{W}_{01}})\circ \HMfrom_{*}(W_{12},\mathfrak{s}_{W_{12}}).
$$

(3) Parallel results hold for $\HMbar$ and $\HMto$.
\end{thm} 
Theorem \ref{thm: composition law for HM} is a family version of \cite[Proposition 26.1.2]{KronheimerMrowkaMonopole} and the proof here is adapted from there. All these adaptions are inessential (although some are not completely obvious). For the sake of completeness, we will sketch the major steps of the proof in the family setting, emphasizing all the adaptions we need. 

We will be dealing with two similar but different situations (either $\widetilde{W}_{01}$ or $\widetilde{W}_{12}$ is a product family). So we first set up a general framework without assuming that any of them is a product. Then we specify to one of the two cases and discuss extra constrains coming from the product structure.
\subsubsection{The general setup}
We use $N_{0},N_{1}$ to denote collar neighborhoods of the boundary component $Y_{0},\ Y_{1}$ in the fibers of $\widetilde{W}_{01}$ and we use $N'_{1},N_{2}$ to denote collar neighborhoods of the boundary component $Y_{1},\ Y_{2}$ in the fibers of  $\widetilde{W}_{12}$. We equip $\widetilde{W}_{01}$ and $\widetilde{W}_{12}$ with family metrics $g_{\widetilde{W}_{01}}$ and $g_{\widetilde{W}_{12}}$ which are isomorphic to $[0,1]\times Y_{*}$ when restricts to these collar neighborhoods. For any $S>0$, we consider the family cobordism  
$
\widetilde{W}(S)$ obtained by inserting a neck of length $2S$ in the middle. So the fiber of $\widetilde{W}(S)$ over $b$ is given by
$$W(S)_{b}:=W_{01,b}\cup ([-S,S]\times Y_{1})\cup W_{12,b}.
$$
By further attaching cylindrical ends $(-\infty,0]\times Y_{0}$ and $[0,\infty)\times Y_{1}$ on each fiber, we obtain the family  
$
\widetilde{W}(S)^{*}
$.

\begin{cons}\label{construction: perturbations on composed cobordism} We specify a set of perturbations on $\widetilde{W}(S)^{*}$ in three steps. (c.f. Construction \ref{construction: perturbations on cobordism}.) 
\begin{enumerate}[(i)]
    \item We pick smooth families $\{\omega^{01}_{b}\}$,  $\{\omega^{12}_{b}\}$ of self-dual 2-forms $\widetilde{W}_{01}$ and $\widetilde{W}_{12}$, with the requirement that $\omega^{01}_{b}$ and $\omega^{12}_{b}$ vanishes on $N_{0}\cup N_{1}$ and $N_{1}'\cup N_{2}$ respectively.
    \item Then we pick tame, admissible perturbations $\mathfrak{q}_{0},\ \mathfrak{q}_{1},\ \mathfrak{q}_{2}$ on $Y_{0},\ Y_{1},\ Y_{2}$. We add the corresponding 4-dimensional perturbations $\hat{\mathfrak{q}}_{0},\ \hat{\mathfrak{q}}_{1},\ \hat{\mathfrak{q}}_{2}$ on the cylindrical parts $(-\infty,0]\times Y_{0},\ [-S,S]\times Y_{1},\ [0,\infty)\times Y_2$ of each fiber.
    \item We pick families of 3-dimensional perturbations $\{\mathfrak{p}'_{0,S,b}\},\ \{\mathfrak{p}'_{1,S,b}\},\ \{\mathfrak{p}''_{1,S,b}\},\ \{\mathfrak{p}'_{2,S,b}\}$ on $Y_{0},\ Y_{1},\ Y_{1},\ Y_{2}$. Similar to (\ref{eq: perturbation near boundary}), we then add the following 4-dimensional perturbations 
\begin{equation*}
\begin{split}
    \beta(t)\hat{\mathfrak{q}}_{0}+\beta_{0}(t)\hat{\mathfrak{p}}'_{0,S,b},\quad \beta(t)\hat{\mathfrak{q}}_{1}+\beta_{0}(t)\hat{\mathfrak{p}}'_{1,S,b}\\
    \beta(t)\hat{\mathfrak{q}}_{1}+\beta_{0}(t)\hat{\mathfrak{p}}''_{1,S,b},\quad \beta(t)\hat{\mathfrak{q}}_{2}+\beta_{0}(t)\hat{\mathfrak{p}}'_{2,S,b}
    \end{split}
\end{equation*}
on $N_{0},\ N_{1},\ N'_{1}$ and $N_{2}$ respectively.
\end{enumerate}
\end{cons}

Let $[\mathfrak{a}]$ and $[\mathfrak{b}]$ be a pair of critical points for $(Y_{0},\mathfrak{s}_{Y_{0}})$ and $(Y_{2},\mathfrak{s}_{Y_{2}})$. For any homotopy class 
$$z\in \bigcup_{\mathfrak{s}\in \mathbb{S}}\pi_{0}(\mathcal{B}^{\sigma}([\mathfrak{a}],\widetilde{W}^{*}_{02},\mathfrak{s},[\mathfrak{b}]),$$
we consider the moduli space
\begin{equation}\label{eq: paramatrized moduli space}
\mathcal{M}_{z}([\mathfrak{a}],[\mathfrak{b}])=\bigcup_{S\in [0,+\infty)}\{S\}\times \mathcal{M}_{z}([\mathfrak{a}],\widetilde{W}(S)^{*},[\mathfrak{b}]),
\end{equation}
which admits an obvious projection map to $[0,+\infty)$. By the same argument as \cite[Proposition 26.1.3]{KronheimerMrowkaMonopole}, one can prove that under suitable perturbations, the moduli space $\mathcal{M}_{z}([\mathfrak{a}],[\mathfrak{b}])$ is a smooth manifold with boundary for any $z,[\mathfrak{a}],[\mathfrak{b}]$ and the boundary is the fiber over $S=0$. 

To compactify $\mathcal{M}_{z}([\mathfrak{a}],[\mathfrak{b}])$, we first compactify each fiber $$\{S\}\times \mathcal{M}_{z}([\mathfrak{a}],\widetilde{W}(S)^*,[\mathfrak{b}])$$ into $$\{S\}\times \mathcal{M}^{+}_{z}([\mathfrak{a}],\widetilde{W}(S)^*,[\mathfrak{b}]),$$ by adding in broken trajectories on the fibers $\{W(S)_{b}\}$ in the homotopy class $z$. The result is still not compact due the non-compactness of the base $[0,+\infty)$. So we further add the fiber over $S=\infty$, which we denote by  $\mathcal{M}^{+}_{z}([\mathfrak{a}],\widetilde{W}(\infty)^*,[\mathfrak{b}])$. An element in this space is a quintuple $([\boldsymbol{\breve{\gamma}}_{0}],[\gamma_{01}],[\boldsymbol{\breve{\gamma}}_{1}],[\gamma_{12}],[\boldsymbol{\breve{\gamma}}_{2}])$, where 
\begin{equation*}
    \begin{split}
        &[\boldsymbol{\breve{\gamma}}_{0}]\in \breve{\mathcal{M}}^{+}([\mathfrak{a}],[\mathfrak{a}_{0}]),\\
     &[\gamma_{01}]\in \mathcal{M}([\mathfrak{a}_{0}],W_{01,b}^{*},[\mathfrak{a}_{1}]),
     \\
     &[\boldsymbol{\breve{\gamma}}_{1}]\in \breve{\mathcal{M}}^{+}([\mathfrak{a}_{1}],[\mathfrak{b}_{1}]),
     \\
     &[\gamma_{12}]\in \mathcal{M}([\mathfrak{b}_{1}],W_{12,b}^{*},[\mathfrak{b}_{2}]),
     \\
      &[\boldsymbol{\breve{\gamma}}_{2}]\in \breve{\mathcal{M}}^{+}([\mathfrak{b}_{2}],[\mathfrak{b}]).
    \end{split}
\end{equation*}
Here $[\mathfrak{a}_{i}]$ and $[\mathfrak{b}_{i}]$ are critical points on $Y_{i}$ and $\breve{\mathcal{M}}^{+}(-)$ denotes the moduli space of broken trajectories on $Y_{i}$.  (We allow an unbroken constant trajectory in this moduli space.) We require that the homotopy classes compose to give $z$ and that $[\gamma_{01}]$ and $[\gamma_{12}]$ are over the same point $b\in \Q$. Note that on the fiber over $S=\infty$, we are using the same perturbations $\{\omega^{01}_{b}\}$,  $\{\omega^{12}_{b}\}$, $\mathfrak{q}_{0},\ \mathfrak{q}_{1},\ \mathfrak{q}_{2}$ as finite $S$. And we set $\{\mathfrak{p}'_{0,\infty,b}\},\ \{\mathfrak{p}'_{1,\infty,b}\},\ \{\mathfrak{p}''_{1,\infty,b}\},\ \{\mathfrak{p}'_{2,\infty,b}\}$ to be limits of $\{\mathfrak{p}'_{0,S,b}\},\ \{\mathfrak{p}'_{1,S,b}\},\ \{\mathfrak{p}''_{1,S,b}\},\ \{\mathfrak{p}'_{2,S,b}\}$ as $S\rightarrow \infty$. (Perturbations for finite $S$ are chosen such that these limits exist.)

We use $\mathcal{M}^{+}_{z}([\mathfrak{a}],[\mathfrak{b}])$ to denote the space
$$
\mathcal{M}^{+}_{z}([\mathfrak{a}],[\mathfrak{b}])=\bigcup_{S\in [0,+\infty]}\{S\}\times \mathcal{M}^{+}_{z}([\mathfrak{a}],\widetilde{W}(S)^{*},[\mathfrak{b}]).
$$
When $\Q$ is a point, the topology on $\mathcal{M}^{+}_{z}([\mathfrak{a}],[\mathfrak{b}])$ is defined in \cite[Page 539]{KronheimerMrowkaMonopole}. The topology for a general $\Q$ can be defined similarly. The following theorem is a family version of \cite{KronheimerMrowkaMonopole}. The proof is essential identical so we omit it.
\begin{thm}
For any $[\mathfrak{a}]$, $[\mathfrak{b}]$, $z$, the moduli space $\mathcal{M}^{+}_{z}([\mathfrak{a}],[\mathfrak{b}])$ is compact.
\end{thm}
There is the following finiteness result, which is a family version of \cite[Lemma 26.1.5]{KronheimerMrowkaMonopole}. Again, the proof is essentially the same and omitted. 
\begin{lem}\label{lem: finiteness}
For any $[\mathfrak{a}]$, $d_{0}$ and $i_{0}$, there are only finitely many pairs $(z,[\mathfrak{b}])$ such that
\begin{enumerate} [(i)]
    \item the moduli space $\mathcal{M}^{+}_{z}([\mathfrak{a}],[\mathfrak{b}])$ is non-empty;
    \item the dimension of the moduli space is at most $d_0$;
    \item if $\mathfrak{s}_{Y_{2}}$ is torsion, then $\gr^{\mathbb{Q}}([\mathfrak{b}])\geq i_{0}$.
\end{enumerate}
\end{lem}

The moduli space $\mathcal{M}^{+}_{z}([\mathfrak{a}],[\mathfrak{b}])$ is a stratified space and there are two types of strata:
\begin{enumerate}[(I)]
    \item The first type of strata do not intersect the fiber over $S=\infty$. These strata are products of moduli spaces of unbroken solutions. A typical example is 
$$
\breve{\mathcal{M}}_{z_{0}}([\mathfrak{a}],[\mathfrak{a}_{0}])\times \mathcal{M}_{z_{1}}([\mathfrak{a}_{0}],[\mathfrak{b}_{2}])\times \breve{\mathcal{M}}_{z_{2}}([\mathfrak{b}_{2}],[\mathfrak{b}])
$$
for some critical points $[\mathfrak{a}_{0}]$ and $[\mathfrak{b}_{2}]$ on $Y_{0}$ and $Y_{2}$ respectively. 
\item The second type of strata are contained in the fiber over $S=\infty$. These strata involve fiber products over $\Q$. A typical exam of such strata is
$$ \mathcal{M}_{z_{0}}([\mathfrak{a}],\widetilde{W}_{01}^{*},[\mathfrak{a}_{1}])
\times_ {\Q} (\breve{\mathcal{M}}_{z_{1}}([\mathfrak{a}_{1}],[\mathfrak{b}_{1}])\times \Q)\times_{\Q} \mathcal{M}_{z_{2}}([\mathfrak{b}_{1}],\widetilde{W}_{12}^{*},[\mathfrak{b}]).$$
for some critical points $[\mathfrak{a}_{1}]$ and $[\mathfrak{b}_{1}]$ on $Y_{1}$.
\end{enumerate}
\begin{defi}\label{defi: parametrized moduli space regular}
We say the moduli space $\mathcal{M}^{+}_{z}([\mathfrak{a}],[\mathfrak{b}])$ is regular if the following conditions are both satisfied:

\begin{enumerate}
    \item All factors of all these strata are regular (see \cite[Definition 24.4.9]{KronheimerMrowkaMonopole}).
    \item The projection map 
    $$
    \mathcal{M}_{z_{0}}([\mathfrak{a}_{0}],\widetilde{W}_{01}^{*},[\mathfrak{a}_{1}])\times \mathcal{M}_{z_{1}}([\mathfrak{b}_{1}],\widetilde{W}_{12}^{*},[\mathfrak{b}_{2}])\rightarrow \Q\times \Q
    $$
is transverse to the diagonal whenever $\mathcal{M}_{z_{0}}([\mathfrak{a}_{0}],\widetilde{W}_{01}^{*},[\mathfrak{a}_{1}])$ and $\mathcal{M}_{z_{1}}([\mathfrak{b}_{1}],\widetilde{W}_{12}^{*},[\mathfrak{b}_{2}])$ appear as two factors of a same strata of type (II).
\end{enumerate}
\end{defi}

The following result is a family version of \cite[Proposition 26.1.6]{KronheimerMrowkaMonopole}. The proof is essential identical.
\begin{pro}\label{pro: boundary strata}
Assume $\mathcal{M}^{+}_{z}([\mathfrak{a}],[\mathfrak{b}])$ is regular. And suppose that the moduli space $\mathcal{M}_{z}([\mathfrak{a}],[\mathfrak{b}])$ contains irreducible solutions and has dimension $d+1$. Then $\mathcal{M}^{+}_{z}([\mathfrak{a}],[\mathfrak{b}])$ is a space stratified by manifolds and has $\mathcal{M}_{z}([\mathfrak{a}],[\mathfrak{b}])$ as its top stratum. In the fiber over $S=\infty$, the strata of dimension $d$ are the top strata of the following pieces of $\mathcal{M}^{+}_{z}([\mathfrak{a}],\widetilde{W}(\infty)^*,[\mathfrak{b}])$:
\begin{enumerate}[(i)]
    \item $\widetilde{\mathcal{M}}_{01}\times_{\Q}\widetilde{\mathcal{M}}_{12}$,
    \item $\breve{\mathcal{M}}_{0}\times (\widetilde{\mathcal{M}}_{01}\times_{\Q}\widetilde{\mathcal{M}}_{12})$,
    \item $\widetilde{\mathcal{M}}_{01}\times_{\Q}(\breve{\mathcal{M}}_{1}\times \Q)\times_{\Q}\widetilde{\mathcal{M}}_{12}$,
    \item $(\widetilde{\mathcal{M}}_{01}\times_{\Q}\widetilde{\mathcal{M}}_{12})\times \breve{\mathcal{M}}_{2}$,
\end{enumerate}
where $\widetilde{\mathcal{M}}_{ij}$ denotes a typical moduli space on $\widetilde{W}_{ij}$ and $\breve{\mathcal{M}}_{i}$ denotes a typical moduli space (of unparametrized, possibly broken trajectories) on $\mathbb{R}\times Y_{i}$. In the last three cases, the middle one  in the three  moduli spaces is boundary-obstructed. In The moduli space has a codimension-1 $\delta$-structure along the strata in the last three cases, and is a manifold with boundary along the stratum in the first case. 
\end{pro}

\subsubsection{Proof of Theorem \ref{thm: composition law for HM}} Our argument will be an adaption of the argument in \cite[Section 26.2]{KronheimerMrowkaMonopole} to the family setting. We focus on Part (1) as the other parts are essential the same. For rest of the proof, we assume $(\widetilde{W}_{01},\mathfrak{s}_{\widetilde{W}_{01}})$ is a product family associated to a single cobordism $(W_{01},\mathfrak{s}_{W_{01}})$. 

\begin{cons}\label{cons: gluing perturbation}
We start by choosing appropriate perturbations in this setting:
\begin{enumerate}
    \item Since we are dealing with induced maps on $\HMto$, by Remark \ref{rmk: perturbation 2-forms}, we may set the perturbation 2-forms $\{\omega_{01,b}\}$ and $\{\omega_{12,b}\}$ to be identically zero. 
    \item We choose perturbations $p''_{1,\infty}$ (resp. $p''_{2,\infty}$) and on $Y_{1}$ (resp. $Y_{2}$) such that the moduli space 
$
\mathcal{M}_{z}([\mathfrak{a}],W^{*}_{12},[\mathfrak{b}])
$ is regular for any $[\mathfrak{a}]$, $[\mathfrak{b}]$ and $z$.
\item We define the family perturbations $\{p''_{1,\infty,b}\}$ and $\{p'_{2,\infty,b}\}$ by setting $p''_{1,\infty,b}=p''_{1,\infty}$ and $p'_{2,\infty,b}=p'_{2,\infty}$ for any $b\in \Q$. Then the family moduli space $\mathcal{M}_{z}([\mathfrak{a}],\widetilde{W}^{*}_{12},[\mathfrak{b}])$ is automatically regular and we have the diffeomorphism 
\begin{equation}\label{eq: moduli space product}
 \mathcal{M}_{z}([\mathfrak{a}],\widetilde{W}^{*}_{12},[\mathfrak{b}])\cong \Q\times 
\mathcal{M}_{z}([\mathfrak{a}],W^{*}_{12},[\mathfrak{b}]).   
\end{equation}
\item We also choose family perturbations $\{p'_{0,\infty,b}\}$ and $\{p'_{1,\infty,b}\}$ on $Y_{0}$ and $Y_{1}$ such that the moduli space 
$
\mathcal{M}_{z}([\mathfrak{a}],\widetilde{W}^{*}_{12},[\mathfrak{b}])
$
is always regular.
\item On the fiber over $S=0$, we pick family perturbations $\{p'_{0,0,b}\}$, $\{p'_{1,0,b}\}$, $\{p''_{1,0,b}\}$ and $\{p'_{2,0,b}\}$ such that the moduli space 
$
\mathcal{M}_{z}([\mathfrak{a}],\widetilde{W}^{*}_{02},[\mathfrak{b}])
$
is regular for any $[\mathfrak{a}],\ [\mathfrak{b}]$ and $z$. (Recall that $\widetilde{W}(0)^{*}=\widetilde{W}^{*}_{02}$.)
\item By the same argument as \cite[Proposition 24.4.10]{KronheimerMrowkaMonopole}, we can find family perturbations $\{p'_{0,S,b}\}$, $\{p'_{1,S,b}\}$, $\{p''_{1,S,b}\}$ and $\{p'_{2,S,b}\}$ that extend the corresponding perturbations at $S=0$ and $S=\infty$ and make the compactified parametrized moduli space $\mathcal{M}^{+}_{z}([[\mathfrak{a}]],[\mathfrak{b}])$ regular for any $[\mathfrak{a}]$, $[\mathfrak{b}]$ and $z$. Note that the second condition in Definition \ref{defi: parametrized moduli space regular} is automatically satisfied because of equation (\ref{eq: moduli space product}).
\end{enumerate}
\end{cons}
Again because of (\ref{eq: moduli space product}), the four types of strata in Proposition \ref{pro: boundary strata} become the following strata
\begin{equation}\label{eq: boundary strata product form}
\begin{split}
   (i)\ \mathcal{M}_{01}\times\widetilde{\mathcal{M}}_{12},\quad &(ii)\ \breve{\mathcal{M}}_{0}\times \mathcal{M}_{01}\times\widetilde{\mathcal{M}}_{12}\\
   (iii)\ \mathcal{M}_{01}\times\breve{\mathcal{M}}_{1}\times\widetilde{\mathcal{M}}_{12},\quad &(iv)\  \mathcal{M}_{01}\times\widetilde{\mathcal{M}}_{12}\times \breve{\mathcal{M}}_{2}.
   \end{split}
\end{equation}
Here $\mathcal{M}_{01}$ denotes a typical moduli space on the single cobordism $W_{01}$.

We have the following Proposition regarding the boundary orientations of these strata. 

\begin{pro}\label{pro: orientation boundary strata}
Suppose the moduli space $\mathcal{M}_{z}([\mathfrak{a}],[\mathfrak{b}])$ in Proposition \ref{pro: boundary strata} contains irreducible solutions and has dimension $d+1$. Then the four types of strata (\ref{eq: boundary strata product form}) differ from their boundary orientations by the signs 
\begin{enumerate}[(i)]
    \item $1\cdot (-1)^{d_{01}\cdot \text{dim}(\Q)}$
    \item $(-1)^{d_{0}+d_{01}+1}\cdot (-1)^{(d_{0}+d_{01}-1)\cdot \text{dim}(\Q)}$
    \item $(-1)^{d_{1}}\cdot (-1)^{(d_{01}+d_{1}-1)\cdot \text{dim}(\Q)}$
    \item $(-1)\cdot (-1)^{d_{01}\cdot \text{dim}(\Q)}$
\end{enumerate}
accordingly, where $d_{0}$ and $d_{1}$ are the dimensions of the parametrized moduli spaces $\mathcal{M}_0$ and $\mathcal{M}_{1}$, (so the dimension of $\breve{\mathcal{M}}_{i}$  equals $d_{i}-1$) and $d_{01}$ is the dimension of $\mathcal{M}_{01}$.
\end{pro}
\begin{proof}
This is a family version of \cite[Proposition 26.1.7]{KronheimerMrowkaMonopole}. One can give a prove by unwinding the gluing construction and adapting the corresponding arguments in \cite{KronheimerMrowkaMonopole}. But the argument is long and tedious. Instead, we give an indirect proof here. 

The discrepancy between the orientation of the boundary strata and the boundary orientation is caused by switching various factors in the gluing construction. Therefore, there are universal formulas that express the sign differences in terms of $\operatorname{dim}\Q$ and dimension of the various moduli spaces. To find out these formulas, we consider the special case $\widetilde{W}_{12}=\Q\times {W}_{12}$, which also implies that $\widetilde{W}_{02}=\Q\times {W}_{02}$. 
Therefore, we have orientation preserving homeomorphisms
$$
\widetilde{M}_{12}\cong \Q\times M_{12},\quad \widetilde{M}_{02}\cong \Q\times M_{02}.  
$$
(Recall that we orient the moduli spaces by putting the $\Q$-coordinates first.)
By \cite[Proposition 26.1.7]{KronheimerMrowkaMonopole} , the orientation of the boundary strata 
$
\mathcal{M}_{01}\times\mathcal{M}_{12}\times \breve{\mathcal{M}}_{2}
$ differs from the boundary orientation of $\bigcup_{S\in [0,+\infty]}\{S\}\times \mathcal{M}^{+}_{z}([\mathfrak{a}],W(S)^{*},[\mathfrak{b}])$ by the sign $(-1)$. As a result, the orientation of the boundary strata 
$$\mathcal{M}_{01}\times\widetilde{\mathcal{M}}_{12}\times \breve{\mathcal{M}}_{2}
\cong \mathcal{M}_{01}\times \Q\times \mathcal{M}_{12}\times \breve{\mathcal{M}}_{2}$$ differs from the boundary orientation of the parametrized moduli space $$
\mathcal{M}^{+}_{z}([\mathfrak{a}],[\mathfrak{b}])\cong \Q\times (\bigcup_{S\in [0,+\infty]}\{S\}\times \mathcal{M}^{+}_{z}([\mathfrak{a}],W(S)^{*},[\mathfrak{b}]))
$$
by the sign $(-1)\cdot (-1)^{d_{01}\cdot \text{dim}(\Q)}$. This finishes the proof of (iv). And the other cases are similar.
\end{proof}

Now we define maps 
\begin{equation*}
\begin{split}
K^{o}_{o}: C^{o}_{\bullet}(Y_{0},\mathfrak{s}_{Y_{0}})\rightarrow C^{o}_{\bullet}(Y_{2},\mathfrak{s}_{Y_{2}})\\
K^{o}_{s}: C^{o}_{\bullet}(Y_{0},\mathfrak{s}_{Y_{0}})\rightarrow C^{s}_{\bullet}(Y_{2},\mathfrak{s}_{Y_{2}})\\
K^{s}_{o}: C^{o}_{\bullet}(Y_{0},\mathfrak{s}_{Y_{0}})\rightarrow C^{o}_{\bullet}(Y_{2},\mathfrak{s}_{Y_{2}})
\\
K^{s}_{s}: C^{o}_{\bullet}(Y_{0},\mathfrak{s}_{Y_{0}})\rightarrow C^{o}_{\bullet}(Y_{2},\mathfrak{s}_{Y_{2}})
\end{split}
\end{equation*}
by counting points in the zero-dimensional moduli spaces $\mathcal{M}_{z}([\mathfrak{a}],[\mathfrak{b}])$. More precisely, we define
\begin{equation}\label{eq: operators K}
\begin{split}
K^{o}_{o}(e_{[\mathfrak{a}]})=\sum_{[\mathfrak{b}]\in \mathfrak{C}^{o}(Y_{2},\mathfrak{s}_{Y_2})}\sum_{z}\# \mathcal{M}_{z}([\mathfrak{a}],[\mathfrak{b}])\cdot e_{[\mathfrak{b}]}\\
K^{o}_{s}(e_{[\mathfrak{a}]})=\sum_{[\mathfrak{b}]\in \mathfrak{C}^{s}(Y_{2},\mathfrak{s}_{Y_2})}\sum_{z}\# \mathcal{M}_{z}([\mathfrak{a}],[\mathfrak{b}])\cdot e_{[\mathfrak{b}]}\\
K^{s}_{o}(e_{[\mathfrak{a}]})=\sum_{[\mathfrak{b}]\in \mathfrak{C}^{o}(Y_{2},\mathfrak{s}_{Y_2})}\sum_{z}\# \mathcal{M}_{z}([\mathfrak{a}],[\mathfrak{b}])\cdot e_{[\mathfrak{b}]}\\
K^{s}_{s}(e_{[\mathfrak{a}]})=\sum_{[\mathfrak{b}]\in \mathfrak{C}^{s}(Y_{2},\mathfrak{s}_{Y_2})}\sum_{z}\# \mathcal{M}_{z}([\mathfrak{a}],[\mathfrak{b}])\cdot e_{[\mathfrak{b}]}
\end{split}
\end{equation}
Here we are summing over all homotopy classes $$z\in \bigcup_{\mathfrak{s}\in \mathbb{S}}\pi_{0}(\mathcal{B}^{\sigma}([\mathfrak{a}],\widetilde{W}^{*}_{02},\mathfrak{s},[\mathfrak{b}]),$$ such that $\mathcal{M}_{z}([\mathfrak{a}],[\mathfrak{b}])$ is $0$-dimensional. The sum may not be finite but the finiteness result Lemma \ref{lem: finiteness} ensures that these maps are still well-defined because we have passed to the completion. Similarly, we can define operators $\widebar{K}^{u}_{s},\ \widebar{K}^{u}_{u},\ \widebar{K}^{s}_{u},\ \widebar{K}^{s}_{s}$ by counting reducible solutions only.

As in \cite{KronheimerMrowkaMonopole}, we package these operators and obtain an operator 
$$
\widecheck{K}: \widecheck{C}_{\bullet}(Y_{0},\mathfrak{s}_{Y_{0}})\rightarrow \widecheck{C}_{\bullet}(Y_{2},\mathfrak{s}_{Y_{2}})
$$
defined by
\begin{equation}
    \widecheck{K}=\left[\begin{array} {cc}
 K^{o}_{o}  & -K^{u}_{o}\widebar{\partial}^{s}_{u}-m^{u}_{o}\widebar{m}^{u}_{s}-\partial^{u}_{o}\widebar{K}^{u}_{s}  \\
 K^{o}_{s}  & \widebar{K}^{s}_{s}-K^{u}_{s}\widebar{\partial}^{s}_{u}-m^{u}_{s}\widebar{m}^{s}_{u}-\partial^{u}_{s}\widebar{K}^{s}_{u}
\end{array}\right].
\end{equation}
Here $m^{u}_{o}$ and $m^{u}_{s}$ are defined in (\ref{eq: m operator}) using the single cobordism $(W_{01},\mathfrak{s}_{W_{01}})$ with perturbations $\mathfrak{p}''_{1,\infty}$ and $\mathfrak{p}'_{2,\infty}$ introduced in Construction \ref{cons: gluing perturbation} (3). And  $\overline{m}^{u}_{s}$ and $\overline{m}^{s}_{u}$ are defined in (\ref{eq: m-bar operator}) using the family cobordism $(\widetilde{W}_{12},\mathfrak{s}_{\widetilde{W}_{12}})$, with perturbations $\{\mathfrak{p}'_{0,\infty,b}\}$ and $\{\mathfrak{p}'_{1,\infty,b}\}$ introduced in Construction \ref{cons: gluing perturbation} (4).

Using the perturbations specified in Construction \ref{cons: gluing perturbation}, we can define the operator
$$
\widecheck{m}_{01}:\widecheck{C}_{\bullet}(Y_{0},\mathfrak{s}_{Y_{0}})\rightarrow \widecheck{C}_{\bullet}(Y_{1},\mathfrak{s}_{Y_{1}})
$$
for the cobordism $(W_{01},\mathfrak{s}_{W_{01}})$, as well as the operator 
$$
\widecheck{m}_{12}:\widecheck{C}_{\bullet}(Y_{1},\mathfrak{s}_{Y_{1}})\rightarrow \widecheck{C}_{\bullet}(Y_{2},\mathfrak{s}_{Y_{2}})
$$
for the (family) cobordism   $(\widetilde{W}_{12},\mathfrak{s}_{\widetilde{W}_{12}})$. Furthermore, for any $\mathfrak{s}\in \mathbb{S}$, we define the operator 
$$
\widecheck{m}_{02}:\widecheck{C}_{\bullet}(Y_{0},\mathfrak{s}_{Y_{0}})\rightarrow \widecheck{C}_{\bullet}(Y_{2},\mathfrak{s}_{Y_{2}})
$$
for the family $\spinc$-cobordism $(\widetilde{W}_{02},\mathfrak{s})$. By summing these operators for all $\mathfrak{s}\in \mathbb{S}$, we get the operator 
$$
\widecheck{m}_{02,\mathbb{S}}:\widecheck{C}_{\bullet}(Y_{0},\mathfrak{s}_{Y_{0}})\rightarrow \widecheck{C}_{\bullet}(Y_{2},\mathfrak{s}_{Y_{2}})
$$

\begin{pro}\label{pro: gluing chain homotopy}
The following equality holds:
\begin{equation}\label{eq: gluing chain homotopy}
\widecheck{K}\widecheck{\partial}+\widecheck{\partial}\widecheck{K}=\widecheck{m}_{02,\mathbb{S}}-\widecheck{m}_{01}\widecheck{m}_{12}.
\end{equation}
\end{pro}
\begin{proof}
The proof is essentially identical to \cite[Lemma 26.2.2]{KronheimerMrowkaMonopole}. We consider parametrized moduli space $\mathcal{M}^{+}_{z}([\mathfrak{a}],[\mathfrak{b}])$ of dimension 1. By counting points in various boundary strata of this moduli space, we obtain relations between the operators $\overline{\partial}^{*}_{*}, \partial^{*}_{*}, m^{*}_{*}, m^{*}_{*}, K^{*}_{*}$ (see \cite[Lemma 26.2.3]{KronheimerMrowkaMonopole}), which implies (\ref{eq: gluing chain homotopy}). Note that we are using the family cobordism $\widetilde{W}_{12}$ instead of a single cobordism $W_{12}$. This could potentially change the sign in these relations. However, since we are only considering $1$-dimensional moduli space  $\mathcal{M}^{+}_{z}([\mathfrak{a}],[\mathfrak{b}])$ here. The boundary strata are all $0$-dimensional. As a result, the second factors in Proposition \ref{pro: orientation boundary strata} is always 1. This allows us to use the same sign as  \cite{KronheimerMrowkaMonopole} and repeat their calculation without worrying about sign-issue.
\end{proof}

Proposition \ref{pro: gluing chain homotopy} directly implies that
$$
\sum_{\mathfrak{s}\in \mathbb{S}}\HMto_{\bullet}(\widetilde{W}_{02},\mathfrak{s})= \HMto_{\bullet}(W_{01},\mathfrak{s}_{W_{01}})\circ \HMto_{\bullet}(\widetilde{W}_{12},\mathfrak{s}_{\widetilde{W}_{12}}).
$$
Since there is a natural inclusion $\HMto_{*}(Y_{i},\mathfrak{s}_{Y_{i}})\hookrightarrow \HMto_{\bullet}(Y_{i},\mathfrak{s}_{Y_{i}})$, the proof of Theorem \ref{thm: composition law for HM} (1) is finished. The other parts of Theorem \ref{thm: composition law for HM} are proved by modifying the corresponding proofs of \cite{KronheimerMrowkaMonopole} in exactly the same way. 

\subsubsection{Composition law for the mixed map}
Next, we discuss the composition law for  $\overrightarrow{HM}$. We start with our assumption in the current subsection:
\begin{assum}\label{assum: W12 product}
$(\widetilde{W}_{12},\mathfrak{s}_{\widetilde{W}_{12}})$ is a product family associated to a  $\spinc$-cobordism  $(W_{12},\mathfrak{s}_{W_{12}})$ with $b^{+}_{2}(W_{12})>0$. 
\end{assum}
We fix an orientation on $H^{+}_{2}(W_{12};\mathbb{R})$ so that when $b^{+}_{2}(W_{12})=1$, we can talk about the maps
$$
\HMarrow_{*}(W_{12},\mathfrak{s}_{W_{12}})_{+},\  \HMarrow_{*}(W_{12},\mathfrak{s}_{W_{12}})_{-}: \HMfrom(Y_{1},\mathfrak{s}_{Y_{1}})\rightarrow \HMto(Y_{2},\mathfrak{s}_{Y_{2}})
$$
associated to the positive/negative chamber. If we further have $b^{+}_{2}(W_{02})=b^{+}_{2}(W_{12})=1$, then this orientation on  $H^{+}_{2}(W_{12};\mathbb{R})$ also induces an orientation on $H^{+}_{2}(W_{02};\mathbb{R})$. So we have maps 
$$\HMarrow_{*}(\widetilde{W}_{02},\mathfrak{s})_{+},\ \HMarrow_{*}(\widetilde{W}_{02},\mathfrak{s})_{-}: \HMfrom(Y_{0},\mathfrak{s}_{Y_{0}})\rightarrow \HMto(Y_{2},\mathfrak{s}_{Y_{2}}) .$$ 

When $b^{+}_{2}(W_{12})>1$ (resp. $b^{+}_{2}(W_{02})>1$), we have the map $\HMarrow_{*}(W_{12},\mathfrak{s})$ (resp. $\HMarrow_{*}(\widetilde{W}_{02},\mathfrak{s})$) defined using the unique canonical chamber.

\begin{thm}\label{thm: gluing for hm-arrow}
Under Assumption \ref{assum: W12 product}, we have the following results:
\begin{enumerate}
    \item Suppose $b^{+}_{2}(W_{12})>1$. Then 
$$
\sum_{\mathfrak{s}\in \mathbb{S}}\HMarrow_{*}(\widetilde{W}_{02},\mathfrak{s})=\HMfrom_{*}(\widetilde{W}_{01},\mathfrak{s}_{\widetilde{W}_{01}})\circ \HMarrow_{*}(W_{12},\mathfrak{s}_{W_{12}}).
$$
    \item  Suppose $b^{+}_{2}(W_{02})> b^{+}_{2}(W_{12})=1$. Then 
$$
\sum_{\mathfrak{s}\in \mathbb{S}}\HMarrow_{*}(\widetilde{W}_{02},\mathfrak{s})=\HMfrom_{*}(\widetilde{W}_{01},\mathfrak{s}_{\widetilde{W}_{01}})\circ \HMarrow_{*}(W_{12},\mathfrak{s}_{W_{12}})_{\pm}.
$$
\item Suppose $b^{+}_{2}(W_{02})=b^{+}_{2}(W_{12})=1$. Then 
$$
\sum_{\mathfrak{s}\in \mathbb{S}}\HMarrow_{*}(\widetilde{W}_{02},\mathfrak{s})_{\pm}=\HMfrom_{*}(\widetilde{W}_{01},\mathfrak{s}_{\widetilde{W}_{01}})\circ \HMarrow_{*}(W_{12},\mathfrak{s}_{W_{12}})_{\pm}.
$$
\end{enumerate}
\end{thm}
We will prove Theorem \ref{thm: gluing for hm-arrow} through the rest of this subsection. We will focus on (1) because the proof of the other two cases are essential the same by working with the positive or the negative canonical chamber. 

We first show that the sums in above formulas are always finite.
\begin{lem}
For any chamber $\xi\in [\Q,V^{+}(W_{02})]$, there are only finitely many $\mathfrak{s}\in \mathbb{S}$ such that the map $\HMarrow(\widetilde{W}_{02},\mathfrak{s})_{\xi}$ is nontrivial.
\end{lem}

\begin{proof}
The proof is similar to Lemma \ref{lem: finiteness for HM-to from and bar}. The current case is actually simpler because there is no reducible solutions on $\widetilde{W}_{02}$ by our choice of perturbation 2-forms $\{\omega_{02,b}\}$. Hence no solutions on  $\widetilde{W}_{02}$ can limit to a boundary stable solution on $Y_0$ or a boundary unstable solution on $Y_2$.
\end{proof}
For the  rest of this section, we assume $b^{+}_{2}(W_{12})>1$. We fix a single metric $g_{W_{12}}$ on $W_{12}$ and a family metric  $\{g_{W_{01,b}}\}$ on $\widetilde{W}_{01}$. Note that we have chosen $g_{W_{12}}$ such that the collar neighborhood $N_{i}$ of the boundary component $Y_{i}$ is isometric to $[0,1]\times Y_{i}$. For $T>0$, we let $W^{T}_{12}$ be obtained from $W_{12}$ by inserting two cylinders $[1,T]\times Y_{1}$ and $[1,T]\times Y_{2}$ between $W_{12}\setminus (N_{1}\cup N_{2})$ and $N_{1}\cup N_{2}$. 

For $S>0$, we consider a new family of cobordism $\widetilde{W}(S,T)$ from $Y_{0}$ to $Y_{2}$, whose fiber over $b$ is given by
$$
W(S,T)_{b}:=W_{01,b}\cup ([-S,S]\times Y_{1})\cup W^{T}_{12}.
$$
As before, we use $\widetilde{W}(S,T)^{*}$ to denote the family obtained by attaching cylindrical ends to each fiber.  
We denote the family metric on $\widetilde{W}(S,T)$ by $g_{\widetilde{W}(S,T)}=\{g_{W(S,T)_{b}}\}$.

To construct suitable perturbations on $\widetilde{W}(S,T)^{*}$, we slightly modify Construction \ref{cons: gluing perturbation} as follows: 
\begin{itemize}
    \item In step (1), we pick a self-dual 2-form  $w_{12}$ supported on $W_{12}\setminus (N_{1}\cup N_{2})$. We require that $(\omega_{12},g_{W_{12}})$ is admissible with respect to $\mathfrak{s}_{W_{12}}$ if  $c_{1}(\mathfrak{s}_{Y_{1}})$ and $c_{1}(\mathfrak{s}_{Y_{2}})$ are torsion. Then we use $\omega_{12}$ as a perturbation 2-form on the fiber $W(S,T)_{b}$ for any $S,T,b$.
\item We repeat the rest steps without change. We emphasise that we don't apply \textbf{any} perturbations on the newly-added cylinders $[1,T]\times Y_{1}$ and $[1,T]\times Y_{2}$.
\end{itemize}
\begin{lem}\label{lem: admissible for large T}
Suppose $c_{1}(\mathfrak{s}_{Y_{0}})$ and $c_{1}(\mathfrak{s}_{Y_{2}})$ are both torsion. Then for any $\mathfrak{s}\in \mathbb{S}$, there exists $T_{\mathfrak{s},0}\gg 0$ such that for any  $S\geq 0$ and $T\geq T_{\mathfrak{s},0}$, the pair $(g_{\widetilde{W}(S,T)},\omega_{12})$ is admissible with respect to $\mathfrak{s}$ and the corresponding chamber is the canonical chamber $\xi_{c}$.
\end{lem}
\begin{proof}
Because $(\widetilde{W}(S,T),g_{\widetilde{W}(S,T)})$ is isometric to $(\widetilde{W}(S+T),g_{\widetilde{W}(S+T)})$, we may assume $S=0$. We pick a family $\spinc$ connection $\{A_{b,T}\}$ on the fibers of $\widetilde{W}(T)^{*}$ that satisfies the following conditions:
\begin{itemize}
    \item When restricted to $W_{01,b}$, $A_{b,T}$ is independent with $T$;
    \item When restricted to $W_{12,b}$, $A_{b,T}$ is independent with $T$ and $b$;
    \item When restricted to the cylinders $(-\infty,0]\times Y_{0},\ [-T,T]\times Y_{1}$ and $[0,+\infty)\times Y_{1}$, $A_{b,T}$ is pulled back from a connection $A_{Y_{i}}$ on $Y_{i}$ that is independent with $b,T$. We assume that $A^{t}_{Y_{i}}$ is flat for $i=0,2$. This is possible because our assumption that $c_{1}(\mathfrak{s}_{Y_{1}})$ and $c_{1}(\mathfrak{s}_{Y_{2}})$ are torsion.
\end{itemize}
We first suppose $c_{1}(\mathfrak{s}_{Y_1})$ is torsion. Then we can pick $A_{Y_{1}}$ such that $A^{t}_{Y_{1}}$ is also flat. Hence $F_{A^{t}_{b,T}}$ is supported in $W_{01,b}\cup W_{12,b}$. 

Because  $(\omega_{12},g_{W_{12}})$ is admissible, there exists an embedded surface $\Sigma$ in $W_{12}$ such that $\Sigma\cdot \Sigma\geq 0$ and that
$$
\int_{W_{12,b}} (\ii F^{+}_{A^{t}_{b,T}} +\omega_{12})\wedge \kappa_{\Sigma}>0.
$$
Here $\kappa_{\Sigma}$ is the square integrable, harmonic $2$-form on $W^{*}_{12}$ that corresponds to the image of $[\Sigma]$ under the map
$$
H_{2}(W_{12};\mathbb{R})\xrightarrow{\operatorname{PD}} H^{2}(W_{12},\partial W_{12};\mathbb{R})\rightarrow H^{2}(W_{12};\mathbb{R}).$$
Similarly, we can treat $\Sigma$ as a surface in $W(T)_{b}$ and consider the corresponding harmonic two form $\kappa_{\Sigma,b,T}$ on $W(T)^{*}_{b}$. As $T$ goes to infinity, both norms
$$
\|\kappa_{\Sigma,b,T}|_{W_{01,b}}\|_{L^{2}}\text{ and } \|\kappa_{\Sigma,b,T}|_{W_{12,b}}-\kappa_{\Sigma}|_{W_{12,b}}\|_{L^{2}}$$ converge to zero uniformly in $b$. Therefore, when $T$ is large enough, we have
$$
\int_{W(T)^{*}_{b}} (\ii F^{+}_{A^{t}_{b,T}} +\omega_{12})\wedge \kappa_{\Sigma,b,T}=\int_{W_{01,b}} \ii F^{+}_{A^{t}_{b,T}}\wedge \kappa_{\Sigma,b,T} +\int_{W_{12,b}} (\ii F^{+}_{A^{t}_{b,T}} +\omega_{12})\wedge \kappa_{\Sigma,b,T}>0
$$
for any $b\in \Q$. This implies that $(g_{\widetilde{W}(0,T)},\omega_{12})$ is admissible. By Lemma \ref{lem: topology of positive cone}, the corresponding chamber is the canonical chamber.

Next we assume $c_{1}(\mathfrak{s}_{Y_{1}})$ is nontorsion. Then there exists an embedded surface $\Sigma'$ in $Y_{1}$ such that $\langle c_{1}(\mathfrak{s}_{Y_{1}}),[\Sigma']\rangle>0$. We let $\kappa_{\Sigma',b,T}$ be the square integrable harmonic $2$-form on $W(T)^{*}_{b}$ that corresponds to $\operatorname{PD}([\Sigma'])$. Then as $T\rightarrow \infty$, we have 
$$
\|\kappa_{\Sigma',b,T}|_{W_{01,b}\cup W_{12,b}}\|_{L^2}\rightarrow 0,\quad  \|\kappa_{\Sigma',b,T}|_{[-T,T]\times Y_{1}}-\frac{dt\wedge \alpha}{2T}\|_{L^{2}}\rightarrow 0.
$$
Here $\alpha$ is the harmonic $1$-form on $Y_{1}$ that represents the Poincare dual of $\Sigma'$. Therefore, we have $$
\int_{W(T)^{*}_{b}} (\ii F^{+}_{A^{t}_{b,T}} +\omega_{12})\wedge \kappa_{\Sigma',b,T}\rightarrow \frac{1}{4T}\int _{[-T,T]\times Y_{1} }  \ii (F_{A^{t}_{Y_{1}}}+*F_{A^{t}_{Y_{1}}})\cup (dt\wedge \alpha) =\frac{c_{1}(\mathfrak{s}_{Y_{1}})[\Sigma']}{2}>0.
$$
Rest of proof proceeds similarly.
\end{proof}
We take a pair of critical points 
\begin{equation*}
    \begin{split}
\mathfrak{a}\in \mathfrak{C}^{o}(Y_{0},\mathfrak{s}_{Y_0})\cup \mathfrak{C}^{u}(Y_{0},\mathfrak{s}_{Y_0}),\\   
\mathfrak{b}\in \mathfrak{C}^{o}(Y_{2},\mathfrak{s}_{Y_2})\cup \mathfrak{C}^{s}(Y_{2},\mathfrak{s}_{Y_2}).
    \end{split}
\end{equation*}
and a homotopy class 
$$z\in \pi_{0}(\mathcal{B}^{\sigma}([\mathfrak{a}],\widetilde{W}^{*}_{02},\mathfrak{s},[\mathfrak{b}]).$$ Then we consider the parametrized moduli space 
\begin{equation}
\mathcal{M}_{T,z}([\mathfrak{a}],[\mathfrak{b}])=\bigcup_{S\in [0,+\infty)}\{S\}\times \mathcal{M}_{z}([\mathfrak{a}],\widetilde{W}(S,T)^{*},[\mathfrak{b}])
\end{equation}
\begin{lem}
For any $\mathfrak{s}\in \mathbb{S}$,
there exists $T_{\mathfrak{s},1}\gg 0$ such that for all $T\geq T_{\mathfrak{s},1}$, the moduli space $\mathcal{M}_{T,z}([\mathfrak{a}],[\mathfrak{b}])$ has no reducible for any $\mathfrak{a},\mathfrak{b}$ and $z\in \pi_{0}(\mathcal{B}^{\sigma}([\mathfrak{a}],\widetilde{W}^{*}_{02},\mathfrak{s},[\mathfrak{b}])$.
\end{lem}
\begin{proof}
This lemma is a family version of \cite[Lemma 27.3.5]{KronheimerMrowkaMonopole} and it follows directly from Lemma \ref{lem: admissible for large T}. \end{proof}
We fix a choice of $T_{\mathfrak{s}}\geq \max(T_{0,\mathfrak{s}},T_{1,\mathfrak{s}})$. Using these new moduli spaces $\mathcal{M}_{T,z}([\mathfrak{a}],[\mathfrak{b}])$, we define 
operators $K^{o}_{o},\ K^{o}_{s},\ K^{u}_{s},\ K^{u}_{o}$ exactly as in in (\ref{eq: operators K}). We then define the operator 
\begin{equation*}
    \overrightarrow{K}_{\mathfrak{s}}=\left[\begin{array} {cc}
 K^{o}_{o}  & K^{o}_{s}  \\
 K^{u}_{o}  & K^{u}_{s}
\end{array}\right]:\quad \widehat{C}_{*}(Y_{0},\mathfrak{s}_{Y_{0}})\rightarrow \widecheck{C}_{*}(Y_{2},\mathfrak{s}_{Y_{2}}).
\end{equation*}
There is no need for completion in this case because the rational gradings of the generators of $\widecheck{C}_{*}(Y_{2},\mathfrak{s}_{Y_{2}})$ are bounded below.

Using the perturbations over $S=\infty$, we define the map
$$
\widehat{m}_{01}:\widehat{C}_{*}(Y_{0},\mathfrak{s}_{Y_{0}})\rightarrow \widehat{C}_{*}(Y_{1},\mathfrak{s}_{Y_{1}})
$$
for the family cobordisms $(\widetilde{W}_{01},\mathfrak{s}_{\widetilde{W}_{01}})$, and the map 
$$\overrightarrow{m}_{12}:\widehat{C}_{*}(Y_{1},\mathfrak{s}_{Y_{1}})\rightarrow \widecheck{C}_{*}(Y_{2},\mathfrak{s}_{Y_{2}})$$
for the single cobordism
$(W_{12},\mathfrak{s}_{W_{02}})$. Using the perturbations over $S=0$, we define the map
$$\overrightarrow{m}_{02,\mathfrak{s}}:\widehat{C}_{*}(Y_{0},\mathfrak{s}_{Y_{0}})\rightarrow \widecheck{C}_{*}(Y_{2},\mathfrak{s}_{Y_{2}})$$
for the family cobordism $(\widetilde{W}_{02},\mathfrak{s})$ with any $\mathfrak{s}\in \mathbb{S}$. We set
$$
\overrightarrow{m}_{02,\mathbb{S}}=\sum_{\mathfrak{s}\in \mathbb{S}}\overrightarrow{m}_{02,\mathfrak{s}}\quad \text{and}\quad \overrightarrow{K}_{\mathbb{S}}=\sum_{\mathfrak{s}\in \mathbb{S}}\overrightarrow{K}_{\mathfrak{s}}.
$$

The following proposition can be proved by modifying the  proof of \cite[formula (27.7)]{KronheimerMrowkaMonopole} 
 in a similar way as last subsection. 
\begin{pro}\label{pro: gluing chain homotopy hmarrow} We have $-\overrightarrow{K}_{\mathbb{S}}\circ \widehat{\partial}-\widecheck{\partial}\circ \overrightarrow{K}_{\mathbb{S}}+\overrightarrow{m}_{02,\mathbb{S}}+\overrightarrow{m}_{12}\circ\widehat{m}_{01}=0$.
\end{pro}
Suppose either $\mathfrak{s}_{Y_{0}}$ or $\mathfrak{s}_{Y_{2}}$ is non-torsion. Then there is no chamber structure in defining $\HMarrow_{*}(\widetilde{W}_{02},\mathfrak{s})
$. So Proposition \ref{pro: gluing chain homotopy hmarrow} directly implies Theorem \ref{thm: gluing for hm-arrow} (1).

Suppose $\mathfrak{s}_{Y_{0}}$ and $\mathfrak{s}_{Y_{2}}$ are both torsion. By Lemma \ref{lem: admissible for large T}, the operator $\overrightarrow{m}_{02,\mathfrak{s}}$ corresponds to the canonical chamber $\xi_{c}$ and thus induces the map
$
\HMarrow_{*}(\widetilde{W}_{02},\mathfrak{s})
$. So Proposition \ref{pro: gluing chain homotopy hmarrow} again implies Theorem \ref{thm: gluing for hm-arrow} (1).

\subsection{A computation of some cobordism-induced maps} Computing the cobordism induced maps $\widehat{HM}_{*}(\widetilde{W},\mathfrak{s}_{\widetilde{W}})$ and $\widecheck{HM}_{*}(\widetilde{W},\mathfrak{s}_{\widetilde{W}})$ are difficult in general. But  $\HMbar_{*}(\widetilde{W},\mathfrak{s}_{\widetilde{W}})$ is much easier to handle because it only involves counting reducible solutions. For the same reason, $\HMbar_{*}(\widetilde{W},\mathfrak{s}_{\widetilde{W}})$ contains significant less amount of information about the cobordism $\widetilde{W}$. However, when $Y_{0},Y_{1}$ are $L$-spaces, $\HMbar_{*}(\widetilde{W},\mathfrak{s}_{\widetilde{W}})$ actually determines $\widehat{HM}_{*}(\widetilde{W},\mathfrak{s}_{\widetilde{W}})$ and $\widecheck{HM}_{*}(\widetilde{W},\mathfrak{s}_{\widetilde{W}})$ through the exact sequence (\ref{exact triangle}). The main result of this section (Theorem \ref{thm: computing cobordism induced map}) computes the map $\widehat{HM}_{*}(\widetilde{W},\mathfrak{s}_{\widetilde{W}})$ in this situation. 

We start by stating the assumptions imposed in the current section: 
\begin{assum}\label{assum: homology product family of cobordism} Let $Y_{0}$ and $Y_{1}$ be two rational homology 3-spheres and let $\widetilde{W}/\Q$ be a family cobordism from $Y_0$ to $Y_1$. We assume that $\widetilde{W}/\Q$ is a homology product and the fiber $W$ satisfies $
b_{1}(W)=0$. Let $\mathfrak{s}_{\widetilde{W}}$ be a family $\spinc$ structure on $\widetilde{W}$ which is a product near the boundary.
\end{assum}
We consider the cobordism induced map
$$
\HMbar_{*}(\widetilde{W},\mathfrak{s}_{\widetilde{W}}):\HMbar_{*}(Y_{0},\mathfrak{s}_{Y_{0}})\rightarrow \HMbar_{*}(Y_{1},\mathfrak{s}_{Y_{1}}).
$$
We use $\widebar{e}^{i}_{k}$ to denote the canonical generator of $\HMbar_{*}(Y_{i},\mathfrak{s}_{Y_{i}})$ with rational grading $k$. It  is given by a reducible critical point $[(B_{i},0,\psi^{i}_{*})]$. (See Notation \ref{notation: generators}.) 

In what follows, we will use $\Ind\{P_{*}\}\in K^{0}(\Q)$ to denote the virtual index bundle of a family of complex Fredholm operator $\{P_{b}\}_{b\in \Q}$ and use $\ind \in \mathbb{Z}$ to denote the complex numerical index of a single operator $P$.

\begin{pro}\label{pro: induced map on HM-bar} Under Assumption \ref{assum: homology product family of cobordism}, we have the following results.
\begin{enumerate}
    \item Suppose $b^{+}_{2}(W)>0$. Then $\HMbar_{*}(\widetilde{W},\mathfrak{s}_{\widetilde{W}})$ is trivial.
    \item Suppose $\dimension\Q$ is odd. Then $\HMbar_{*}(\widetilde{W},\mathfrak{s}_{\widetilde{W}})$ is trivial.
    \item Suppose $b^{+}_{2}(W)=0$ and $\dimension\Q$ is even. Then we have 
$$
\HMbar_{*}(\widetilde{W},\mathfrak{s}_{\widetilde{W}})(\widebar{e}^{0}_{n})=c_{\frac{\dimension\Q}{2}}(-\Ind\{\slashed{D}_{A^{0}_{b},\hat{\mathfrak{q}}_{b}}^{+}(W_b)\})[\Q]\cdot \widebar{e}^{1}_{n+d}, 
$$
where 
\begin{equation}
d=\frac{c_{1}(\mathfrak{s}_{W})^2-\sigma(W)}{4}+\dimension\Q,
\end{equation}
and $\{\slashed{D}_{A_{b},\hat{\mathfrak{q}}_{b}}^{+}(W_b)\}    
$ denotes the family perturbed Dirac operator 
\begin{equation}\label{eq: perturbed Dirac unweighted}
    \slashed{D}^{+}_{A^{0}_{b},\hat{\mathfrak{q}}_{b}}:L^{2}_{k}(W_{b}^{*},S^{+})\rightarrow L^2_{k-1}(W_{b}^{*},S^{-}) \quad \text{for}\quad b\in \Q,
\end{equation}
as defined in (\ref{eq: perturbed 4d Dirac equation}). Here we choose any smooth family $\{A^{0}_{b}\}_{b\in \mathcal{Q}}$ of $\spinc$ connections on $\widetilde{W}^{*}$ that equals the pull back of $B_{0}$ and $B_{1}$ on the cylindrical ends. 
\end{enumerate}
\end{pro}

To prove Proposition \ref{pro: induced map on HM-bar}, we start with some general discussions in index theory. Consider a smooth family of Fredholm operator $\{P_{b}:H_1\rightarrow H_2\}_{b\in \mathcal{Q}}$ between Hilbert spaces. They together give a map 
$$\widetilde{P}:\Q\times H_{1}\rightarrow H_{2}.$$ Suppose we have $$2\ind(P_{b})=2-\dimension\Q$$ and suppose $\vec{0}$ is a regular value of $\widetilde{P}$. Then for any 
$b\in\Q$, the composition  
$$
L_{b}:T_{b}\Q\xrightarrow{\mathcal{D}_{(b,0)}} H_{2}\xrightarrow{\text{projection}}\operatorname{coker}P_{b}
$$
is surjective. This implies the inequality 
$$
\dimension \Q\geq \dimension \operatorname{coker}P_{b}= \dimension \operatorname{ker}P_{b}-2\ind(P_{b}).
$$
Therefore, the complex dimension of $\ker P_{b}$ is either $0$ or $1$. When $\ker P_{b}$ is 1-dimensional, the map $L_{b}$ is an isomorphism. For such $b$, we assign it with a positive sign if $L_{b}$ sends the orientation on $T_{b}\Q$ to the complex orientation of $\operatorname{coker}P_{b}$ and a negative sign otherwise.

\begin{lem}\label{lem: kernel counting}
 The number of $b\in Q$ such that $\ker P_{b}\neq 0$, when counted with sign specified above, equals $c_{\frac{\dimension\Q}{2}}(-\Ind\{P_{*}\})[\Q]$.
\end{lem}
\begin{proof}
By compactness of $\Q$, we can find a finite dimensional subbundle $F$ of $H_{2}$ such that $P_{b}$ is  transverse to $F$ for every $b$. Then $E:=\widetilde{P}^{-1}(F)$ is a finite dimensional sub-bundle of $\Q\times H_{1}$. By enlarging $F$ if necessary, we may assume that $E$ is a trivial bundle $\underline{\mathbb{C}}^{N}$. The restriction $\widetilde{P}|_{E}$ defines sections $s_{1},s_{2},\cdots ,s_{N}$ of $F$. Note that $\ker P_{b}\neq 0$ if and only if $\ker (\widetilde{P}|_{E})_{b}\neq 0$ and this happens exactly when $s_{1},\cdots ,s_{N}$ are linearly dependent at $b$. By the classical obstruction theoretic definition of Chern class, the number of such $b$ equals $c_{\dimension(F)+1-N}(F)[\Q]$. The proof is finished by observing that 
$$
[-\Ind\{P_{*}\}]=[F-\underline{\mathbb{C}}^{N}]\in K^{0}(\Q).
$$
\end{proof}
\begin{proof}[Proof of Proposition \ref{pro: induced map on HM-bar}] (1) Suppose $b^{+}_{2}(W)>0$. Then there exists a family of self-dual 2-forms $\{\omega_{b}\}_{b\in \Q}$ supported in the interior of $W_{b}$ such that $\{g_{\widetilde{W}/\Q},\{\omega_{b}\}\}$ is admissible with respect to $\mathfrak{s}_{\widetilde{W}}$. By Lemma \ref{lem: no reducible}, the perturbed Seiberg-Witten equations has no reducible solutions. Therefore, the map    $\HMbar_{*}(\widetilde{W},\mathfrak{s}_{\widetilde{W}})$ is trivial.

(2) Suppose $\Q$ is odd dimensional. By (1), we may also assume $b^{+}_{2}(W)=0$. Then for any pair of reducible critical points $[\mathfrak{a}]$ and $[\mathfrak{b}]$ on $Y_{0}$ and $Y_{1}$. Dimension of the moduli space 
$
\mathcal{M}^{\text{red}}([\mathfrak{a}],\widetilde{W}^{*},[\mathfrak{b}])
$
equals the sum of $\dimension\Q$ and the index of the linearized Seiberg-Witten operators. Since $b_{1}(W)=b^{+}_{2}(W)=0$, the curvature equation doesn't contribute to the index and the index is an even number because the Dirac operator is complex linear. Therefore, when $\Q$ is odd-dimensional, this moduli space can never be $0$-dimensional. Hence the the map    $\HMbar_{*}(\widetilde{W},\mathfrak{s}_{\widetilde{W}})$ is trivial.   

(3) Now suppose $\Q$ is even-dimensional and $b^{+}_{2}(W)=0$. Let $
[(B_{0},0,\psi^{0}_{j_{0}})]$ and $[(B_{1},0,\psi^{1}_{j_{1}})]$
be the critical points that give $\widebar{e}^{0}_{n}$ and $\widebar{e}^{1}_{n+d}$
respectively. Since the map $\HMbar(\widetilde{W},\mathfrak{s}_{\widetilde{W}})$ increase the rational grading by $d$ (see \cite[(28.3)]{KronheimerMrowkaMonopole}), the moduli space 
$$
\mathcal{M}^{\text{red}}:=\mathcal{M}^{\text{red}}([(B_{0},0,\psi^{0}_{j_{0}})],\widetilde{W}^{*},\mathfrak{s}_{\widetilde{W}},[(B_{1},0,\psi^{1}_{j_{1}})])
$$
is an oriented 0-dimensional manifold and we have 
$$
\HMbar_{*}(\widetilde{W},\mathfrak{s}_{\widetilde{W}})(\widebar{e}^{0}_{n})=\#\mathcal{M}^{\text{red}} \cdot \widebar{e}^{1}_{n+d}.
$$
Because $b_{1}(W)=b^{+}_{2}(W)=0$, on each fiber $W^*_{b}$, the curvature equation $F^{+}_{A^{t}}=i\omega_{b}$ has a unique solution that converges to $B_{0}$ and $B_{1}$ on the two ends. Therefore, by choosing the perturbations $\mathfrak{q}_{i}$ and $\mathfrak{q}'_{i}$ small enough, we may assume that for any $b\in \Q$, the the perturbed Seiberg-Witten equations on $W^{*}_{b}$ has a unique finite energy downstairs reducible solution $[(A_{b},0)]$ up to gauge transformation. Upstairs, a reducible solution of the blown-up Seiberg-Witten equations is given by $[(A_{b},\phi,0)]$, where $\phi$ is a nonzero local-$L^{2}_{k}$ spinor that satisfies the perturbed Dirac equation
$
\slashed{D}^{+}_{A_{b},\hat{\mathfrak{q}}_{b}}(\phi)=0
$. (See formula (\ref{eq: perturbed 4d Dirac equation}).)

To discuss the when $[(A_{b},\phi,0)]$ limits to $[(B_{i},0,\psi^{i}_{j_{i}})]$ on the two ends, we consider the weighted Sobolev norm for a spinor $\psi$ on $W_{b}^{*}$:
$$
\|\psi\|_{L^{2}_{p,\lambda,\lambda'}}:=\sqrt{\sum^{p}_{i=0}\|\nabla ^{(p)}_{A_{b}}\tau \cdot \psi\|^{2}_{L^{2}}}.
$$
Here $\lambda,\lambda'\in \mathbb{R}$ are two real numbers. And $\tau_{\lambda,\lambda'}:W^{*}_{b}\rightarrow \mathbb{R}^{>0}$ is a smooth function that satisfies
$$
\tau_{\lambda,\lambda'}(x,t)=\begin{cases} e^{-\lambda t}& \mbox{if } (x,t)\in (-\infty,0]\times Y_{0}\\ e^{\lambda't} & \mbox{if } (x,t)\in [0,\infty)\times Y_{1}.
\end{cases}
$$

On cylindrical ends, the equation $\slashed{D}^{+}_{A_{b},\hat{\mathfrak{q}}_{b}}(\phi)=0$ can be written as 
$$
\frac{d}{dt}\varphi(t)+\slashed{D}_{B_{i}(t),\mathfrak{q}_{i}}\varphi(t)=0,
$$
where $B_{i}(t)=A_{b}\mid _{\{t\}\times Y_{i}}$ and $\varphi(t)=\phi|_{\{t\}\times Y_{i}}$. Using the spectral decomposition for $\slashed{D}_{B_{i}(t),\mathfrak{q}_{i}}$, one can prove that upstairs, $[(A_{b},\phi,0)]$ converges to $[(B_{0},0,\psi^{0}_{j})]$ and $[(B_{1},0,\psi^{1}_{j'})]$ for some $j\leq j_{0}$ and $j'\geq j_{1}$ if and only if
$$
\phi\in L^{2}_{k,-\lambda_{j_0}-\epsilon,\lambda_{j_{1}}-\epsilon}(W^{*}_{b}).
$$
Here $\epsilon$ can be any positive number small enough.
(See \cite[Proposition 14.6.1]{KronheimerMrowkaMonopole} for a similar argument on the cylinder $\mathbb{R}\times Y$.)  However, for any $j\leq j_{0}$ and  any $j'\geq j_{1}$, the moduli space 
$$
\mathcal{M}^{\text{red}}([(B_{0},0,\psi^{0}_{j})],\widetilde{W}^{*},[(B_{1},0,\psi^{1}_{j'})])
$$
has a negative expected dimension unless both equality are achieved. By our transversality assumption, these moduli spaces are all empty. So $[(A_{b},0,\phi)]$ belongs to $\mathcal{M}^{\text{red}}$ if and only if $\phi$ is a nonzero vector in the kernel of the operator
\begin{equation}\label{eq: perturbed Dirac weighted}
\slashed{D}^{+}_{A_{b},\hat{\mathfrak{q}}_{b}}:L^{2}_{k,-\lambda_{j_0}-\epsilon,\lambda_{j_{1}}-\epsilon}(W^{*}_{b};S^{+})\rightarrow L^{2}_{k-1,-\lambda_{j_0}-\epsilon,\lambda_{j_{1}}-\epsilon}(W^{*}_{b};S^{-}) \quad \text{for}\quad b\in \Q.
\end{equation}
This leads to our explicit description of the moduli space 
$$
\mathcal{M}^{\text{red}}\cong \bigcup_{b\in \Q}\mathbb{P}(\ker \slashed{D}^{+}_{A_{b},\hat{\mathfrak{q}}_{b}})
$$
where $\mathbb{P}(\cdot)$ denotes the associated complex project space. (Note that $(A_{b},0,\phi)$ and $(A_{b},0,\phi')$ are gauge equivalent if $\phi=z\phi'$ for some $z\neq 0$.)

At a point $[(A_{b},0,\phi)]\in \mathcal{M}^{\text{red}}$, deformation of the moduli space is controlled by the Fredholm operator 
\begin{equation*}
\begin{split}
R:T_{b}\Q\oplus L^{2}_{k,-\lambda_{j_0}-\epsilon,\lambda_{j_{1}}-\epsilon}(W^{*}_{b};S^{+})\rightarrow L^{2}_{k-1,-\lambda_{j_0}-\epsilon,\lambda_{j_{1}}-\epsilon}(W^{*}_{b};S^{+})\oplus \mathbb{C}\\
(\overrightarrow{v},\psi)\mapsto (\mathcal{D}_{(b,\phi)}\slashed{D}^{+}_{A_{b},\hat{\mathfrak{q}}_{b}}(\overrightarrow{v},\psi),\int_{W^{*}_{b}} \tau ^2\langle \psi,\phi\rangle d\text{vol} )
\end{split}    
\end{equation*}
Since $\mathcal{M}^{\text{red}}$ is regular and $0$-dimensional, we see that $R$ is surjective and has index $0$. This implies that the family of Dirac operators $\{\slashed{D}^{+}_{A_{b},\hat{\mathfrak{q}}_{b}}\}_{b\in \Q}$ in (\ref{eq: perturbed Dirac weighted}) satisfies the condition of Lemma \ref{lem: kernel counting}. As a result, we get 
$$
\# \mathcal{M}^{\text{red}}= c_{\frac{\dimension\Q}{2}}(-\Ind\{\slashed{D}^{+}_{A_{b},\hat{\mathfrak{q}}_{b}}\})[\Q].
$$
The proof is completed by the following lemma.
\end{proof}
\begin{lem}\label{lem: relating family index}
The index bundles of the families (\ref{eq: perturbed Dirac unweighted}) and (\ref{eq: perturbed Dirac weighted}) differ from each other by a trivial bundle over $\Q$. Hence they have the same Chern class.
\end{lem}
\begin{proof}
We first deform the connection $A_{b}$ to the connection $A^{0}_{b}$ via the linear homotopy 
$$
A^{s}_{b}:=s\cdot A_{b}+ (1-s)\cdot A^{0}_{b},\quad s\in [0,1].
$$
Note that for any $s$, the restrictions $A^{s}_{b}|_{t\times Y_{1}}$ and $A^{s}_{b}|_{-t\times Y_{0}}$ limits to $B_{1}$ and $B_{0}$ respectively. Since $\lambda_{j_{i}}\pm \epsilon$ is not an eigenvalue of the Dirac operator $\slashed{D}_{B_{i},\mathfrak{q}_{i}}$. For any $s\in [0,1]$ and $b\in \Q$, the operator
$$
\slashed{D}^{+}_{A^{s}_{b},\hat{\mathfrak{q}}_{b}}:L^{2}_{k,-\lambda_{j_0}-\epsilon,\lambda_{j_{1}}-\epsilon}(W^{*}_{b};S^{+})\rightarrow L^{2}_{k-1,-\lambda_{j_0}-\epsilon,\lambda_{j_{1}}-\epsilon}(W^{*}_{b};S^{-})
$$
is Fredholm. By invariance of the index bundle under continuous deformation, we see that (\ref{eq: perturbed Dirac weighted}) has the same index bundle as the family 
\begin{equation}\label{eq: perturbed Dirac weighted, base connection}
\slashed{D}^{+}_{A^{0}_{b},\hat{\mathfrak{q}}_{b}}:L^{2}_{k,-\lambda_{j_0}-\epsilon,\lambda_{j_{1}}-\epsilon}(W^{*}_{b};S^{+})\rightarrow L^{2}_{k-1,-\lambda_{j_0}-\epsilon,\lambda_{j_{1}}-\epsilon}(W^{*}_{b};S^{-}) \quad \text{for}\quad b\in \Q.
\end{equation}
The only difference (\ref{eq: perturbed Dirac unweighted}) and (\ref{eq: perturbed Dirac weighted, base connection}) is that we are using different weights on the Sobolev spaces. To compare their index bundles, we consider the following four Dirac operators, treated as trivial family over $\Q$: 
\begin{equation}\label{eq: Dirac excision I}
\slashed{D}^{+}_{\hat{B}_{0},\hat{\mathfrak{q}}_{0}}:     L^{2}_{k,-\lambda_{j_{0}}-\epsilon,0}(\mathbb{R}\times Y_{0};S^{+})\rightarrow L^{2}_{k-1,-\lambda_{j_{0}}-\epsilon,0}(\mathbb{R}\times Y_{0};S^{-})
\end{equation}
\begin{equation}\label{eq: Dirac excision II}
\slashed{D}^{+}_{\hat{B}_{0},\hat{\mathfrak{q}}_{0}}:     L^{2}_{k,0,0}(\mathbb{R}\times Y_{0};S^{+})\rightarrow L^{2}_{k-1,0,0}(\mathbb{R}\times Y_{0};S^{-})
\end{equation}
\begin{equation}\label{eq: Dirac excision III}
\slashed{D}^{+}_{\hat{B}_{1},\hat{\mathfrak{q}}_{1}}:     L^{2}_{k,0,\lambda_{j_{1}}-\epsilon}(\mathbb{R}\times Y_{1};S^{+})\rightarrow L^{2}_{k-1,0,\lambda_{j_{1}}-\epsilon}(\mathbb{R}\times Y_{1};S^{-})
\end{equation}
\begin{equation}\label{eq: Dirac excision IV}
\slashed{D}^{+}_{\hat{B}_{1},\hat{\mathfrak{q}}_{1}}:     L^{2}_{k,0,0}(\mathbb{R}\times Y_{1};S^{+})\rightarrow L^{2}_{k-1,0,0}(\mathbb{R}\times Y_{1};S^{-}).
\end{equation}
Here $\hat{B}_{i}$ and $\hat{\mathfrak{q}}_{i}$ are 4-dimensional connection and perturbation pulled back from $B_{i}$ and $\mathfrak{q}_{i}$ respectively. By the excision theorem of index, we see that the direct sum of the index bundles for the families (\ref{eq: perturbed Dirac weighted, base connection}), (\ref{eq: Dirac excision II}) and (\ref{eq: Dirac excision IV}) equals the corresponding bundle for the families (\ref{eq: perturbed Dirac unweighted}), (\ref{eq: Dirac excision I}) and (\ref{eq: Dirac excision III}). This shows that the index bundles of the families  (\ref{eq: perturbed Dirac unweighted}) and (\ref{eq: perturbed Dirac weighted, base connection}) differ from each other by a linear combination of index bundles of (\ref{eq: Dirac excision I}) and (\ref{eq: Dirac excision II}), (\ref{eq: Dirac excision III}) and (\ref{eq: Dirac excision IV}), which are all trivial bundles over $\Q$. This finishes the proof. 
\end{proof}
In practical, it's convenient to convert the index bundle of the family (\ref{eq: perturbed Dirac unweighted}) into a family of Dirac operators on closed 4-manifolds, so that one can apply the Atiyah-Singer index theorem. The following lemma can be proved similarly as Lemma \ref{lem: relating family index} using the excision principle of indices.  
\begin{lem}\label{lem: index excision}
For $i=0,1$, let $(M_{i},\mathfrak{s}_{M_{i}})$ a be $\spinc$ 4-manifold bounded by $(Y_{i},\mathfrak{s}_{Y_i})$. Consider the new family of \emph{closed} $\spinc$ 4-manifolds $(\widetilde{W}',\mathfrak{s}_{\widetilde{W}'})$ obtained by gluing together $\widetilde{W}$, $M_{0}\times \Q$ and $-M_{1}\times \Q$. Then the index bundle of the family Dirac operators 
$$
\slashed{D}_{A'_{b}}^{+}(W'_b,\mathfrak{s}_{W'_{b}}): L^{2}_{k}(W'_b;S^{+})\rightarrow L^{2}_{k-1}(W'_b;S^{-})\quad \text{for }b\in \Q
$$
differ from the index bundle of the family (\ref{eq: perturbed Dirac unweighted}) by a trivial bundle and hence have the same Chern classes. Here we can choose any smooth family of base connection $\{A'_{b}\}$ on $W'_{b}$.
\end{lem}

Now we state and proof the main result of this section. 
\begin{thm}\label{thm: computing cobordism induced map}
Let $(\widetilde{W},\mathfrak{s}_{\widetilde{W}})$ be a $\spinc$-family cobordism from $S^{3}$ to an $L$-space $Y$ which satisfies Assumption (\ref{assum: homology product family of cobordism}). Let $\widehat{1}$ be the standard generator of $\widehat{HM}_{*}(S^{3})$ as a $\mathbb{Q}[U]$-module. Then one has
\begin{equation}\label{eq: cobodism induced map on HMhat}
\widehat{HM}_{*}(\widetilde{W},\mathfrak{s}_{\widetilde{W}})(\widehat{1})=\langle c_{\frac{\dimension\Q}{2}}(-\Ind\{\slashed{D}_{A^{0}_{b},\hat{\mathfrak{q}}_{b}}^{+}(W_b)\}),[\Q]\rangle \cdot \hat{e}_{d-1}
\end{equation}
if $\dimension\Q$ is even and $b^{+}_{2}(\widetilde{W})=0$. 
Otherwise, $\widehat{HM}_{*}(\widetilde{W},\mathfrak{s}_{\widetilde{W}})(\widehat{1})=0$.
\end{thm}
\begin{proof}
We set $Y_{0}=S^{3}$ and $Y_{1}=Y$. Cobordism induced maps are compatible with the exact sequence (\ref{exact triangle}). In particular, we have 
\begin{equation}\label{eq: cobordism induced map commutes with connecting morphism}
p_{*,Y}\circ \widehat{HM}(\widetilde{W},\mathfrak{s}_{\widetilde{W}})=\HMbar(\widetilde{W},\mathfrak{s}_{\widetilde{W}})\circ p_{*,S^{3}}.
\end{equation}
First, we assume that $\dimension\Q$ is odd or $b^{+}_{2}(\widetilde{W})>0$. By Proposition \ref{pro: induced map on HM-bar}, the  map $\HMbar(\widetilde{W},\mathfrak{s}_{\widetilde{W}})$ is trivial. Since $Y$ is an $L$-space, the map $p_{*,Y}$ is injective. So $\widehat{HM}(\widetilde{W},\mathfrak{s}_{\widetilde{W}})$ is trivial as well. 

Next, we assume that $\dimension\Q$ is even, $b^{+}_{2}(\widetilde{W})=0$, and $-2h(Y,\mathfrak{s}_Y)< d-1$. Since  $\gr^{\mathbb{Q}}(\widehat{1})=-1$ and $\widehat{HM}_{*}(\widetilde{W},\mathfrak{s}_{\widetilde{W}})$ increases the rational grading by $d$, we have 
$$
\gr^{\mathbb{Q}}(\widehat{HM}_{*}(\widetilde{W},\mathfrak{s}_{\widetilde{W}})(\widehat{1}))=d-1.
$$
Since $Y$ is an $L$-space, $\widehat{HM}_{*}(Y,\mathfrak{s}_{Y})$ is supported in $\gr^{\mathbb{Q}}$-grading $\leq -2h(Y,\mathfrak{s}_{Y})$. Therefore, $\widehat{HM}(\widetilde{W},\mathfrak{s}_{\widetilde{W}})_{*}(\widehat{1})=0$. On the other hand, by (\ref{eq: cobordism induced map commutes with connecting morphism}) and the fact that $p_{*,S^{3}}$ is injective. We conclude that the map $\HMbar(\widetilde{W},\mathfrak{s}_{\widetilde{W}})$ is trivial as well. By Proposition \ref{pro: induced map on HM-bar}, this implies $c_{\frac{\dimension\Q}{2}}(-\Ind\{\slashed{D}_{A^{0}_{b},\hat{\mathfrak{q}}_{b}}^{+}(W_b)\})[\Q]=0$. So in this case, both sides of (\ref{eq: cobodism induced map on HMhat}) are zero.

Lastly, we assume $\dimension\Q$ is even, $b^{+}_{2}(\widetilde{W})=0$, and $-2h(Y,\mathfrak{s}_Y)\geq d-1$. Then we have $\widehat{HM}(\widetilde{W},\mathfrak{s}_{\widetilde{W}})_{*}(\widehat{1})=x\cdot \widehat{e}_{d-1}(Y)$ for some constant $x$. Note that $
p_{*,S^{3}}(\widehat{1})=\overline{e}^{0}_{-2}$ and $ p_{*,Y}(\widehat{e}_{d-1}(Y))=\overline{e}^{1}_{d-2}
$. By (\ref{eq: cobordism induced map commutes with connecting morphism}), we have 
$$\HMbar(\widetilde{W},\mathfrak{s}_{\widetilde{W}})(\overline{e}^{0}_{-2})=x\cdot \overline{e}^{1}_{d-2}
$$
By Proposition \ref{pro: induced map on HM-bar}, we see that $x=c_{\frac{\dimension\Q}{2}}(-\Ind\{\slashed{D}_{A^{0}_{b},\hat{\mathfrak{q}}_{b}}^{+}(W_b)\})[\Q]$. This proves the last case.
\end{proof}

\section{A Gluing theorem for the family Seiberg-Witten invariant} 
In this section, we prove Theorem \ref{thm: FSW gluing}. This is a gluing formula which expresses the family Seiberg-Witten invariants in terms of Floer theoretic data. Then we prove Theorem \ref{thm: family switching}, which relates the  Seiberg-Witten invariants for certain families with the corresponding invariant for their fibers. This formula will be the main tool in our applications.
 
\subsection{Expressing family Seiberg invariants using cobordism induced maps} Let $\widetilde{M}/\Q$ be a smooth family of closed 4-manifolds which is a homology product. Let $\tau_{i}:\Q\rightarrow \widetilde{M}$ ($i=0,1$) be two disjoint smooth sections. We assume the vertical tangent bundle $T\widetilde{M}/\Q$ is trivial when restricted to the image of $\tau_{i}$ and we fix a trivilization. Let $\widetilde{W}/\Q$ be the family of 4-manifolds with boundary, obtained by removing 
tubular neighborhoods of $\tau_{i}(\Q)$. Then we can identity components $\partial \widetilde{W}$ with the unit sphere bundles of $(T\widetilde{M}/\Q)|_{\tau_{i}(Q)}$. Using the trivializations we chose, we get a natural boundary parametrization of $\widetilde{W}$, making it into a family of cobordism from $S^{3}$ to $S^{3}$. Let $\widetilde{\mathfrak{s}}$ be a family $\spinc$ structure on $\widetilde{M}$. We assume that $\widetilde{\mathfrak{s}}$ is the trivial $\spinc$ structure when restricted to $\tau_{i}(\Q)$, then $\mathfrak{s}_{\widetilde{W}}$ (the restriction of $\widetilde{\mathfrak{s}}$ to $\widetilde{W}$) will satisfy Assumption \ref{assum: homology product family of cobordism}. Hence we can talk about cobordism induced maps. 

\begin{thm}\label{thm: family SW in terms of Floer}
Suppose $d(\widetilde{\mathfrak{s}},\widetilde{M})=0$. Then for any chamber $\xi\in [\Q,V^{+}(W)]$, we have the formula 
\begin{equation}\label{eq: FSW in terms of HMvec}
\FSW_{\xi}(\widetilde{M},\widetilde{\mathfrak{s}})= \langle \HMarrow_{*}(\widetilde{W},\mathfrak{s}_{\widetilde{W}})_{\xi}(\widehat{1}),\widecheck{1}\rangle.    
\end{equation}
Here $\widehat{1}$ and $\widecheck{1}$ denotes the standard generator of  $\widehat{HM}_{*}(S^{3})\cong \mathbb{Q}[U]$ and $\widecheck{HM}^{*}(S^{3})\cong \mathbb{Q}[U]$ respectively, and $\langle-,-\rangle$ denotes the natural pairing (\ref{eq: pairing}).
\end{thm}
\begin{proof}
This is a family version of \cite[Proposition 27.4.1]{KronheimerMrowkaMonopole}. It is essentially a neck-stretching argument. For any $T> 0$, we consider a decomposition of $\widetilde{M}$ into five pieces 
$$
\widetilde{M}_{T}:=(B^{4}\times \Q)\cup (S^{3}\times [-T,0]\times \Q)\cup \widetilde{W}\cup (S^{3}\times [0,T]\times \Q)\cup (-B^{4}\times \Q),
$$
where $-B^{4}$ denotes the orientation reversal of $B^{4}$. We specify the family metrics and perturbations on each piece as follows:
\begin{itemize}
    \item On the cylindrical pieces $S^{3}\times [-T,0]\times b$ and $S^{3}\times [0,T]\times b$, the metrics are the cylindrical metrics for the standard round metric on $S^{3}$. 
    \item On $B^{4}\times \{b\}$ and $-B^{4}\times \{b\}$, we put positive scalar curvature metrics which are cylindrical near the boundary.
    \item We put a smooth family of metrics $g_{\widetilde{W}/\Q}=\{g_{b}\}$ on the fibers of $\widetilde{W}$. They also need to be cylindrical near the boundaries. 
    \item We pick a family of compactly supported 2-forms $\widetilde{w}=\{\omega_{b}\}$ on $W_{b}$ to perturb the curvature equation. We require that $(g_{\widetilde{W}/\Q},\widetilde{w})$ is admissible with respect to $\mathfrak{s}_{\widetilde{W}}$ and gives the chamber $\xi$. Note that $\omega_{b}$ depends on $b$ but doesn't depend on $T$.
    \item We pick a small 3-dimensional perturbation $\mathfrak{q}$ on $S^{3}$ as in Section \ref{sec: definition of HM} and denote the corresponding 4-dimensional perturbation by $\hat{\mathfrak{q}}$. 
    \item On $B^{4}\times \{b\}$ and $-B^{4}\times \{b\}$, we put 4-dimensional perturbations $\hat{\mathfrak{q}}^{\pm}$ which equals $\hat{\mathfrak{q}}$ near the boundary. $\hat{\mathfrak{q}}^{\pm}$ are independent with $b$ and $T$.  
    \item Denote the perturbations on $S^{3}\times [-T,0]\times \{b\}$ and $S^{3}\times [0,T]\times \{b\}$ by $\hat{\mathfrak{q}}^{-}_{T,b}$ and $\hat{\mathfrak{q}}^{+}_{T,b}$ respectively. Then we require that 
    \begin{equation}\label{eq: perturbation near the ball}
    \hat{\mathfrak{q}}^{-}_{T,b}|_{[-T,-\frac{2T}{3}]\times S^{3}}=\hat{\mathfrak{q}},\quad \hat{\mathfrak{q}}^{+}_{T,b}|_{[\frac{2T}{3},T]\times S^{3}}=\hat{\mathfrak{q}}.
\end{equation}
    and 
    $$
    \hat{\mathfrak{q}}^{-}_{T,b}|_{[-\frac{T}{3},0]\times S^{3}}=0,\quad \hat{\mathfrak{q}}^{+}_{T,b}|_{[0,\frac{T}{3}]\times S^{3}}=0.
    $$
    Note that $\hat{\mathfrak{q}}^{\pm}_{T,b}$ will depend on both $T$ and $b$.
\end{itemize}
Let $\mathcal{M}_{T}$ denote the moduli space of perturbed Seiberg-Witten equation on the family $\widetilde{M}_{T}$. Then by choosing generic perturbations as above, one can ensure that the parametrized moduli space 
$$
\mathcal{M}:=\mathop{\bigcup}\limits_{T\in [T_{0},+\infty)} \mathcal{M}_{T}
$$
is regular, where $T_{0}$ is a large enough number. As a result, $\mathcal{M}$ is a 1-dimensional manifold because $d(\widetilde{\mathfrak{s}},\widetilde{M})=0$. To compactify $\mathcal{M}$, we need add in the corresponding moduli space when $T=\infty$. Namely, the broken trajectories on 
$$
((B^{4})^{*}\times Q)\cup \widetilde{W}^{*}\cup ((-B^{4})^{*}\times Q)
$$
Here $(B^{4})^{*}$ is obtained by attaching a cylindrical end to $B^{4}$. This trajectory belongs to the moduli space 
$$
\mathcal{M}^{i,j}_{\infty}=\mathcal{M}((B^{4})^{*}\times Q, [\mathfrak{a}_{i}])\times_{Q}\mathcal{M}([\mathfrak{a}_{i}],\widetilde{W}^{*},[\mathfrak{a}_{j}])\times_{Q}\mathcal{M}([\mathfrak{a}_{j}],(-B^{4})^{*}\times Q).
$$
Note that (\ref{eq: perturbation near the ball}) implies that the limiting perturbation on $
(B^{4})^{*}
$ and $(-B^{4})^{*}$ are independent with $b$, therefore, we have
$$
\mathcal{M}((B^{4})^{*}\times Q, [\mathfrak{a}_{i}])=\mathcal{M}((B^{4})^{*}, [\mathfrak{a}_{i}])\times Q
$$
and 
$$
\mathcal{M}([\mathfrak{a}_{j}],(-B^{4})^{*}\times Q)=\mathcal{M}([\mathfrak{a}_{j}],(-B^{4})^{*})\times Q.
$$
That implies
$$
\mathcal{M}^{i,j}_{\infty}=\mathcal{M}((B^{4})^{*}, [\mathfrak{a}_{i}])\times \mathcal{M}([\mathfrak{a}_{i}],\widetilde{W}^{*},[\mathfrak{a}_{j}])\times \mathcal{M}([\mathfrak{a}_{j}],(-B^{4})^{*}).
$$
By computing the expected dimension, we see that $\mathcal{M}^{i,j}_{\infty}$ is empty unless $i=-1$ and $j=0$. Furthermore, as proved in \cite{KronheimerMrowkaMonopole}, when the perturbations are small enough, $\mathcal{M}((B^{4})^{*}, [\mathfrak{a}_{-1}])$ and $\mathcal{M}([\mathfrak{a}_{0}],(-B^{4})^{*})$ both contain a single point (with positive sign). Hence we further have 
$$
\mathcal{M}^{-1,0}_{\infty}=\mathcal{M}([\mathfrak{a}_{-1}],\widetilde{W}^{*},[\mathfrak{a}_{0}]).
$$
To this end, we see that we can compactify $\mathcal{M}$ into a 1-dimensional manifold $\mathcal{M}^{+}$ with boundary by adding in points in $\mathcal{M}([\mathfrak{a}_{-1}],\widetilde{W}^{*},[\mathfrak{a}_{0}])$. The signed count of boundary points of $\mathcal{M}^{+}$ must be zero. That implies 
$$
\#\mathcal{M}_{T_0}=\# \mathcal{M}([\mathfrak{a}_{-1}],\widetilde{W}^{*},[\mathfrak{a}_{0}]).
$$
This is exactly (\ref{eq: FSW in terms of HMvec}) because $\widecheck{1}$ and $\widehat{1}$ correspond to the critical points $[\mathfrak{a}_{0}]$ and $[\mathfrak{a}_{-1}]$ respectively.
\end{proof}
\subsection{A gluing formula for the family Seiberg-Witten invariants} Now we are ready to  prove Theorem \ref{thm: FSW gluing}. For readers' convenience, we restate our setting as follows:

Let $\widetilde{M}_{1}$ be a family of 4-manifolds with boundary $Y$. We assume $\widetilde{M}_{1}$ is a homology product and a product near the boundary. We take a  family $\spinc$ structure  $\mathfrak{s}_{\widetilde{M_{1}}}$ on $\widetilde{M}_{1}$ which is a product near the boundary and we use $\mathfrak{s}_{Y}$ to denote its resection to the boundary of each fiber. As in previous section, we use the product structure near the boundary to specify a section $s:\Q\rightarrow \widetilde{M}_{1}$. By fiberwisely removing small balls surrounding the image of $s$, we obtain a family $\spinc$ cobordism $(\widetilde{W}_{01},\mathfrak{s}_{\widetilde{W}_{01}})$ from $S^{3}$ to $Y$.

On the other side, we take a single 4-manifold $M_{2}$ with boundary $-Y$. Let $\mathfrak{s}_{M_{2}}$ be a $\spinc$ that is isomorphic to $\mathfrak{s}_{Y})$ on the boundary. We remove a small ball from $M_{2}$ to obtain a $\spinc$ cobordism $(W_{12}, \mathfrak{s}_{W_{12}})$ from $Y$ to $S^{3}$. 

We then form a family of closed 4-manifold 
\begin{equation}\label{eq: closed family for gluing}
\widetilde{M}:=\widetilde{M}_{1}\cup_{\Q\times Y}(\Q\times M_{2})    
\end{equation}
with fiber $M=M_{1}\cup_{Y}M_{2}$

\begin{thm}[Theorem \ref{thm: FSW gluing} restated]
Let $\widetilde{M}$ be defined as above. Then we have the following results regarding the family Seiberg-Witten invariants of $\widetilde{M}$.
\begin{enumerate}
    \item Suppose $b^{+}_{2}(M_{2})>1$. Then we have 
    $$
\sum_{\widetilde{\mathfrak{s}}}\FSW(\widetilde{M},\widetilde{\mathfrak{s}})= \langle\HMfrom_{*}(\widetilde{W}_{01},\mathfrak{s}_{\widetilde{W}_{01}})(\widehat{1})\ ,\  \HMarrow^{*}(W_{12},\mathfrak{s}_{W_{12}})(\widecheck{1})\rangle $$
\item Suppose $b^{+}_{2}(M)>b^{+}_{2}(M_{2})=1$. Then we have 
    $$
\sum_{\widetilde{\mathfrak{s}}}\FSW(\widetilde{M},\widetilde{\mathfrak{s}})= \langle\HMfrom_{*}(\widetilde{W}_{01},\mathfrak{s}_{\widetilde{W}_{01}})(\widehat{1})\ ,\  \HMarrow^{*}(W_{12},\mathfrak{s}_{W_{12}})_{\pm}(\widecheck{1})\rangle .$$
\item Suppose $b^{+}_{2}(M)=b^{+}_{2}(M_{2})=1$. Then we have 
    $$
\sum_{\widetilde{\mathfrak{s}}}\FSW(\widetilde{M},\widetilde{\mathfrak{s}})_{\pm}= \langle\HMfrom_{*}(\widetilde{W}_{01},\mathfrak{s}_{\widetilde{W}_{01}})(\widehat{1})\ ,\  \HMarrow^{*}(W_{12},\mathfrak{s}_{W_{12}})_{\pm}(\widecheck{1})\rangle .$$
\end{enumerate}
In these formulas, we are summing over all family $\spinc$ structures $\widetilde{\mathfrak{s}}$ that satisfy
$$
\widetilde{\mathfrak{s}}|_{\widetilde{M}_{1}}=\mathfrak{s}_{\widetilde{M}_{1}},\ \widetilde{\mathfrak{s}}|_{M_{2}}=\mathfrak{s}_{M_{2}}\text{ and }d(\widetilde{\mathfrak{s}},\widetilde{M})=0.
$$
\end{thm}
\begin{proof}
This is a direct corollary of Theorem \ref{thm: family SW in terms of Floer} and Theorem \ref{thm: gluing for hm-arrow}. Note that when $b_{1}(Y)>0$, there can be a family $\spinc$ structure $\widetilde{\mathfrak{s}}$ such that $\widetilde{\mathfrak{s}}|_{\widetilde{M}_{1}}=\mathfrak{s}_{\widetilde{M}_{1}}$, $\widetilde{\mathfrak{s}}|_{M_{2}}=\mathfrak{s}_{M_{2}}$ but $d(\widetilde{\mathfrak{s}},\widetilde{M})\neq 0$. Such a $\spinc$ structure has no contribution on either side of the above equalities.
\end{proof}
.
\begin{thm}[Theorem \ref{thm: family switching} restated]\label{thm: FSW=SW}
Let $\widetilde{M}$ be a family of closed 4-manifolds with an even dimensional, connected base $\Q$. Suppose $\pi_{1}(\Q)$ acts trivially on $H_{*}(M;\mathbb{Z})$ and suppose $\widetilde{M}$ has a decomposition of the form (\ref{eq: closed family for gluing}) with $Y$ being an $L$-space. Let $\widetilde{\mathfrak{s}}$  be a family $spin^{c}$ structure on $\widetilde{M}$ and let $\mathfrak{s}$ be a $spin^{c}$ structure on the fiber $M$ which satisfy $\widetilde{\mathfrak{s}}|_{M_{2}}=\mathfrak{s}|_{M_{2}}$
and 
$
d(\mathfrak{s},M)\geq d(\widetilde{\mathfrak{s}},\widetilde{M})=0. 
$ Then the following results hold.
\begin{enumerate}
    \item Suppose $b^{+}_{2}(M_{1})=0$ and $b^{+}_{2}(M_{2})>1$. Then
    \begin{equation}\label{eq: FSW=SW 1}
    \FSW(\widetilde{M},\widetilde{\mathfrak{s}})=\langle c_{\frac{\operatorname{dim}(\Q)}{2}}(-\operatorname{Ind}(\slashed{D}^{+}(\widetilde{M},\widetilde{\mathfrak{s}}))\otimes L^{-1}),[\Q]\rangle \cdot \SW(M,\mathfrak{s})
\end{equation}
\item Suppose $b^{+}_{2}(M_{1})=0$ and $b^{+}_{2}(M_{2})=1$. Then \begin{equation}\label{eq: FSW=SW 2}
    \FSW_{\pm}(\widetilde{M},\widetilde{\mathfrak{s}})=\langle c_{\frac{\operatorname{dim}(\Q)}{2}}(-\operatorname{Ind}(\slashed{D}^{+}(\widetilde{M},\widetilde{\mathfrak{s}}))\otimes L^{-1}),[\Q]\rangle \cdot \SW_{\pm}(M,\mathfrak{s})
\end{equation}
\item Suppose $b^{+}_{2}(M_{1})>0$ and $b^{+}_{2}(M)>1$. Then $\FSW(\widetilde{M},\widetilde{\mathfrak{s}})=0.$ 
\item Suppose $b^{+}_{2}(M_{1})>0$ and $b^{+}_{2}(M)=1$. Then  $\FSW_{\pm}(\widetilde{M},\widetilde{\mathfrak{s}})=0$.
\end{enumerate}
Here $L$ is the unique complex line bundle over $\Q$ such that $(\widetilde{\mathfrak{s}}|_{\Q\times M_{2}})\otimes p^{*}(L^{-1})$ 
is a pull back of $(M_{2}, \widetilde{\mathfrak{s}}|_{M_{2}})$. And $\operatorname{Ind}(\slashed{D}^{+}(\widetilde{M},\widetilde{\mathfrak{s}})$ denotes the index bundle of the family Dirac operator.
\end{thm}
\begin{proof}
Suppose $b^{+}_{2}(M_1)>0$. By Theorem \ref{thm: computing cobordism induced map}, we have $\HMfrom_{*}(\widetilde{W}_{01},\mathfrak{s}_{\widetilde{W}_{01}})(\widehat{1})=0$. By Theorem \ref{thm: FSW gluing}, we have $\FSW(\widetilde{M},\widetilde{\mathfrak{s}})=0$ when $b^{+}_{2}(M)>1$ and we have $\FSW_{\pm}(\widetilde{M},\widetilde{\mathfrak{s}})=0$ when $b^{+}_{2}(M)=1$. This proves the case (3) and (4).

Next, we assume $b^{+}_{2}(M_{2})>1$. By Lemma \ref{twisting equality}, replacing $\widetilde{\mathfrak{s}}$ with $(\widetilde{\mathfrak{s}})\otimes p^{*}(L^{-1})$ does not change the family Seiberg-Witten invariant. Therefore, it suffice to show that 
\begin{equation}
    \FSW_{\xi_{c}}(\widetilde{M},\widetilde{\mathfrak{s}})=\langle c_{\frac{\operatorname{dim}(\Q)}{2}}(-\operatorname{Ind}(\slashed{D}^{+}(\widetilde{M},\widetilde{\mathfrak{s}})),[\Q]\rangle \cdot \SW(M,\mathfrak{s})
\end{equation}
under the additional assumption that $\widetilde{\mathfrak{s}}|_{\Q\times M_{2}}$ is pulled back from $\widetilde{\mathfrak{s}}|_{M_{2}}$.
Note that since $Y$ is a rational homology 3-sphere, the family $\spinc$-structure $\widetilde{\mathfrak{s}}$ is uniquely determined by its restrictions to $\widetilde{M}_{1}$ and $\Q\times M_{2}$. Therefore, by applying Theorem \ref{thm: FSW gluing} to the family $(\widetilde{M},\widetilde{\mathfrak{s}})$, we get
\begin{equation}\label{eq: switching 1}
\FSW(\widetilde{M},\widetilde{\mathfrak{s}})= \langle\HMfrom_{*}(\widetilde{W}_{01},\widetilde{\mathfrak{s}}|_{\widetilde{W}_{01}})(\widehat{1})\ ,\  \HMarrow^{*}(W_{12},\widetilde{\mathfrak{s}}|_{W_{12}})(\widecheck{1})\rangle.
\end{equation}

For a single $4$-manifold $(M,\mathfrak{s})$, there is the following gluing formula  
\begin{equation}\label{eq: switching 2}
\textrm{SW}(M,\mathfrak{s})= \langle U^{\frac{d(\mathfrak{s},M)}{2}}\cdot \HMfrom_{*}(W_{01},\mathfrak{s}|_{W_{01}})(\widehat{1})\ ,\  \HMarrow^{*}(W_{12},\widetilde{\mathfrak{s}}|_{W_{12}})(\widecheck{1})\rangle.
\end{equation}
which generalizes Theorem \ref{thm: FSW gluing} (see \cite[Proposition 27.4.1]{KronheimerMrowkaMonopole}).
Applying Theorem \ref{thm: computing cobordism induced map} to the single cobordism $(W_{01},\mathfrak{s}|_{W_{01}})$, we get 
\begin{equation}\label{eq: switching 3}
\HMfrom_{*}(W_{01},\mathfrak{s}|_{W_{01}})(\widehat{1})=\hat{e}_{d'-1},    
\end{equation}
where $d'=\frac{c^{2}_{1}(\mathfrak{s}|_{M_1})-\sigma(M_1)}{4}$. Applying Theorem \ref{thm: computing cobordism induced map} again to the family cobordism $(\widetilde{W}_{01},\widetilde{\mathfrak{s}}|_{\widetilde{W}_{01}})$, we get that 
\begin{equation}\label{eq: switching 4}
\HMfrom_{*}(\widetilde{W}_{01},\widetilde{\mathfrak{s}}|_{\widetilde{W}_{01}})(\widehat{1})=\langle c_{\frac{\dimension\Q}{2}}(-\Ind\{\slashed{D}_{A^{0}_{b},\hat{\mathfrak{q}}_{b}}^{+}(W_{01,b})\}),[\Q]\rangle \cdot \widehat{e}_{d'-d(\mathfrak{s},M)-1}.
\end{equation}
We combine (\ref{eq: switching 1}), (\ref{eq: switching 2}), (\ref{eq: switching 3}),  (\ref{eq: switching 4}) and we note that $$U^{\frac{d(\mathfrak{s},M)}{2}}\cdot \widehat{e}_{d'-1}=\widehat{e}_{d'-d(\mathfrak{s},M)-1}.$$
This proves (\ref{eq: FSW=SW 1}). The proof of (\ref{eq: FSW=SW 2})  is completely analogous.
\end{proof}
\section{ Families of 4-manifolds constructed from ADE-type singularities}\label{section: ADE}
In this section, we construct  families of 4-manifolds over $S^{2}$ which arise when one resolves an ADE-type singularity using a hyperh\"ahler family of  complex structures near the singularity. See  \cite{Giansiracusa09,Giansiracusa2021,Kronheimer1989} for some previous studies on these natural families. 

Recall that an ADE singularity is a singularity that locally models on the quotient $\mathbb{C}^{2}/\Gamma$, where $\Gamma$ is a finite subgroup of $\SU(2)$. Such a group correspond to a single laced Dykin diagram $G$ of type A, type D and type E. An ADE singularity carries a minimal resolution $$p: \widetilde{\mathbb{C}^{2}/\Gamma} \rightarrow \mathbb{C}^{2}/\Gamma.$$ The preimage $p^{-1}(0)$ is the union of a collection of $-2$-spheres $\{S_{1},\cdots S_{n}\}$, which one-to-one correspond to the vertices $\{v_{1},\cdots ,v_{n}\}$ of $G$. $S_{i}$ and $S_{j}$ intersects transversely at a single point when there is an edge connecting $v_{i}$ and $v_{j}$ and they don't intersect otherwise.  (See, for example, \cite{Reid1997} for more detailed discussions.)

We use $D(V)$ (resp. $S(V)$) to denote the unit disk (resp. unit sphere) in a normed vector space $V$. Let $M_{1}=p^{-1}(D(\mathbb{C}^{2})/\Gamma)$. Then $M_{1}$ is a 4-manifold whose boundary $Y=S(\mathbb{C}^{2})/\Gamma$ is a rational homology sphere that carries a positive scalar curvature. Note that $M_{1}$ can be obtained by plumbing disk bundles over $S_{i}$. Therefore, $M_{1}$ deformation retracts to $p^{-1}(0)$. In particular, it is simply connected and $H_{2}(M_{1})$ is freely generated by $\{S_{i}\}$. The intersection form of $M_{1}$ is the matrix associated to the Dynkin diagram $G$. Hence it is not hard to check that $M_{1}$ is negative definite. 

Now we construct a smooth family over $\Q=S^{2}$ with fiber $M_{1}$. We identify $\mathbb{C}^{2}$ with $\mathbb{H}$ by sending $(z_{1},z_{2})$ to $z_{1}+\jj\cdot z_{2}$ and identify $S(\mathbb{H})$ with $\SU(2)$ by sending $w_{1}+ \jj w_{2}$ to $\big(\begin{smallmatrix}
  w_{1} & -\bar{w_{2}}\\
  w_{2} & \bar{w_{1}}
\end{smallmatrix}\big)$. Then the standard action $SU(2)$ on $\mathbb{C}^{2}$ can be equivalently described as the left multiplication of $S(\mathbb{H})$ on $\mathbb{H}$. Consider the space
\begin{equation}\label{eq: complex structure base}
\Q= \{a\ii+b\jj +c\kk \mid a,b,c\in \mathbb{R},\  a^{2}+b^{2}+c^{2}=1\}=S^{2}. 
\end{equation}
with base point $b_{0}=\ii$.
We use the ordered basis $\{\ii,\jj,\kk\}$ to orent the space they span and to orient $\Q$. 
For each $b\in \Q$, the right multiplication by $b$ provides a complex structure $J_{b}$ on $\mathbb{H}$. Since right multiplications commute with left multiplications, the action of $\Gamma$ preserve all these complex structures. Therefore, we can take the product space $\Q\times (D(\mathbb{H})/\Gamma)$ and simultaneously resolve the singularities in the fibers $\{\{b\}\times (D(\mathbb{H})/\Gamma)\}_{b\in \Q}$ using the complex structures $\{J_{b}\}_{b\in \Q}$. This gives us a smooth family $\widetilde{M}_{1}/\Q$ of 4-manifolds whose boundary carries a natural trivialization. We use $M_{1,b}$ to denote the fiber over $b$.  

Given any smooth 4-manifold $M_{2}$ with fiber $-Y$, we can complete the family $\widetilde{M}_{1}/\Q$ into a family 
\begin{equation}\label{eq: closed ADE family}
\widetilde{M}=\widetilde{M_{1}}\cup_{\Q\times Y} \widetilde{M_{2}}    
\end{equation}
of closed 4-manifolds.
Here $\widetilde{M_{2}}$ is the product family $\Q\times M_{2}$. We call $\widetilde{M}/\Q$ the ADE family with fiber $M=M_{1}\cup M_{2}$.

The clutching function of $\widetilde{M}_{1}/\Q$ is an element in $\pi_{1}(\Diff(M_{1},\partial M_{1}))$. We give an alternative description of this element. Let $S^{1}=\{e^{\ii \theta}\}\subset S(\mathbb{H})$. Consider the $S^{1}$-action on $\mathbb{H}$ by the right multiplication. Since this action preserves both the action of $\Gamma$ and the complex structure $J_{\ii}$, it lifts to an $S^{1}$-action 
\begin{equation}\label{eq: circle action}
\phi:S^{1}\times M_{1}\rightarrow M_{1}.     
\end{equation}
We want to modify $\phi_{1}$ so that it fix the boundary pointwisely. We take a null-homotopy of the inclusion $S^{1}\hookrightarrow S(\mathbb{H})$
\begin{equation}\label{eq: homotopy of circle action}
h:S^{1}\times [\frac{1}{2},1]\rightarrow S(\mathbb{H})
\end{equation}
that satisfies 
$
h(e^{\ii \theta },t)=e^{\ii \theta}
$ for any $t\in [\frac{1}{2},\frac{2}{3}]$ and $
h(e^{\ii \theta },t)=1
$ for any $t\in [\frac{5}{6},1]$. Then we consider
the map $
\phi'_{1}: S^{1}\times M_{1}\rightarrow M_{1}
$
defined by 
\begin{equation}\label{eq: explicit clutching function}
\phi'(e^{\ii \theta}, x)=
\begin{cases}
\phi_{1}(e^{\ii \theta},x) &\mbox{if } |p(x)|<\frac{1}{2}\\
p^{-1}(p(x)\cdot h(e^{\ii \theta},t)) &\mbox{if } |p(x)|\in [\frac{1}{2},1]
\end{cases}.
\end{equation}
Because $\phi'(e^{\ii \theta},-)$ fixes the boundary pointwisely, we obtain a loop  $\gamma$ in $\Diff(M_{1},\partial M_{1})$. Because $\pi_{2}(S(\mathbb{H}))=0$, the homotopy class of $\gamma$ doesn't depend on $h$. We mentioned that for an $A_{1}$-sigularity, $p^{-1}(0)$ consists of a single $-2$-sphere $S$ and $\gamma$ is exactly the generalized Dehn twist along $S$ introduced in Section \ref{section: introduction}.

\begin{rmk}\label{rmk: clutching}
We specify our convention of clutching function as follows: Given a topological space $F$ and a loop $\eta=\{\eta_{\theta}:F\rightarrow F\}_{\theta\in [0,2\pi]}$ of homeomorphisms based at the identity. We can form a fiber bundle $\widetilde{F}$ over $S^{2}=D(\mathbb{C})/ S(\mathbb{C})$, whose total space is defined as 
$$
(D(\mathbb{C})\times F)/(e^{\ii \theta},x)\sim (1,\eta_{\theta}(x)).
$$
We say the clutching function of this bundle $\widetilde{F}$ is represented by $\eta$. Under this convention, suppose we set $F=\mathbb{C}$ and let $\eta$ be the circle action of weight $k$. Then $\widetilde{F}$ has Euler number $-k$.  
\end{rmk}
\begin{lem}\label{lem: clutching function}
$\gamma$ represents the clutching function of the family $\widetilde{M}_{1}/\Q$.
\end{lem}
\begin{proof}
For any $\theta\in [0,2\pi]$ and $\xi \in [0,\pi]$, we let $$b_{\theta,\xi}= \cos \xi \cdot \ii +\sin \xi \cos \theta \cdot \jj +\sin \xi \sin \theta \cdot  \kk.$$ Consider a map $f_{\theta,\xi}:D(\mathbb{H})\rightarrow D(\mathbb{H})$ defined by 
$$
f_{\theta,\xi}(x)=
\begin{cases}
x\cdot b_{\theta,\frac{\xi}{2}} &\mbox{if } |x|<\frac{1}{2}\\
x\cdot (\tau(|x|)+\sqrt{1-\tau(|x|)^{2}} b_{\theta,\frac{\xi}{2}}) &\mbox{if } |x|\in [\frac{1}{2},1]
\end{cases},
$$
where $\tau: [\frac{1}{2},1]\rightarrow [0,1]$ is a non-decreasing function that equals $0$ near $\frac{1}{2}$ and equals $1$ near $1$. Note that we have 
$$
b^{-1}_{\theta,\frac{\xi}{2}}\cdot  \ii \cdot b_{\theta,\frac{\xi}{2}} = b_{\theta,\xi}.
$$
Therefore, near $0$, the map $f_{\theta,\xi}$ sends $J_{i}$ to $J_{b_{\theta,\xi}}$. As a result, it induces a diffeomorphism $$\widetilde{f}_{\theta,\xi}:M_{1,\ii}\rightarrow M_{1,b_{\theta,\xi}}.$$ By definition, the clutching function is given by the family of diffeomorphisms $$\{\widetilde{f}_{0,\pi}^{-1} \circ \widetilde{f}_{\theta,\pi}: M_{1,\ii}\rightarrow M_{1,\ii}\}_{\theta\in [0,2\pi]}.$$
Note that for $|x|<\frac{1}{2}$, we have 
$$
f_{0,\pi}^{-1} \circ f_{\theta,\pi}(x)= x\cdot b_{\theta,\frac{\pi}{2}}\cdot b^{-1}_{0,\frac{\pi}{2}}=x (\cos \theta+ \sin\theta \cdot \ii).
$$
Therefore, the loop in $\Diff(M_{1},\partial M_{1})$ given by the family  $\{\widetilde{f}_{0,\pi}^{-1} \circ \widetilde{f}_{\theta,\pi}\}$ is exactly $\gamma$ for a suitable choice of $h$.
\end{proof}
The following lemma proves the existence of fiberwise almost complex structures on $\widetilde{M}$.
\begin{lem}\label{almost complex structure}
Let $J$ be an almost complex structure on $M=M_{b_{0}}$. Suppose  $c_{1}(TM_{1},J)=0$. Then there exists a fiberwise almost complex structure $\{J_{b}\}_{b\in \Q}$ on fibers of $\widetilde{M}$ that extends $J$.

\end{lem}
\begin{proof}
Since $M_{1}=\widetilde{D(\mathbb{H})/\Gamma }$ deformation retracts to $\cup_{1\leq i\leq n} S_{i}$, an almost complex structure on $M_{1}$ is determined by its first Chern class. Therefore, after a suitable homotopy, we may assume that $J|_{M_{1}}$ is inherited from the complex structure $J_{\ii}$ on $\mathbb{H}$.

Let $Y=\partial M_{1}=S(\mathbb{H})/\Gamma $. Then for each $b\in\Q$, the complex structure  $J_{b}$ on $\mathbb{H}$ induces an almost complex structure on $TM_{1}|_{Y}$, denoted by $J_{Y,b}$. Then $J|_{Y}=J_{Y,i}$.

Next, we consider the fiberwise almost complex structure on the following two families.
\begin{itemize}
\item Since the ADE family $\widetilde{M_{1}}/\Q$ is obtained fiberwisely resolving $\Q\times (D(\mathbb{H})/\Gamma )$ using $\{J_{b}\}_{b\in \Q}$, there is an induced fiberwise almost complex structure $\{J^{1}_{b}\}_{b\in \Q}$ on $\widetilde{M_{1}}/\Q$ such that $J^{1}_{b}|_{Y}=J_{Y,b}$.
\item Take the product family with fiber $M_{1}$, denoted by $\widetilde{M'_{1}}/\Q$. As proved in \cite{Kronheimer1989}, the space $M^{*}_{1}:=\widetilde{\mathbb{H}/\Gamma}$ is an ALE space. In particular, it admits a family of complex structuers $\{J_{M^{*}_{1},b}\}_{b\in \Q}$ such that  $J_{M^{*}_{1},b}$ limits to $J_{b}$ near infinity. By restricting these complex structures to $M_{1}$ and applying a suitable homotopy near the boundary, we obtain a fiberwise almost complex structure $\{J^{2}_{b}\}_{b\in \Q}$ on $\widetilde{M}'_{1}/\Q$ such that $J^{2}_{b}|_{Y}=J_{Y,b}$.
\end{itemize}
We glue together the families $\widetilde{M_{1}}/\Q$ and $\widetilde{M'_{1}}/\Q$ along their common fiber at $b_0$ and form a family $\widetilde{M_{1}}\cup_{M_{1}}\widetilde{M'_{1}}$ over $\Q\vee\Q$.  Then we pull back this bundle via a map $\rho: \Q\rightarrow \Q\vee \Q$ that is of degree 1 on the first component and degree $-1$ on the second component. We denote the resulting family by $\widetilde{M}''_{1}/\Q$ and denote the pulled-back fiberwise almost complex structure by $\{J^{3}_{b}\}_{b\in \Q}$. Note that $J^{3}_{b}|_{Y}= J_{Y,\rho'(b)}$, where $\rho':\Q\rightarrow \Q$ is the composition $$\Q\xrightarrow{\rho}\Q\vee \Q\xrightarrow{\operatorname{Id}_{\Q}\vee \operatorname{Id}_{\Q}}\Q.$$ Since $\rho'$ is null-homotopic, we may homotope $\{J^{3}_{b}\}_{b\in \Q}$ near the boundary to obtain a fiberwise almost complex structure  $\{J^{4}_{b}\}_{b\in \Q}$ on $\widetilde{M}''_{1}/\Q$ such that $J^{4}_{b}|_{Y}=J_{Y,\ii}$ for any $b\in \Q$. Note that  $\widetilde{M}''_{1}/\Q$ and  $\widetilde{M}_{1}/\Q$ are isomorphic. So we can also treat $\{J^{4}_{b}\}_{b\in \Q}$ as a fiberwise almost complex structure on $\widetilde{M}_{1}/\Q$. Let $J_{M,b}$ be obtained by gluing together $J^{4}_{b}$ and $J|_{M_{2}}$ along $Y$. Then $\{J_{M,b}\}_{b\in \Q}$ is a fiberwise almost complex structure on $\widetilde{M}/\Q$ that extends $J$.

\end{proof}

Next, we give an explicit description of the circle action (\ref{eq: circle action}). To do this, we need an alternative construction of $M_{1}$ by plumbing disk bundles. 

We first recall the plumbing construction as follows: Let $D^{2}\hookrightarrow V\rightarrow S$ and $D^{2}\hookrightarrow V'\rightarrow S'$ be two disk bundles over surfaces $S,S'$. Take points $x\in S$ and $x'\in S'$. Let $B$ (resp. $B'$) be a small neighborhood of $x$ (resp. $x'$) in $S$ (resp. $S'$).  
Fix trivializations 
$$
V|_{B}\cong D^{2}\times B,\quad  V|_{B'}\cong D^{2}\times B'.
$$
The plumbing of $V, V'$ along $x, x'$ is defined by gluing $V$ and $V'$ together by  identifying $V|_{B}$ with $V_{i}|_{B_{i}}$ in a way that switches the fiber and base directions. We use the notation $(V,x)\sim (V',x')$ to denote this procedure. 

Recall that $p^{-1}(0)\in \widetilde{\mathbb{C}^{2}/\Gamma}$ is a union of $-2$-spheres $S_{1},\cdots ,S_{n}$ that intersects each other transversely. Let $V_{i}$ be the normal disk bundle of $S_{i}$ in $\widetilde{\mathbb{C}^{2}/\Gamma}$. Since the action (\ref{eq: circle action}) preserve $S_{i}$ as a set, it induces an action $$r_{i}:S^{1}\times V_{i}\rightarrow V_{i}.$$ 
There is a homeomorphism from $M_{1}$ to the plumbing of $\{V_{i}\}_{1\leq i\leq n}$ (along the intersection points $\{S_{i}\cap S_{j}\}_{1\leq i<j\leq n}$). Furthermore, under this homeomorphism, the circle action (\ref{eq: circle action}) is obtained by 
plumbing together $\{r_{i}\}_{1\leq i\leq n}$.

To give an explicit description of $r_{i}$, we fix a complex analytic diffeomorphism  $$l_{i}: \mathbb{C}\cup \{\infty\}=\mathbb{CP}^{1}\xrightarrow{\cong } S_{i}$$ for each $i$. These diffeomorphisms are chosen in a way such that the following conditions hold: 
\begin{itemize}
    \item For a type-A singularity, we have $l_{i}(\infty)=l_{i+1}(0)$ for $1\leq i\leq n-1$. In other words, $S_{i}$ and $S_{i+1}$ intersects each other at a single point $l_{i}(\infty)=l_{i+1}(0)$. 
    \item For a type-D singularity, we have $l_{i}(\infty)=l_{i+1}(0)$ for $1\leq i\leq n-2$ and $l_{n-2}(1)=l_{n}(0)$.
    \item For a type-E singularity, we have $l_{i}(\infty)=l_{i+1}(0)$ for $1\leq i\leq n-2$ and $l_{n-3}(1)=l_{n}(0)$.
    \item For any $i$, both $l_{i}(0)$ and $l_{i}(\infty)$ are fixed under the action $r_{i}$. Since $r_{i}$ must fix points in $S_{i}\cap S_{j}$ (for $j\neq i$), this condition is redundant unless the  vertex in the  Dynkin diagram that corresponds to $S_{i}$ is 1-valent.
\end{itemize}
Using $l_{i}$, we can identify $V_{i}$ with the unit disk bundle associated to $\mathcal{O}(-2)$. There are two natural circle actions. The action $$r^{i}_{(0,1)}:S^{1}\times V_{i}\rightarrow V_{i}$$ is given by the scalar multiplication on the fibers. And the action $$r^{i}_{(1,-1)}:S^{1}\times V_{i}\rightarrow V_{i} $$ is induced by the standard circle action on $\mathbb{CP}^{1}=\mathbb{C}\cup {\infty}$ (Here we identify $\mathcal{O}(-2)$ with $T\mathbb{CP}^{1}$ via an isomorphism that \emph{reverses} the fiber orientation). Both actions have $0\in \mathbb{CP}^{1}$ as a fixed point, and the weights are $(0,1)$ and $(1,-1)$ respectively. (Here the first and second component denote weights in the base and the fiber direction respectively.) For any $a,b\in \mathbb{Z}$, we can consider the composition 
$$
r^{i}_{(a,b)}:S^{1}\times V_{i}\rightarrow V_{i}
$$
defined by 
$$
r^{i}_{(a,b)}(z,x):= r^{i}_{(0,1)}(z^{a+b},r^{i}_{(1,-1)}(z^{a},x)),\quad\text{for } z\in S^{1},\ x\in V.$$
This is a circle action that has $0$ as a fix point of weight $(a,b)$. 
Note that for any $z\in S^{1}$, $r_{i}(z,-):V_{i}\rightarrow V_{i}$ is a bundle automorphism that covers an analytical map $S_{i}\rightarrow S_{i}$ which fixes both $0$ and $\infty$. Such action is completely determined by the weights at $0$ in the fiber and base direction. Therefore, $r_{i}=r^{i}_{(a_{i},b_{i})}$ for some $a_{i},b_{i}\in \mathbb{Z}$. In other words, the circle action $\phi$ on $M_{1}$ is obtained by plumbing the circle actions $\{r^{i}_{(a_{i},b_{i})}\}$ together for some $a_{i},b_{i}$. Our next task is to find the exact values of them. 

\begin{lem}\label{lem: glue plumbing}
(1) The weight of $r_{(a,b)}$ at $\infty$ is $(-a,2a+b)$. Again, we put the base direction in the first component and the fiber direction in the second component\\
(2) Suppose $V_{i}$ and $V_{j}$ are plumbed together along $l_{i}(0)\in S_{i}$ and $l_{j}(\infty)\in S_{j}$. Then we have  $(a_{j},b_{j})=(2a_{i}+b_{i},-a_{i})$.
\\
(3) Suppose $V_{i}$ and $V_{j}$ are plumbed together along $l_{i}(1)\in S_{i}$ and $l_{j}(0)\in S_{j}$. Then $(a_{i},b_{i})=(0,k)$ an $(a_{j},b_{j})=(k,0)$ for some integer $k$.
\end{lem}
\begin{proof}
(1) From the geometric description of $r^{i}_{(0,1)}$ and $r^{i}_{(1,-1)}$, we see that their weights at $\infty$ are $(0,1)$ and $(-1,1)$ respectively. So the weight of $r_{(a,b)}$ at $\infty$ is $(-a,2a+b)$. 

(2) Since the plumbing construction switches the base and the fiber direction, we must have $(b',a')=(2a+b,-a)$ if $r_{(a,b)}$ and $r_{(a',b')}$ can be plumbed together. 

(3) Since $r^{i}_{(a_{i},b_{i})}$ fixes $l(0),l(1)$ and $l(\infty)$, it must be trivial when restricted to $S_{i}$. This implies that $(a_{i},b_{i})=(0,k)$ and $(a_{j},b_{j})=(k,0)$ for some $k$. 
\end{proof}

\begin{pro}\label{values of a,b}The values of $\{(a_{i},b_{i})\}_{1\leq i\leq n}$ are given as follows:
\begin{enumerate}
    \item For an $A_{n}$-singularity, we have $(a_{i},b_{i})=(-n+2i-1,n-2i+3)$ for $1\leq i\leq n$.
    \item For a $D_{n}$-singularity, we have $(a_{i},b_{i})=(-2n+2i+4,2n-2i-2)$ for $1\leq i\leq n-1$ and $(a_{n},b_{n})=(2,0)$.
    \item For an $E_{n}$-singularity, we have $(a_{i},b_{i})=(-2n+2i+6,2n-2i-4)$ for $1\leq i\leq n-1$ and $(a_{n},b_{n})=(2,0)$.
\end{enumerate}
\end{pro}
\begin{proof}
Since $\overline{p^{-1}(D(\mathbb{C}^{2})\setminus \{0\})}=M_{1}$, the action $\phi$ is determined by its restriction to $p^{-1}(D(\mathbb{C}^{2})\setminus \{0\})$. Furthermore, suppose we identify $p^{-1}(D(\mathbb{C}^{2})\setminus \{0\})$ with $\partial M_{1}\times (0,1]$ by sending  $p^{-1}(tx)$ to $(x,t)$. Then $\phi$ is trivial on the second component. Therefore, to understand $\phi$, it suffices to understand its restriction to $\partial M_{1}$, which is just the circle action $$
\phi': S^{1}\times (S(\mathbb{C}^{2})/\Gamma)\rightarrow (S(\mathbb{C}^{2})/\Gamma)
$$
defined by $$\phi'(e^{\ii \theta},[(z_{1},z_{2})])=[(e^{\ii \theta}z_{1}, e^{\ii \theta}z_{2})].$$
This circle action gives $S(\mathbb{C}^{2})/\Gamma$ the structure of a Seifert fibered space. Each fiber (regular or singular) is an orbit by itself. The isotropy group for a regular fiber is the group
$
\Gamma\cap \{\big(\begin{smallmatrix}
  e^{\ii \theta } & 0\\
  0 & e^{\ii \theta }
\end{smallmatrix}\big)\}
$. This group is trivial for an $A_{n}$-singularity with $n$ even and is $\mathbb{Z}/2$ otherwise.

For a $D_{n}$ singularity, $v_{n-2}$ is the unique 3-valent vertex in the Dynkin diagram. Take a fiber of the corresponding disk bundle $V_{n-2}$. Then the boundary of this fiber is a regular fiber of  $S(\mathbb{C}^{2})/\Gamma$. Since any regular fiber is an orbit with isotropy group $\mathbb{Z}/2$, we see that $(a_{n-2},b_{n-2})=(0,2)$. (We don't have $(0,-2)$ here because we are using the complex orientation to orient the fibers.) The value of $(a_{i},b_{i})$ for $i\neq n-2$ is determined by Lemma \ref{lem: glue plumbing}. This finishes the proof of (2). The proof of (3) is completely analogues. And $v_{n-3}$ is the unique 3-valent vertex in this case. 

Now we prove (1). For an $A_{n}$-singularity, we have  $$\Gamma=\{\big(\begin{smallmatrix}
  e^{\frac{2k\pi \ii }{n+1}}  & 0\\
  0 & e^{-\frac{2k\pi \ii }{n+1}}
\end{smallmatrix}\big)\}_{k\in \mathbb{Z}},$$ and $\partial M_{1}=S(\mathbb{C}^{2})/\Gamma$ is the lens space $L(n+1,-1)$. Suppose $n=1$. Then there are no singular fibers. Same argument as before shows that $(a_{1},b_{1})=(0,2)$. Now suppose $n\geq 2$, then there are two singular fibers $\{[(e^{\ii\theta},0)]\}_{\theta}$ and $\{[(0,e^{\ii\theta})]\}_{\theta}$. For points on these two fibers, the isotropy group is $\mathbb{Z}/(n+1)$, which is strictly larger than the isotropy group of points on regular fibers, which is either trivial or $\mathbb{Z}/2$. Now we consider the fiber of $V_{1}$ over $l_{1}(0)$ and the fiber of $V_{n}$ over $l_{n}(\infty)$. Denote their boundaries by $F_{1}$ and $F_{n}$. By our assumption that $l_{1}(0)$ and $l_{1}(\infty)$ are fixed points of $r$, $F_{1}$ and $F_{n}$ are orbits of the circle action on $\partial M_{1}$. Furthermore, from the plumbing description of $\phi$, we see that they are the only orbits whose isotropy group can possibly be larger then orbits near by. (Boundaries of the fibers over $l_{i}(0)$ with $i>1$ or $l_{i}(\infty)$ with $i<n$ also have this property, but they are contained in the interior of $M_{1}$.) Therefore, we see that $F_{0}$, $F_{1}$ are exactly the two singular fibers and points on them have isotropy group $\mathbb{Z}/(n+1)$. In particular, this implies that the weights of $\phi$ at $l_{1}(0)$ and $l_{n}(\infty)$ are $(c,n+1)$ and $(c',n+1)$ for some $c,c'\in \mathbb{Z}$. By Lemma \ref{lem: glue plumbing} (1), this implies  $$b_{1}=2a_{n}+b_{n}=n+1.$$ Putting this equation with the equations $$(a_{i+1},b_{i+1})=(2a_{i}+b_{i},-a_{i})\quad \text{for }1\leq i\leq n-1$$
provided by Lemma \ref{lem: glue plumbing} (2), one can solve that $(a_{i},b_{i})=(-n+2i-1,n-2i+3)$. This finishes the proof. 
\end{proof}

Since the circle action $\phi$ fixes $S_{i}$ as a subset, we obtain a 4-dimensional submanifold $\widetilde{S}_{i}$ of $\widetilde{M}_{1}$ by putting together $S_{i}$ in all fibers. A simple application of the Serre spectral sequence shows that $H_{4}(\widetilde{M}_{1};\mathbb{Z})$ is a free abelian group generated by $\{[\widetilde{S}_{i}]\}_{1\leq i\leq n}$. 

\begin{lem}\label{lem: self-intersection} (1) The triple self-intersection number  $\widetilde{S}_{i}\cdot \widetilde{S}_{i}\cdot \widetilde{S}_{i}=8$ for any $i$. (2) Suppose $S_{i}$ and $S_{j}$ intersects each other at $x\in S_{i}$. And suppose the weight of the circle $\phi$ action on $T_{x}S_{i}$ equals $k$. Then $\widetilde{S}_{i}\cdot \widetilde{S}_{j}\cdot \widetilde{S}_{j}=-k$. (3) We have $\langle p_{1}(T\widetilde{M}_{1}/\Q),[\widetilde{S}_{i}]\rangle=8$. 
\end{lem}
\begin{proof} We have a fiber bundle $S_{i}\hookrightarrow \widetilde{S}_{i}\rightarrow S^{2}$ and the clutching function is the restriction of $\phi$. A fixed point $p$ of $\phi$ gives a  section of the bundle $\widetilde{S}_{i}$. 
The image of this section, denoted by $\Sigma_{p}$, is a 2-dimensional submanifold of $\widetilde{S}_{i}$. The normal bundle of $\Sigma_{p}$ in $\widetilde{S}_{i}$, denoted by $N_{\Sigma_{p}}\widetilde{S}_{i}$, has Euler number $-k$, where $k$ is the weight of $\phi$ on $T_{p}S_{i}$. (See Remark \ref{rmk: clutching}.)

To prove (1), we set $p=l_{i}(0)$. Since the weight of $\phi$ on $T_{p}S_{i}$ equals $a_{i}$, we have 
\[
\widetilde{S}_{i}\cong \begin{cases}
S^{2}\times S^{2}
&\text{ if $a_{i}$ is even} \\ \mathbb{CP}^{2}\# \widebar{\mathbb{CP}}^{2} &\text{ if $a_{i}$ is odd}
\end{cases}.
\]
In any case, $\{[S_{i}], [\Sigma_{p}]\}$ is a basis of $H_{2}(\widetilde{S}_{i};\mathbb{Q})$. Their intersection numbers in $\widetilde{S}_{i}$ are given by 
\begin{equation}\label{eq: double intersection}
S_{i}\cdot S_{i}= 0,\quad  \Sigma_{p}\cdot  \Sigma_{p}=-a_{i},\quad S_{i} \cdot  \Sigma_{p}=1. 
\end{equation}
Let $N$ be the normal bundle of  $\widetilde{S}_{i}$ in $\widetilde{M_{1}}$. Then $N|_{S_{i}}$ has Euler number $-2$ (because $S_{i}$ is a $-2$ sphere) and  $N|_{\Sigma_{p}}$ has Euler number $-b_{i}$. By (\ref{eq: double intersection}), we get 
$$
e(N)=\operatorname{PD}(-2 [\Sigma_{p}]-(2a_{i}+b_{i})[S_{i}]).
$$
Therefore, we get 
$$
\widetilde{S}_{i}\cdot \widetilde{S}_{i}\cdot \widetilde{S}_{i}=\langle e(N)^{2},[\widetilde{S}_{i}]\rangle=4(a_{i}+b_{i})=8.
$$

To prove (2), we set $p=x$. Note that $\widetilde{S}_{i}$ and $\widetilde{S}_{j}$ intersect each other transversely at $\Sigma_{p}$. So 
$\widetilde{S}_{i}\cdot \widetilde{S}_{j}\cdot \widetilde{S}_{j}$ equals the self-intersection of $\Sigma_{p}$ in $\widetilde{S}_{i}$, which is $-k$.

Now we prove (3). We have a decomposition 
$$
(T\widetilde{M}_{1}/\Q)|_{\widetilde{S}_{i}}=(T\widetilde{S}_{i}/\Q)\oplus N_{\widetilde{S}_{i}}\widetilde{M}_{1}.
$$
Here $T\widetilde{S}_{i}/\Q$ denotes the vertical tangle bundle of the fiber bundle $S_{i}\hookrightarrow \widetilde{S}_{i}\rightarrow \Q$. This implies that 
$$
p_{1}((T\widetilde{M}_{1}/\Q)|_{\widetilde{S}_{i}})=e(N_{\widetilde{S}_{i}}\widetilde{M}_{1})^{2}+e(T\widetilde{S}_{i}/\Q)^{2}.
$$

Set $p=l_{i}(0)$. Then we have  $\chi((T\widetilde{S}_{i}/\Q)|_{S_{i}})=\chi (TS_{i})=2$ and $\chi(T\widetilde{S}_{i}|_{\Sigma_{p}})=-a_{i}$. By the same argument as (2), we see that 
$$
e(T\widetilde{S}_{i}/\Q)=\operatorname{PD}(2 [\Sigma_{p}]+a_{i}[S_{i}]).
$$

This implies 
$$
\langle e(T\widetilde{S}_{i}/\Q)^{2}, [\widetilde{S}_{i}]\rangle=0.
$$
Hence we have 
$$
\langle p_{1}(T\widetilde{M}_{1}/\Q), [\widetilde{S}_{i}]\rangle =\langle e(N_{\widetilde{S}_{i}}\widetilde{M}_{1})^{2}, [\widetilde{S}_{i}]\rangle =\widetilde{S}_{i}\cdot \widetilde{S}_{i}\cdot \widetilde{S}_{i}=8.
$$
\end{proof}

\begin{rmk}\label{rmk: topologically nontrivial}
The ADE family $\widetilde{M}/Q$ is always topologically nontrivial because the  total space $\widetilde{M}$ is not homotopy equivalent to $M\times \Q$. To see this, we note that for any $\alpha \in H_{4}(M\times \Q;\mathbb{Z})$, the triple self-intersection number $\alpha\cdot \alpha\cdot \alpha $ is always divisible by $3$ (this can be seen using the K\"unneth formula). However, this condition is not satisfied  for $\widetilde{M}$ because $\widetilde{S}_{i}\cdot \widetilde{S}_{i}\cdot \widetilde{S}_{i}=8$. Similarly, one can show that the blown-up family is also topologically nontrivial. (See Lemma \ref{lem: self-intersection in blown-up}.) 
\end{rmk}

Given any $r=(t_{1},\cdots, t_{n})\in \mathbb{Z}^{n}$, we let $S_{r}$ be the formal linear combination $t_{1}S_{1}+\cdots+t_{n}S_{n}$. Similarly, we let $\widetilde{S}_{r}=t_{1}\widetilde{S}_{1}+\cdots+t_{n}\widetilde{S}_{n}$. 

\begin{pro}\label{pro: calculation of index}
Let $\widetilde{M}=\widetilde{M}_{1}\cup (\Q\times M_{2})$ be an ADE family with fiber $M=M_{1}\cup M_{2}$. Let $\widetilde{\mathfrak{s}}$ be a family $\spinc$ structure on $\widetilde{M}$ that satisfies 
\begin{equation}\label{eq: assumption on c1}
\operatorname{PD}(c_{1}(\widetilde{\mathfrak{s}}|_{\widetilde{M}_{1}}))=2[\widetilde{S}_{r}]\in     H_{4}(\widetilde{M}_{1},\partial \widetilde{M}_{1};\mathbb{Z})
\end{equation}
for some root $r$. Then 
$$\langle c_{1}(\operatorname{Ind}(\slashed{D}^{+}(\widetilde{M},\widetilde{\mathfrak{s}}))),[\Q]\rangle=\pm 1.$$
\end{pro}
\begin{proof} By the family Atiyah-Singer index theorem, one has 
\begin{equation}\begin{split}
\langle c_{1}(\operatorname{Ind}(\slashed{D}^{+}(\widetilde{M},\widetilde{\mathfrak{s}}))),[\Q]\rangle&=  \langle e^{\frac{c_{1}(\widetilde{\mathfrak{s}})}{2}}\cdot \hat{A}(T\widetilde{M}/\Q),[\widetilde{M}]\rangle\\
&=\langle (1+\frac{c_{1}(\widetilde{\mathfrak{s}})}{2}+\frac{c_{1}(\widetilde{\mathfrak{s}})^{2}}{8}+\frac{c_{1}(\widetilde{\mathfrak{s}})^{3}}{48})\cup (1-\frac{p_{1}(T\widetilde{M}/\Q)}{24}), [\widetilde{M}] \rangle\\
&=
\frac{\langle c_{1}(\widetilde{\mathfrak{s}})^{3}-p_{1}(T\widetilde{M}/\Q) c_{1}(\widetilde{\mathfrak{s}}),[\widetilde{M}]\rangle}{48}.
\end{split}
\end{equation}
By (\ref{eq: assumption on c1}), we have $c_{1}(\widetilde{\mathfrak{s}})=\operatorname{PD}(2[\widetilde{S}_{r}]+\widetilde{\beta})$, where $\widetilde{\beta}$ the image of $\beta\times [\Q]$ under the map $H_{4}(M_{2};\mathbb{Z})\rightarrow H_{4}(M;\mathbb{Z})$ for some $\beta\in H_{2}(M;\mathbb{Z})$. As a result, we get $$
\langle c_{1}(\operatorname{Ind}(\slashed{D}^{+}(\widetilde{M},\widetilde{\mathfrak{s}}))),[\Q]\rangle=\frac{ 4[\widetilde{S}_{r}]\cdot [\widetilde{S}_{r}]\cdot [\widetilde{S}_{r}]- \langle p_{1}(T\widetilde{M}/\Q),[\widetilde{S}_{r}]\rangle} {24}
$$

 We consider the function
$$
f(t_{1},\cdots, t_{n}):=\sum\limits_{1\leq i\leq n}8(4t^{3}_{i}-t_{i})-\sum\limits_{\{\{i,j)\}\mid i\leftrightarrow j \}} 12t_{i}t_{j}(t_{i}w_{j,i}+t_{j}w_{i,j}).
$$
Here we use the $ i\leftrightarrow j$ to denote the condition that there is an edge in the Dynkin diagram that connects the vertices $v_i$ and $v_j$, and we use $w_{i,j}$ denotes the weight of the circle action on the tangent space $T_{S_{i}\cap S_{j}}S_{i}$.

Suppose $r=(t_{1},\cdots, t_{n})$. Then by Lemma \ref{lem: self-intersection}, we have 
$$
4[\widetilde{S}_{r}]\cdot [\widetilde{S}_{r}]\cdot [\widetilde{S}_{r}]- \langle p_{1}(T\widetilde{M}/\Q),[\widetilde{S}_{r}]\rangle=f(t_{1},\cdots, t_{n})
$$ So it suffices to show that 
\begin{equation}\label{eq: calculation for roots}
    f(t_{1},\cdots, t_{n})=\pm 24
\end{equation}
for any root $r=(t_{1},\cdots t_{n})$.
 
By changing the sign of $r$, we may assume that it is a positive root (i.e. $t_{i}\geq 0$ for all $i$.) It is obvious that 
$
f(0,\cdots, 0,1,0,\cdots,0)=24.
$
So (\ref{eq: calculation for roots}) holds for simple roots. Since the Hasse diagram for ADE Dynkin diagram is connected (see, for example, \cite{Eugene98}), 
 for any positive root $r$, there exists a sequence $r_{1},\cdots r_{m}$ of positive roots such that $r_{1}=r$, $r_{m}$ is simple, and 
$$
r_{i}-r_{i+1}=\pm (0,\cdots,1,\cdots, 0).
$$
Therefore, the proof is finished by the following lemma.
\end{proof}
\begin{lem}
Suppose $r=(t_{1},\cdots ,t_{k},\cdots ,t_{n})$ and $r'=(t_{1},\cdots ,t_{k}+1,\cdots ,t_{n})$ are both roots. Then $f(r)=f(r')$.
\end{lem}
\begin{proof}
Note that we always have 
$$
w_{i,j}+w_{j,i}=a_{i}+b_{i}=2.
$$
Also note that since $r\cdot r=r'\cdot r'$, we have 
$$
2t^{2}_{k}- \sum_{\{j\mid j\leftrightarrow k\}}2t_{k}t_{j}=2(t_{k}+1)^{2}- \sum_{\{j\mid j\leftrightarrow k\}}2(t_{k}+1)t_{j},
$$
which implies 
$$
2t_{k}+1= \sum_{\{j\mid j\leftrightarrow k\}}t_{j}.   
$$

It is straightforward to compute that $$
f(r')-f(r)=24(2t_{k}+1)^{2}-\sum_{\{j\mid j\leftrightarrow k\}} 12(t_{j}w_{j,k}(2t_{k}+1)+t^{2}_{j}w_{k,j}). 
$$
Suppose the vertex $v_k$ is 3-valent. Then there are three vertices that are connected to $v_{k}$, and for each of them, we have $w_{k,j}=0$ and $w_{j,k}=2$.
Therefore, we get 
$$
f(r')-f(r)= 24(2t_{k}+1)^{2}-\sum_{\{j\mid j\leftrightarrow k\}} 24 t_{j}(2t_{k}+1)=24(2t_{k}+1)(2t_{k}+1- \sum_{\{j\mid j\leftrightarrow k\}}t_{j}   
)=0.
$$
Suppose the vertex $v_k$ is 2-valent and connected to $v_j$ and $v_{j'}$. Then by Lemma \ref{lem: glue plumbing}, we have
$$
w_{k,j}+w_{k,j'}=a_{i}+(-a_{i})=0.
$$
Let $w=w_{k,j}$. Then $w_{j,k}=2-w$ and $w_{j',k}=2+w$ and we have 
\begin{equation}
\begin{split}
f(r')-f(r)&=
24(2t_{k}+1)^{2}- 12(t_{j}(2-w)(2t_{k}+1)+t^2_{j}w+t_{j'}(2+w)(2t_{k}+1)-t^2_{j'}w)\\&=24(2t_{k}+1)^{2}-12(2(t_{j}+t_{j'})(2t_{k}+1)-w(2t_{k}+1)(t_{j}-t_{j'})+w(t^{2}_{j}-t^2_{j'}))
\\&=24(2t_{k}+1)^{2}-12(2(2t_{k}+1)(2t_{k}+1)-w(t_{j}+t_{j'})(t_{j}-t_{j'})+w(t^{2}_{j}-t^2_{j'}))\\&=0.
\end{split}
\end{equation}
Suppose $v_k$ is 1-valent and connected to $v_j$. Then $t_{j}=2t_{k}+1$ and we have 
$$
f(r')-f(r)= 24(2t_{k}+1)^{2}- 12(2t_{k}+1)^{2}(w_{k,j}+w_{j,k})=0.
$$
\end{proof}
\section{The blown-up family}\label{section: blwon-up}
Let $M\hookrightarrow \widetilde{M}\rightarrow \Q $ be a blown-up family as defined in Section \ref{section: introduction}, with fiber $M=M'\#\overline{\mathbb{CP}^{2}}$. Then we  have a decomposition
$$
\widetilde{M}=\widetilde{M}_{1}\cup \widetilde{M}_{2}
$$
where $\widetilde{M}_{1}$ is a nontrivial family with fiber $M_{1}=\overline{\mathbb{CP}^{2}}\setminus \mathring{D^{4}}$ and $\widetilde{M}_{2}$ is the trivial family with fiber $M_{2}=M'-\mathring{D^{4}}$. We use $S$ to denote the exceptional divisor in the fiber $M$ and use $\widetilde{S}$ to denote submanifold of $\widetilde{M}$ formed by 
the union of the exceptional divisors on all fibers. 
\begin{lem}\label{lem: self-intersection in blown-up}
(1) The triple self-intersection $\widetilde{S}\cdot \widetilde{S}\cdot \widetilde{S}=2$. \\(2) We have $\langle p_{1}(T\widetilde{M}/\Q), [\widetilde{S}]\rangle = 2$.
\end{lem}
\begin{proof} (1) By the same argument as Lemma \ref{lem: clutching function}, we can choose the clutching function for $\widetilde{M}$ so that in a tubular neighborhood of $S$, it is a circle action that fixes $S$ pointwisely and has weight $1$ on its normal bundle in $M$. Therefore, $\widetilde{S}$ is diffeomorphic to $S\times \Q$ and its second homology is generated by $S$ and $\Sigma:=*\times \Q$. By Remark \ref{rmk: clutching}, we have 
$$
\chi(N_{\widetilde{S}}\widetilde{M}|_{\Sigma})=-1
$$
and we also have
$$
\chi((N_{\widetilde{S}}\widetilde{M})|_{S})=\chi(N_{S}M)=S\cdot S=-1.  
$$
This implies that  
$
e(N_{\widetilde{S}}\widetilde{M})=-\operatorname{PD}([S]+[\Sigma])
$
and that
$$
\widetilde{S}\cdot \widetilde{S}\cdot \widetilde{S}=\langle e(N_{\widetilde{S}}\widetilde{M})^{2},[\widetilde{S}]\rangle =2.
$$
(2) We have 
$$
p_{1}(T\widetilde{M}/\Q|_{\widetilde{S}})=e(N_{\widetilde{S}}\widetilde{M})^{2}+e(T\widetilde{S}/\Q)^{2}.
$$
where $T\widetilde{S}/\Q$ denotes the vertical tangent bundle of $S\rightarrow \widetilde{S}\rightarrow \Q$. Since this bundle is trivial, we have $e(T\widetilde{S}/\Q)^{2}=0$. This implies
$$
\langle p_{1}(T\widetilde{M}/\Q|_{\widetilde{S}}),[\widetilde{S}]\rangle =\langle e(N_{\widetilde{S}}\widetilde{M})^{2},[\widetilde{S}]\rangle =2.
$$

\end{proof}
\begin{pro}\label{lem: index calculation 2}
Let $\widetilde{\mathfrak{s}}$ be a family $\spinc$ structure on $\widetilde{M}$ such that  $c_{1}(\widetilde{\mathfrak{s}}|_{\widetilde{M_{1}}})=-3\operatorname{PD}[\widetilde{S}]$. Then 
$$\langle c_{1}(\operatorname{Ind}(\slashed{D}^{+}(\widetilde{M},\widetilde{\mathfrak{s}}))),[\Q]\rangle=-1.$$
\end{pro}
\begin{proof}
The proof is similar to Proposition \ref{pro: calculation of index}. By the Atiyah-Singer index theorem, one has 
\begin{equation}\begin{split}
\langle c_{1}(\operatorname{Ind}(\slashed{D}^{+}(\widetilde{M},\widetilde{\mathfrak{s}}))),[\Q]\rangle&=  \langle e^{\frac{c_{1}(\widetilde{\mathfrak{s}})}{2}}\cdot \hat{A}(T\widetilde{M}/\Q),[\widetilde{M}]\rangle\\
&=\langle (1+\frac{c_{1}(\widetilde{\mathfrak{s}})}{2}+\frac{c_{1}(\widetilde{\mathfrak{s}})^{2}}{8}+\frac{c_{1}(\widetilde{\mathfrak{s}})^{3}}{48})\cup (1-\frac{p_{1}(T\widetilde{M}/\Q)}{24}), [\widetilde{M}] \rangle\\
&=
\frac{-27\widetilde{S}\cdot\widetilde{S}\cdot \widetilde{S}+3\langle p_{1}(T\widetilde{M}/\Q),[\widetilde{S}]\rangle}{48}=-1.
\end{split}
\end{equation}
\end{proof}


\section{Symplectic structures on the ADE family and the blown-up family}\label{section: family symplectic structures existence}
In this section, we study the existence of fiberwise symplectic structures on the ADE families and the blown-up families through the calculations we did the previous two sections.  The following lemma is an improvement of Lemma \ref{lem: extension of spinc b1=0} without the assumption that $b_1(M)=0$.

\begin{lem}\label{lem: extension of spinc structure}
Let $\widetilde{M}/\Q$ be an ADE family or a blown-up family. Then any $\spinc$ structure on the fiber $M$ can be extended to a family $\spinc$ structure on $\widetilde{M}/\Q$. Furthermore, any two such extensions are related to each other by a twist of a complex line  bundle $L$ pulled back from $\Q$
\end{lem}
\begin{proof}
Consider the cohomology Serre spectral sequence for $\widetilde{M}/\Q$ and the differential, $$d_{2}^{0,2}: H^{0}(\Q;H^{2}( M;\mathbb{Z}))\rightarrow H^{2}(\Q;H^{1}( M;\mathbb{Z})).$$
Because $H^{1}(M;\mathbb{Z})\rightarrow H^{1}(M_{2};\mathbb{Z})$ is injective and the corresponding differential for the trivial bundle $\widetilde{M}_{2}/\Q$ is trivial. This gives the short exact sequence 
\begin{equation}\label{eq: short exact sequence 1}
0\rightarrow H^{2}(\Q;\mathbb{Z})\xrightarrow{p^{*}} H^{2}(\widetilde{M};\mathbb{Z})
\xrightarrow{i^{*}} H^{2}(M;\mathbb{Z})\rightarrow 0.
\end{equation}
A similar argument implies the short exact sequence 
\begin{equation}\label{eq: short exact sequence 2}
0\rightarrow H^{2}(\Q;\mathbb{Z}/2)\xrightarrow{p^{*}} H^{2}(\widetilde{M};\mathbb{Z}/2)
\xrightarrow{i^{*}} H^{2}(M;\mathbb{Z}/2)\rightarrow 0.
\end{equation}
We first prove that $\widetilde{M}/\Q$ does have a family $\spinc$ structure. Let $\widetilde{w}_{2}(\widetilde{M}/\Q)$ be the second Stiefel-Whitney class of the vertical tangent bundle $T\widetilde{M}/\Q$. Since any 4-manifold has a $\spinc$ structure , $i^{*}(\widetilde{w}_{2}(\widetilde{M}/\Q))=\omega_{2}(M)$ is the reduction of a class $c\in H^{2}(M;\mathbb{Z})$. By (\ref{eq: short exact sequence 1}), there exists $\widetilde{c}\in H^{2}(\widetilde{M};\mathbb{Z})$ such that $i^{*}(\widetilde{c})=c$. Let $\widetilde{w}$ be the mod-2 reduction of $\widetilde{c}$. Then $i^{*}(\widetilde{w})=i^{*}(\widetilde{w}_{2}(\widetilde{M}/\Q))$. By (\ref{eq: short exact sequence 2}), $\widetilde{w}=\widetilde{w}_{2}(\widetilde{M}/\Q)+p^{*}(w')$ for some $w'\in H^{2}(\Q;\mathbb{Z}/2)$. Let $c'$ be an element in $H^{2}(\Q;\mathbb{Z})$ which reduces to $w'$. Then the mod-2 reduction of $\widetilde{c}-p^{*}(c')$ is exactly $\widetilde{w}_{2}(\widetilde{M}/\Q)$. Therefore, the family $\widetilde{M}/\Q$ has a $\spinc$ structure.

Let $\mathfrak{s}$ be a $\spinc$ structure on the fiber $M$. Take any family $\spinc$ structure $\widetilde{\mathfrak{s}}$. Let  $\beta=\widetilde{\mathfrak{s}}|_{M}-\mathfrak{s}\in H^{2}(M;\mathbb{Z})$. By (\ref{eq: short exact sequence 1}), there exists $\widetilde{\beta}\in H^{2}(\widetilde{M};\mathbb{Z})$ such that $i^{*}(\widetilde{\beta})=\beta$. Then $\widetilde{\mathfrak{s}}+\tilde{\beta}$ is an extension of $\mathfrak{s}$.

Suppose $\mathfrak{s}^{1}_{\widetilde{M}}$ and $\mathfrak{s}^{2}_{\widetilde{M}}$ are two extensions of $\widetilde{\mathfrak{s}}$. Then $i^{*}(\mathfrak{s}^{1}_{\widetilde{M}}-\mathfrak{s}^{2}_{\widetilde{M}})=0$. By (\ref{eq: short exact sequence 1}), there exists some complex line bundle  $L$ over $\Q$ such that $\mathfrak{s}^{1}_{\widetilde{M}}$ and $\mathfrak{s}^{2}_{\widetilde{M}}$ are related by a twist with $p^{*}(L)$.
\end{proof}

\begin{pro}\label{pro: FSW for ADE} Let $\widetilde{M}/\Q$ be an ADE family. Let $\mathfrak{s}$ be a  $\spinc$ structure on the fiber $M$ that satisfies $d(\mathfrak{s},M)=0$ and $ c_{1}(\mathfrak{s}|_{M_{1}})=0$. 
Let $\widetilde{\mathfrak{s}}$ be any family $\spinc$ structure on $\widetilde{M}/\Q$ such that $d(\widetilde{\mathfrak{s}},\widetilde{M})=0$ and $\widetilde{\mathfrak{s}}|_{M_{2}}=\mathfrak{s}|_{M_{2}}$. Then we have 
$$
\FSW(\widetilde{M}/\Q,\widetilde{\mathfrak{s}})=\pm \SW(M,\mathfrak{s})
$$ when $b^{+}_{2}(M)>1$ and we have  $$
\FSW_{+}(\widetilde{M}/\Q,\widetilde{\mathfrak{s}})=\pm  \SW(M,\mathfrak{s}), \quad
\FSW_{-}(\widetilde{M}/\Q, \widetilde{\mathfrak{s}})=\pm  \SW_{-}(M,\mathfrak{s}) 
$$
when $b^{+}_{2}(M)=1$.
\end{pro}
\begin{proof}
We use $\mathfrak{s}'_{M}$ to denote $\widetilde{\mathfrak{s}}|_{M}$. Since $\mathfrak{s}'_{M}|_{M_{2}}=\mathfrak{s}|_{M_{2}}$, we have $\mathfrak{s}'_{M}=\mathfrak{s}+\operatorname{PD}[S_{r}]$ for some $r\in \mathbb{Z}^{n}$. This implies that 
$$
c_{1}(\widetilde{\mathfrak{s}}|_{\widetilde{M}_{1}})=2\operatorname{PD}[\widetilde{S}_{r}]+p^{*}(\alpha)
$$ 
for some  $\alpha\in H^{2}(\Q;\mathbb{Z})$. Furthermore, we have 
$$
[S_{r}]\cdot 
[S_{r}]= [S_{r}]\cdot [S_{r}]+\langle c_{1}(\mathfrak{s}), [S_{r}]\rangle =d(\mathfrak{s}'_{M},M)-d(\mathfrak{s},M)=-2.
$$
So $r$ is a root.  

By twisting $\widetilde{\mathfrak{s}}$ with a line bundle pulled back from $\Q$ (which doesn't change the family Seiberg-Witten invariant according to Lemma \ref{lem: twisting by a line bundle doesn't change FSW}), we may assume that $\widetilde{\mathfrak{s}}|_{\widetilde{M_{2}}}$ is a pull back of $\mathfrak{s}|_{M_{2}}$. This implies $\alpha=0$ and
$$
c_{1}(\widetilde{\mathfrak{s}}|_{\widetilde{M}_{1}})=2\operatorname{PD}[\widetilde{S}_{r}].$$
Suppose $b^{+}_{2}(M)>1$, then by Theorem \ref{thm: FSW=SW}, we have
    $$\FSW(\widetilde{M},\widetilde{\mathfrak{s}}_{\widetilde{M}})=-\langle c_{1}(\operatorname{Ind}\slashed{D}^{+}(\widetilde{M},\widetilde{\mathfrak{s}})),[\Q]\rangle \cdot \SW(M,\mathfrak{s}).
$$
By Proposition \ref{pro: calculation of index}, we have $\langle c_{1}(\operatorname{Ind}\slashed{D}^{+}(\widetilde{M},\widetilde{\mathfrak{s}})),[\Q]\rangle=\pm 1$. This finishes the proof when $b^{+}_{2}(M)>1$. The argument when $b^{+}_{2}(M)=1$ is the same.
\end{proof}

\begin{pro}\label{pro: FSW for blown-up} Let $\widetilde{M}/\Q$ be a blown-up family. Let $\mathfrak{s}$ be a $\spinc$ structure on the fiber $M$ that satisfies $d(\mathfrak{s},M)=0$ and $ \langle c_{1}(\mathfrak{s}),[S]\rangle =1$.
Let $\widetilde{\mathfrak{s}}$ be any family $\spinc$ structure such that $d(\widetilde{\mathfrak{s}},\widetilde{M})=0$ and $\widetilde{\mathfrak{s}}|_{M_{2}}=\mathfrak{s}|_{M_{2}}$. Then we have 
$$
\FSW(\widetilde{M}/\Q,\widetilde{\mathfrak{s}})= \SW(M,\mathfrak{s})
$$ when $b^{+}_{2}(M)>1$ and we have $$
\FSW_{\pm}(\widetilde{M}/\Q,\widetilde{\mathfrak{s}})=\SW_{\pm}(M,\mathfrak{s}).
$$ when $b^{+}_{2}(M)=1$.
\end{pro}
\begin{proof}
The proof is identical to Proposition \ref{pro: FSW for ADE}, except that we use Proposition \ref{lem: index calculation 2} instead of Lemma \ref{pro: calculation of index}.
\end{proof}

Now we suppose the fiber $M$ has a symplectic structure $\omega$ with the canonical class  $K\in H^{2}(M;\mathbb{Z})$ and the canonical $\spinc$ structure $\mathfrak{s}_{J}$. By Lemma \ref{lem: extension of spinc structure}, there exists a unique family $\spinc$ structure that restricts to $\mathfrak{s}_{J}$ on the fiber $M$ and restricts to a pull-back of $\mathfrak{s}_{J}|_{M_{2}}$ on the product piece $\widetilde{M}_{2}=\Q\times M_{2}$. We denote this family $\spinc$ stucture by $\widetilde{\mathfrak{s}}_{J}$.

\begin{proof}[Proof of Theorem \ref{thm: ADE family doesn't carry symplectic structure}] Suppose (1) is not true and such $a$ exists. Then $a=[S_{r}]$ for some root $r$. By replacing $r$ with $-r$, we may assume that 
\begin{equation}\label{eq: symplectic form all negative}
\langle [\omega_{b}], [S_{r}]\rangle \leq0,\quad \forall b\in \Q.    
\end{equation}
Consider the family $\spinc$ structure $\widetilde{\mathfrak{s}}=\widetilde{\mathfrak{s}}_{J}+\operatorname{PD}[\widetilde{S}_{r}]$. Then
By Theorem \ref{pro: family vanishing}, we must have \begin{equation}\label{eq: FSW=0 for symplectic}
\FSW_{\xi_{\omega}}(\widetilde{M},\widetilde{\mathfrak{s}})=0.
\end{equation}
However, Proposition \ref{pro: FSW for ADE} and Theorem \ref{thm: taubes's vanishing} together imply  that 
$$
\FSW_{\xi_{c}}(\widetilde{M},\widetilde{\mathfrak{s}})=\pm \SW(M,\mathfrak{s}_{J})=\pm 1.
$$
So the only possibility is that $\xi_{\omega}\neq \xi_{c}$. This can only happen when $b^{+}_{2}(M)=3$. By (\ref{eq: wall-crossing}), $\xi_{\omega}$ must be the winding chamber. That implies the map (\ref{eq: characteristic map}) has degree $\pm 1$. Since $M$ is minimal, by a result of Taubes \cite[Theorem 0.2]{Taubes1996}, we have $K\cdot K\geq 0$ and $[\omega_{b}]\cdot K\geq 0$ for any $b\in \Q$. By Lemma \ref{lem: topology of positive cone}(ii), $K$ must be torsion.
\end{proof}

\begin{proof}[Proof of Theorem \ref{thm: blown-up family has fiberwise symplectic structure}]
By Theorem \ref{thm: Taubes2} and Theorem \ref{thm: Li-Liu}, we have $\langle [\omega_{b}],[S]\rangle >0$ for any $b\in \Q$. 
Consider the family $\spinc$ structure
$$
\widetilde{\mathfrak{s}}=\widetilde{\mathfrak{s}}_{J}-\operatorname{PD}[\widetilde{S}],
$$
Then by Proposition \ref{pro: family vanishing}, we have 
$$
\FSW_{\xi_{\omega}}(\widetilde{M}/\Q,\widetilde{\mathfrak{s}})=0.
$$
Again, we compare this with Theorem \ref{thm: taubes's vanishing} and Proposition \ref{pro: FSW for blown-up}. We see that the only possibility is that $b^{+}_{2}(M)=3$ and the map (\ref{eq: characteristic map}) has mapping degree $\pm 1$.
\end{proof}

We end this section by the following proposition, which used in Section \ref{section: linearly independence} to show linearly independence of generalized Dehn twists.

\begin{pro}\label{pro: FSW linearly indep} Let $M$ be a smooth 4-manifold with two smoothly embedded spheres $S_{0}$ and $S_{1}$ of self-intersection $-1$ or $-2$. Let $\mathfrak{s}$ be a $\spinc$ structure with $d(\mathfrak{s},M)=0$ such that the adjunction formula 
$$
\langle c_{1}(\mathfrak{s}), [S_{i}]\rangle =S_{i}\cdot S_{i}+2
$$
holds for $i=1,2$. When $S_{0}\cdot S_{0}=-1$ (resp. $S_{0}\cdot S_{0}=-2$), we let $\widetilde{M}/\Q$ be the blown-up family (resp. the $A_{1}$-family) associated to $S_{0}$. Let $\widetilde{\mathfrak{s}}_{1}$ be a family $\spinc$ structure on $\widetilde{M}/\Q$ that restricts to $\mathfrak{s}-\operatorname{PD}(S_{1})$ on the fiber. Then the following conclusion holds.
\begin{enumerate}
    \item Suppose $b^{+}_{2}(M)>1$ and assume $$
\SW(M,\mathfrak{s}-\operatorname{PD}(S_{1})-\operatorname{PD}(S_{0}))=\SW(M,\mathfrak{s}-\operatorname{PD}(S_{1})+\operatorname{PD}(S_{0}))=0.
$$
Then we have 
\begin{equation}\label{eq: FSW vanishing}
\FSW(\widetilde{M}/
\Q,\widetilde{\mathfrak{s}}_{1})=0.    
\end{equation}
\item Suppose $b^{+}_{2}(M)=1$ and assume $$
\SW_{-}(M,\mathfrak{s}-\operatorname{PD}(S_{1})-\operatorname{PD}(S_{0}))=\SW_{-}(M,\mathfrak{s}-\operatorname{PD}(S_{1})+\operatorname{PD}(S_{0}))=0.
$$
Then we have 
\begin{equation}\label{eq: FSW vanishing2}
\FSW_{-}(\widetilde{M}/
\Q,\widetilde{\mathfrak{s}}_{1})=0. \end{equation}
Same result holds for $\SW_{+}$ and $\FSW_{+}$.
\end{enumerate}

\end{pro}

\begin{lem}\label{lem: index calculation}
Let $M, \widetilde{M}, S_{0}, S_{1}, \widetilde{\mathfrak{s}}_{1}$ be as in Proposition \ref{pro: FSW linearly indep}. Let $k=S_{0}\cdot S_{0}$ and let $n=k+2-2S_{0}\cdot S_{1}$.  Suppose $ \widetilde{\mathfrak{s}}_{1}|_{\widetilde{M}_2}$ is pull-back from $M_{2}$. 
Then we have 
$$
\langle c_{1}(\operatorname{Ind}(\slashed{D}^{+}(\widetilde{M},\widetilde{\mathfrak{s}}_{1}))),[\Q]\rangle=\frac{n(n-k)(n+k)}{24k}.
$$
\end{lem}
\begin{proof}
We first assume $k=-2$. By the Atiyah-Singer index theorem, we have 

$$
\langle c_{1}(\operatorname{Ind}(\slashed{D}^{+}(\widetilde{M},\widetilde{\mathfrak{s}}_{1}))),[\Q]\rangle
=\frac{\langle c_{1}(\widetilde{\mathfrak{s}}_{1})^{3}-p_{1}(T\widetilde{M}/\Q)\cdot c_{1}(\widetilde{\mathfrak{s}}_{1},[\widetilde{M}]\rangle}{48}.
$$ 
Consider $c_{1}(\widetilde{\mathfrak{s}}_{1}|_{\widetilde{M}_{1}})\in H^{2}(\widetilde{M}_{1};\mathbb{R})$. Since $\widetilde{\mathfrak{s}}_{1}|_{\widetilde{M}_{1}}$ is pulled back from $\partial \widetilde{M}_{1}$ and $H^{2}(\partial \widetilde{M}_{1};\mathbb{R})=0$, there exists $\widetilde{x}_{1}\in H^{2}(\widetilde{M}_{1},\partial \widetilde{M}_{1};\mathbb{R})$ that is sent to $c_{1}(\widetilde{\mathfrak{s}}_{1}|_{\widetilde{M}_{1}})$. Smilarly, there exists $x_{2}\in H^{2}(M_{2},\partial M_{2} ;\mathbb{R})$ that is sent to $c_{1}(\widetilde{\mathfrak{s}}_{1}|_{\widetilde{M}_{2}})\in H^{2}(M_{2};\mathbb{R})$. We let $\widetilde{x}_{2} \in H^{2}(\widetilde{M}_{2},\partial \widetilde{M}_{2} ;\mathbb{R})$ be the pull-back of $x_{2}$ under the projection map $\widetilde{M}_{2}\rightarrow M_{2}$. Then we have 
\[\begin{split}
\langle c_{1}(\operatorname{Ind}(\slashed{D}^{+}(\widetilde{M},\widetilde{\mathfrak{s}}_{1}))),[\Q]\rangle
&=\frac{\langle \widetilde{x}^{3}_{1}-p_{1}(T\widetilde{M}_{1}/\Q)\cdot \widetilde{x}_{1},[\widetilde{M}_{1}]\rangle}{48}+\frac{\langle \widetilde{x}^{3}_{2}-p_{1}(T\widetilde{M}_{2}/\Q)\cdot \widetilde{x}_{2},[\widetilde{M}_{2}]\rangle}{48}\\
&=\frac{\langle \widetilde{x}^{3}_{1}-p_{1}(T\widetilde{M}_{1})/\Q\cdot \widetilde{x}_{1}),[\widetilde{M}_{1}]\rangle}{48}
\end{split}
\] 
where the second equality follows from the observation that both $\widetilde{x}_{2}$ and $p_{1}(T\widetilde{M}_{2}/\Q)$ are pulled back from $M_{2}$. Since $H_{4}(\widetilde{M}_{1};\mathbb{R})$ is generated by $\widetilde{S}_{0}$, we can have 
$$
\widetilde{x}_{1}=a\cdot \operatorname{PD}(\widetilde{S}_{0})
$$
for some $a\in \mathbb{R}$. Then by Lemma \ref{lem: self-intersection}, we have 
$$
\frac{\langle \widetilde{x}^{3}_{1}-p_{1}(T\widetilde{M}_{1}/\Q)\cdot \widetilde{x}_{1},[\widetilde{M}_{1}]\rangle}{48}=\frac{a^{3}-a}{6}.
$$
Because $\widetilde{S}_{0}\cdot S_{0}=-2$ and $\langle c_{1}(\widetilde{\mathfrak{s}}_{1}), [S_{0}]\rangle =n$, we get $a=\frac{n}{2}$. This finishes the proof when $k=-2$. The case $k=-1$ is similar. 
\end{proof}

\begin{proof}[Proof of Proposition \ref{pro: FSW linearly indep}]
We first assume $b^+(M)>1$. By twisting $\widetilde{\mathfrak{s}}_{1}$ with a line bundle $L$ pulled back from $\Q$, we may assume that $\widetilde{\mathfrak{s}}_{1}|_{\Q\times M_{2}}$ is pulled back from $M_{2}$. 
We may further assume that \begin{equation}\label{eq: FSW linearly indepentent eq 1}
    \frac{n(n-k)(n+k)}{24k}\neq 0
\end{equation}
because otherwise Lemma \ref{lem: index calculation},  
 Theorem \ref{thm: FSW gluing} and Theorem \ref{thm: computing cobordism induced map} will directly imply (\ref{eq: FSW vanishing}).

As the next step, we apply the switching formula (\ref{eq: FSW=SW 1}) to conclude that 
\begin{equation}\label{eq: FSW linearly independent eq2}
\FSW(\widetilde{M}_{0}/
\Q,\widetilde{\mathfrak{s}}_{1})=-\frac{n(n-k)(n+k)}{24k}\cdot \SW(M,\mathfrak{s}-\operatorname{PD}(S_{1})+l\operatorname{PD}(S_{0}))
\end{equation}
for any integer $l$ that satisfies 
$$
d(\mathfrak{s}-\operatorname{PD}(S_{1})+l\operatorname{PD}(S_{0}),M)\geq 0.
$$
We denote this expected dimension by $d(l)$ and compute it as follows:
\begin{equation*}
\begin{split}
d(l)=&d(\mathfrak{s},M)+(l\operatorname{PD}(S_{0})-\operatorname{PD}(S_{1}))^{2}+c_{1}(\mathfrak{s})\cdot (l\operatorname{PD}(S_{0})-\operatorname{PD}(S_{1}))\\
=& l^{2}S_{0}\cdot S_{0} +S_{1}\cdot S_{1}-2lS_{0}\cdot S_{1}+l \langle c_{1}(\mathfrak{s}),[S_{0}]\rangle - \langle c_{1}(\mathfrak{s}),[S_{1}]\rangle\\
=& kl^{2}+nl-2.
\end{split}
\end{equation*}
To finish the proof, we consider four possible cases. Note that $k, n$ have the same parity.
\begin{enumerate}[(i)]
    \item $k=-1$ and $|n|>1$: We let $l=1$ if $n>0$ and let $l=-1$ if $n<0$. Then 
$$
d(l)=|n|-3\geq 0.
$$
 By our assumption, we have $$\SW(M,\mathfrak{s}-\operatorname{PD}(S_{1})\pm  \operatorname{PD}(S_{0}))=0$$ and (\ref{eq: FSW vanishing}) follows from (\ref{eq: FSW linearly independent eq2}). 
\item $k=-1$ and $n=\pm 1$: This is ruled out by (\ref{eq: FSW linearly indepentent eq 1}). 
    \item $k=-2$ and $|n|\geq 4$: We let $l=1$ if $n>0$ and let $l=-1$ if $n<0$. Then $d(l)=|n|-4\geq 0$ and the argument is the same as (i).
    \item $k=-2$ and $n=0,\pm 2$: This is ruled out by (\ref{eq: FSW linearly indepentent eq 1}). 
\end{enumerate}
This finishes the proof when $b^{+}_{2}(M)>1$. The proof when $b^{+}_{2}(M)=1$ is completely analogous.
\end{proof}

\section{Linearly independence of the generalized Dehn twists}\label{section: linearly independence}
In this section, we prove Theorem \ref{thm: infinitely generated diff}, Theorem \ref{thm: linearly independence} and various corollaries of them. The ideas in these proofs are similar: To show a loop in $\operatorname{Diff}(M)$ is noncontractible, we use the loop as a clutching function and show that the corresponding family over $S^2$ has nontrivial family Seiberg-Witten invariants. To show that this loop is not homotopic to a loop in $\operatorname{Symp}(X)$, we use our results in previous sections to conclude that this family do not support a fiberwise symplectic structure (in a fixed cohomology class). 
\subsection{The symplectic case.}

Let $(M,\omega)$ be a symplectic $4$-manifold 
with canonical class $K$. We fix a compatible almost complex structure $J$.

 We assume that the set $R_{X}$ (as defined in (\ref{eq: the set of roots})) is nonempty. For each $\alpha\in R_{X}$, we pick an embedded sphere $S_{\alpha}$ representing the corresponding homology class. 

\begin{lem}
One can choose a suitable orientation on $S_{\alpha}$ so that 
\begin{equation}\label{eq: suitable orientation}
\int_{S_{\alpha}}\omega\geq 0,\quad -\langle K, [S_{\alpha}]\rangle =2+S_{\alpha}\cdot S_{\alpha}.    
\end{equation}
\end{lem}
\begin{proof}
By the definition of $R_{M}$, the self-intersection $S_{\alpha}\cdot S_{\alpha}$ can be $-1$ or $-2$. 

When $S_{\alpha}\cdot S_{\alpha}=-1$, we orient $S_{\alpha}$ so that $\langle K, [S_{\alpha}]\rangle=-1$. Then by Theorem \ref{thm: Taubes2} and Theorem \ref{thm: Li-Liu}, $S_{\alpha}$ is homologous to an embedded holnomorphic curve $C$. Therefore, (\ref{eq: suitable orientation}) is satisfied.

When $S_{\alpha}\cdot S_{\alpha}=-2$, we have $\langle K, [S_{\alpha}]\rangle=0$ so the equality $-\langle K, [S_{\alpha}]\rangle=2+S_{\alpha}\cdot S_{\alpha}$ holds for both orientations. By choosing a suitable orientation, we may assume $\int_{S_{\alpha}}\omega\geq 0$
\end{proof}
From now on, we assume that (\ref{eq: suitable orientation}) is satisfied for all $S_{\alpha}$. For each $\alpha\in R_{M}$, we use the generalized Dehn twist $\gamma_{S_{\alpha}}$ as the clutching function and form a family $\widetilde{M}_{\alpha}/\Q$ over $\Q=S^{2}$. When $S_{\alpha}\cdot S_{\alpha}=-1$,  $\widetilde{M}_{\alpha}/
\Q$ is a blown-up family. When $S_{\alpha}\cdot S_{\alpha}=-2$, $\widetilde{M}_{\alpha}/
\Q$ is an $A_{1}$ family. 
For any formal linear combination
$
\vec{\alpha}=\sum_{i=1}^{n}a_{i}\alpha_{i}\in \mathbb{Z}^{R_{M}}
$, we consider the loop 
$$ 
\gamma_{\vec{\alpha}}:=\gamma^{a_{1}}_{S_{\alpha_{1}}}\cdots \gamma^{a_{n}}_{S_{\alpha_{n}}}$$ in $\operatorname{Diff}(M)$ and we use it as the clutching function to define the family 
$\widetilde{M}_{\vec{\alpha}}/\Q$.

\begin{thm}\label{thm: linearly independence restated}  Consider a closed,  oriented surface $\Q'$ and a smooth map $f:\Q'\rightarrow \Q$ of nonzero degree. Let $\widetilde{M}'_{\vec{\alpha}}/\Q'$ be the smooth family obtained by pulling back $\widetilde{M}_{\vec{\alpha}}/\Q$ via $f$. Suppose $\vec{\alpha}\neq \vec{0}$. Then $\widetilde{M}_{\vec{\alpha}}'/\Q'$ does not carry a fiberwise symplectic structure $\{\omega_{b}\}$ in the fixed cohomology class $[\omega]$. In particular, $\widetilde{M}_{\vec{\alpha}}/\Q$ itself doesn't carry a fiberwise symplectic structure in the fixed cohomology class $[\omega]$.
\end{thm}
\begin{proof} We first assume $b^{+}(M)>1$. Let $\vec{\alpha}=a_{1}\alpha_{1}+\cdots +a_{n}\alpha_{n}$. Without loss of generality, we may assume that 
$$
\int_{S_{\alpha_{1}}}\omega \geq \int_{S_{\alpha_{2}}}\omega,\cdots, \int_{S_{\alpha_{n}}}\omega\geq 0.
$$
Then by Taubes' theorem Theorem \ref{thm: taubes's vanishing}, we have 
$$
\SW(M,\mathfrak{s}_{J}-\operatorname{PD}(S_{\alpha_1})\pm \operatorname{PD}(S_{\alpha_i}))=0\quad \text{ for any } 2\leq i\leq n.
$$
Next, we consider a family $\spinc$ structure $\mathfrak{s}_{\widetilde{M}_{\vec{\alpha}}}$ such that 
\begin{equation*}
\mathfrak{s}_{\widetilde{M}_{\vec{\alpha}}}|_{M}=\mathfrak{s}_{J}-\operatorname{PD}(S_{\alpha_{1}}).
\end{equation*}
By the additivity of the family Seiberg-Witten invariant, we have 
$$
\FSW(\widetilde{M}_{\vec{\alpha}}/\Q,  \mathfrak{s}_{\widetilde{M}_{\vec{\alpha}}})=\sum^{n}_{i=1}a_{i}\cdot \FSW(\widetilde{M}_{\alpha_{i}}/\Q,  \mathfrak{s}_{\widetilde{M}_{\alpha_{i}}}),
$$
where $\mathfrak{s}_{\widetilde{M}_{\alpha_{i}}}$ is a family $\spinc$ structure on $\widetilde{M}_{\alpha_{i}}/\Q$ that restricts to $\mathfrak{s}_{J}-\operatorname{PD}(S_{\alpha_{1}})$ on the fiber $M$. 
By Proposition \ref{pro: FSW linearly indep}, we have $$\FSW(\widetilde{M}_{\alpha_{i}}/\Q,  \mathfrak{s}_{\widetilde{M}_{\alpha_{i}}})=0,\text{ for any }i\geq 2.$$By Proposition \ref{pro: FSW for ADE} and Proposition \ref{pro: FSW for blown-up}, we have $$\FSW(\widetilde{M}_{\alpha_{1}}/\Q,  \mathfrak{s}_{\widetilde{M}_{\alpha_{1}}})=\pm 1.$$
Therefore, we get $$
\FSW(\widetilde{M}_{\vec{\alpha}}/\Q,  \mathfrak{s}_{\widetilde{M}_{\vec{\alpha}}})=\pm a_{1}.$$

Now we consider the family Seiberg-Witten invariant of $\widetilde{M}_{\vec{\alpha}}'/\Q'$. Since $\widetilde{M}_{\vec{\alpha}}'/\Q'$  is pulled back from the homology product $\widetilde{M}_{\alpha}/\Q$, itself is also a homology product.  Let $\mathfrak{s}_{\widetilde{M}'_{\vec{\alpha}}}$ be the pull-back of  $\mathfrak{s}_{\widetilde{M}_{\vec{\alpha}}}$. Then by Lemma \ref{lem: FSW pull-back}, we have 
$$
\FSW(\widetilde{M}'_{\vec{\alpha}}/\Q',  \mathfrak{s}_{\widetilde{M}'_{\vec{\alpha}}})=\operatorname{deg}(f)\cdot \FSW(\widetilde{M}_{\vec{\alpha}}/\Q,  \mathfrak{s}_{\widetilde{M}_{\vec{\alpha}}})=\operatorname{deg}(f)\cdot a_{1}\neq 0
$$

Now suppose $\widetilde{M}'_{\vec{\alpha}}/\Q'$ support a fiberwise symplectic structure $\{\omega_{b}\}$ in the cohomology class $[\omega]$.
By Theorem \ref{thm: taubes's vanishing}, the canonical class of $\omega_{b}$ equals $K$. Therefore, Proposition \ref{pro: family vanishing} tells us that $\FSW(\widetilde{M}'_{\vec{\alpha}}/\Q',  \mathfrak{s}_{\widetilde{M}'_{\vec{\alpha}}})=0$. This leads to the contradiction. 

The case $b^{+}(M)=1$ is completely analogous. One can show that
$$
\FSW_{-}(\widetilde{M}'_{\vec{\alpha}}/\Q,  \mathfrak{s}_{\widetilde{M}'_{\vec{\alpha}}})=\operatorname{deg}(f)\cdot a_{1}\neq 0.
$$
and then apply Proposition \ref{pro: family vanishing}.
\end{proof}


\begin{proof}[Proof of Theorem \ref{thm: linearly independence}] By Theorem \ref{thm: linearly independence restated}, the structure group of $\widetilde{M}_{\vec{\alpha}}/\Q$ can not be reduced to $\Symp(M)$. Hence its clutching function $\gamma_{\vec{\alpha}}$ is not homotopic to a loop in $\Symp(M)$.

\end{proof}

\begin{proof}[Proof of Corollary \ref{thm: non-symplectic loop}] Suppose $M$ is minimal. Then $S\cdot S=-2$. By Lemma \ref{pro: minimal implies adjunction}, we have 
$
-\langle K,[S]\rangle=0
$ and $[S]\in R_{M}$. Suppose $M$ is nonminimal. Then the exceptional curve represents an element in $R_{M}$. In both cases, the set $R_{M}$ is nonempty so the result follows from Theorem \ref{thm: linearly independence}.\end{proof}

\begin{proof}[Proof of Corollary \ref{cor: noncontractible loop of symplectic forms}] (1) For any $\vec{\alpha}\in \mathbb{Z}^{R_{M}}\setminus \{\vec{0}\}$, we consider the loop of symplectic structures 
$
\widetilde{\omega}_{\vec{\alpha}}=\{\omega_{\vec{\alpha}}(t)\}_{t\in S^1}
$
defined by $\widetilde{\omega}_{\vec{\alpha}}(t):=\gamma_{\alpha}(t)^{*}(\omega)$. We will show that $\widetilde{\omega}_{\vec{\alpha}}$ is not null-homologous in the space $\mathbb{S}_{[\omega]}(M)$. Suppose this is not the case. Then there exists a smooth family of cohomologous symplectic forms 
$
\widetilde{\omega}_{1}:=\{\omega_{b}\}_{b\in F}
$ on $M$ that is parameterized by an oriented surface $F$ and restricts to $\widetilde{\omega}_{\vec{\alpha}}$ on $\partial F$. We also consider the constant family of symplectic forms $\widetilde{\omega}_{2}:=\{\omega_{b}\}_{b\in D^{2}}$ with $\omega_{b}=\omega$ for any $b\in D^2$. Let $\Q'=F\cup_{S^{1}}D^2$. Consider the smooth family $\widetilde{M}'_{\vec{\alpha}}/\Q'$ defined by 
$$
\widetilde{M}'_{\vec{\alpha}}=((F\times X)\sqcup (D^2\times X))/\sim,$$
where $\sim$ is generated by $(t,x)\sim (t,\gamma_{\alpha}(t)(x))$ for $t\in \partial D^{2}=\partial F.
$ By gluing together $\widetilde{\omega}_{1}$ and $\widetilde{\omega}_{2}$, we obtain a fiberwise symplectic structure on $\widetilde{M}'_{\vec{\alpha}}/\Q'$ in the cohomology class $[\omega]$. However, $\widetilde{M}'_{\vec{\alpha}}/\Q'$ can also be obtained by pulling back $\widetilde{M}_{\vec{\alpha}}/\Q$ via a degree-$1$ map from $\Q'$ to  $\Q$. So  Theorem \ref{thm: linearly independence restated} implies that such a fiberwise symplectic structure does not exist. This leads to the contradiction.

(2) Now we assume  $\vec{\alpha}\in \mathbb{Z}^{R^{2}_{M}}\setminus \{\vec{0}\}$ and prove that $\widetilde{\omega}_{\vec{\alpha}}$ is indeed null-homotopic (and hence null-homologous) in $\mathcal{S}(M)$ (the space of almost symplectic forms). By linearity, it suffices to focus on the case $\vec{\alpha}=\alpha$ for some $\alpha\in R^{2}_{M}$. 
By Lemma \ref{almost complex structure}, the family $\widetilde{M}_{\vec{\alpha}}/\Q$ admits a fiberwise almost complex structure $\widetilde{J}=\{J_{b}\}_{b\in \Q}$ such that $J_{b_0}=J$. By choosing compatible matrices on the fibers, $\widetilde{J}$ induces a smooth family of almost symplectic forms $\widetilde{\omega}_{4}=\{\omega_{b}\}_{b\in \Q}$ such that $\omega_{b_0}=\omega$. By pulling back  $\widetilde{\omega}_{4}$ via the quotient map
$$
D^{2}\rightarrow D^{2}/\partial D^{2}\cong \Q,
$$
we obtain a smooth family $\widetilde{\omega}_{5}$ of almost symplectic forms that is parametrzied by $D^2$  and restricts to $\widetilde{\omega}_{\vec{\alpha}}$ on $\partial D^2$. So $\widetilde{\omega}_{\vec{\alpha}}$ is null-homotopic in $\mathcal{S}(M)$ when $\vec{\alpha}=\alpha$.

(3) We first assume $b^{+}_{2}(M)>1$. Consider a nontrivial linear combination $
\vec{\alpha}=\sum_{i=1}^{n}a_{i}\alpha_{i}
$
of elements $\alpha_{i}\in R^{1}_{M}$. We may assume
$$
\int_{S_{\alpha_{1}}}\omega \geq \int_{S_{\alpha_{2}}}\omega,\cdots, \int_{S_{\alpha_{n}}}\omega\geq 0.
$$
Suppose the loop $\widetilde{\omega}_{\vec{\alpha}}$ is null-homologous in $\mathbb{S}(M)$. Then smilar to (1), one can construct a smooth family $\widetilde{M}'_{\vec{\alpha}}/\Q'$, parametrized by a surface $\Q'$, such that
it admits a fiberwise symplectic structure  $\{\omega_{b}\}$ with $\omega_{b_0}=\omega$ and that 
$$
\FSW(\widetilde{M}'_{\vec{\alpha}}/\Q',  \mathfrak{s}_{\widetilde{M}'_{\vec{\alpha}}})\neq 0.
$$
for some family $\spinc$-structure $\mathfrak{s}_{\widetilde{M}'_{\vec{\alpha}}}$ that restricts to $\mathfrak{s}_{J}-\operatorname{PD}(S_{\alpha_{1}})$ on the fibers. By Theorem \ref{thm: Taubes2} and Theorem \ref{thm: Li-Liu}, we have $\langle [\omega_{b}],[S_{\alpha}]\rangle >0$ for any $b\in \Q'$. So Proposition \ref{pro: family vanishing} tells us that 
$$
\FSW_{\xi_{\omega}}(\widetilde{M}'_{\vec{\alpha}}/\Q',  \mathfrak{s}_{\widetilde{M}'_{\vec{\alpha}}})= 0.
$$
Since $b^{+}(M)\neq 2, 3$. Every map from $\Q'$ to $V^{+}(M)\simeq S^{b^{+}(M)-1}$ is null-homotopic. So the symplectic chamber $\xi_{\omega}$ is just the canonical chamber. This leads to the contradiction. 

When $b^{+}(M)=1$, we use $\FSW_{-}$ instead of $\FSW$ and then repeat the argument. 

(4) The proof is the same as (3), except that we can set $\Q'=\Q$ in this case. So every map $\Q'$ to $V^{+}(M)\simeq S^{b^{+}(M)-1}$ is null-homotopic as long as $b^{+}(M)\neq 3$.
\end{proof}

\subsection{The smooth case} From now on, we assume that $M$ is a smooth (but not necessarily symplectic) $4$-manifold with $b_{1}(M)=0$. 
\begin{lem}\label{lem: relfecting spinc} Suppose $b^{+}_{2}(M)>1$. Consider the set $$
\mathcal{A}_{M}:=\{\mathfrak{s}\mid \SW(M,\mathfrak{s})\neq 0\}
$$
of Seiberg-Witten basic classes. Then for any smoothly embedded  $-2$-sphere $S$ and any $\mathfrak{s}\in \mathcal{A}_{M}$, we have \begin{equation}\label{eq: relfecting spinc}
\mathfrak{s}+\frac{\langle c_{1}(\mathfrak{s}),[S]\rangle }{2} \cdot S\in \mathcal{A}_{M}
\end{equation}
\end{lem}

\begin{proof}
Given any such $S$, we consider the relfection diffeomorphism
$
r: M\rightarrow M
$. (See the proof of Lemma \ref{pro: minimal implies adjunction} for an explicit description.)   Since $r$ is supported near $S$, we have 
$$
r_{*}(\mathfrak{s})=\mathfrak{s}+k\operatorname{PD}(S)
$$
for some $k$. To compute $k$, we note that the induced action
$$r^{*}:H^{2}(M;\mathbb{Z})\rightarrow H^{2}(M;\mathbb{Z})$$ is the relation along $\operatorname{PD}(S)$. In particular, one has 
$$
\langle c_{1}(r^{*}(\mathfrak{s})),[S]\rangle =-\langle c_{1}(\mathfrak{s}),[S]\rangle
$$
This implies $k=\frac{\langle c_{1}(\mathfrak{s}),[S]\rangle }{2}$. Since $\mathcal{A}_{M}$ is preserved by any orientation preserving diffeomorphism, we have proved (\ref{eq: relfecting spinc}).
\end{proof}
Now we are ready to prove Theorem \ref{thm: infinitely generated diff}, restated as follows.

\begin{thm}
Let $M$ be a smooth closed 4-manifold that contains an infinite collection of smoothly embedded $(-2)$-spheres $\{S_{i}\}_{i\in \mathbb{N}}$, representing different homology classes. Suppose  $b_{1}(M)=0$ and assume one of the following conditions holds:
\begin{enumerate}
    \item $b^{+}_{2}(M)>1$ and   $\SW(M,\mathfrak{s}_M)\neq 0$ for some $\spinc$-structure $\mathfrak{s}$ with $d(\mathfrak{s},M)=0$. 
    \item $b^{+}_{2}(M)=1$ and there exists a $\spinc$ structure $\mathfrak{s}$ such that $d(\mathfrak{s},M)=0$ and that 
    $
    \langle c_{1}(\mathfrak{s}),[S_{i}]\rangle =0 
    $ for any $i$. (We don't impose any constraint on the Seiberg-Witten invariant in this case.)
\end{enumerate}
Then $\pi_{1}(\operatorname{Diff}(M))$ contains a $\mathbb{Z}^{\infty}$-summand, generated by generalized Dehn twists along a subcollection $\{S_{n_{i}}\}_{i\in \mathbb{N}}$.
\end{thm}

\begin{proof} (1) Since $\mathcal{A}_{M}$ is a finite set, after passing to an infinite subcollection, we may assume that for any $\mathfrak{s}\in \mathcal{A}_{M}  $, any $i\neq j$ and any $k\neq 0$, we have 
\begin{equation}\label{eq: infinite summand eq1}
\mathfrak{s}+k\operatorname{PD}(S_{i})\notin  \mathcal{A}_{M}   \end{equation}
and 
\begin{equation}\label{eq: infinite summand eq2}
\mathfrak{s}-\operatorname{PD}(S_{i})\pm \operatorname{PD}(S_{j})\notin  \mathcal{A}_{M}.   
\end{equation}

By (\ref{eq: infinite summand eq1}) and Lemma \ref{lem: relfecting spinc}, we have 
$\langle c_{1}(\mathfrak{s}), [S_{i}]\rangle =0$ for any $i\in \mathbb{N}$. Denote $\mathfrak{s}-\operatorname{PD}(S_{i})$ by $\mathfrak{s}_{i}$. Then $d(\mathfrak{s}_{i},M)=-2$. By Proposition \ref{pro: finiteness b+>1}, we have a group homomorphism 
$$
\mathop{\bigoplus}_{i\in \mathbb{N}}\SW(-,\mathfrak{s}_{i}):\pi_{1}(\Diff(M))\rightarrow \mathbb{Z}^{\infty}.
$$
It suffices to show that this homomorphism is surjective. To see this, we consider the generalized Dehn twist $\gamma_{S_{i}}$ around $S_{i}$. Then 
Proposition \ref{pro: FSW for ADE}, Proposition \ref{pro: FSW linearly indep} and (\ref{eq: infinite summand eq2})
together imply that 
$$
\SW(\gamma_{S_{i}},\mathfrak{s}_{j})=\begin{cases}
\pm \SW(M,\mathfrak{s})\neq 0 &\text{ if }i=j\\
0 &\text{ if }i\neq j
\end{cases}
$$ This finishes the proof.

(2) We pick $H\in H^{2}(M;\mathbb{Z})$ with $H\cdot H>0$ and use it to fix a homology orientation on $M$. Let $$\mathcal{A}_{M}:=\{\mathfrak{s}\mid \SW_{H}(M,\mathfrak{s})\neq 0\},
$$
where $\SW_{H}(M,\mathfrak{s})$ is the small perturbation Seiberg-Witten invariant defined as 
$$
\SW_{H}(M,\mathfrak{s}):= \begin{cases} \SW_{+}(M,\mathfrak{s})\quad &\text{if }c_{1}(\mathfrak{s})\cdot H\geq 0\\ 
\SW_{-}(M,\mathfrak{s}) \quad &\text{if }c_{1}(\mathfrak{s})\cdot H< 0
\end{cases}. 
$$
Then compactness of the Seiberg-Witten moduli space implies that $\mathcal{A}_{M}$ is a finite set. Since $\SW_{+}(M,\mathfrak{s})$ and $\SW_{-}(M,\mathfrak{s})$ differ by $\pm 1$, at least one of them is nonzero. By changing the sign of $H$ if necessary, we may assume $\SW_{-}(M,\mathfrak{s})\neq 0$. 

Since $S_{i}$ all represent distinct homology classes of self-intersection $-2$, it's straightforward to see that $$\lim\limits_{i\rightarrow +\infty}|\operatorname{PD}(S_{i})\cdot H|= +\infty.$$ Therefore, by passing to an infinite subcollection and choosing suitable orientations on $S_{i}$, we may assume that the following conditions all hold.
\begin{itemize}
    \item For any $i\in \mathbb{N}$, we have $
c_{1}(\mathfrak{s}-\operatorname{PD}(S_{i}))\cdot H<0.
$. 
   \item For any $j>i$, we have 
 \begin{equation}\label{eq: infinite summand 1}
c_{1}(\mathfrak{s}-\operatorname{PD}(S_{j})\pm \operatorname{PD}(S_{i}))\cdot H<0.
 \end{equation}
 \item For any $j> i$, we have  
 \begin{equation}\label{eq: infinite summand 2}
 \mathfrak{s}-\operatorname{PD}(S_{j})\pm \operatorname{PD}(S_{i})\notin \mathcal{A}_{M}.
\end{equation}

\end{itemize}
We let $\mathfrak{s}_{i}=\mathfrak{s}-\operatorname{PD}(S_{i})$. By Proposition \ref{pro: finiteness b+=1}, we have a group homomorphism 
\begin{equation}\label{eq: homomorphism to Z-infinity}
\mathop{\bigoplus}_{i\in \mathbb{N}}\SW_{H}(-,\mathfrak{s}_{i}):\pi_{1}(\Diff(M))\rightarrow \mathbb{Z}^{\infty}.
\end{equation}

To show that this homomorphism is surjective, we again consider the generalized Dehn twist $\gamma_{S_{i}}$ around $S_{i}$. By Proposition \ref{pro: FSW for ADE}, we have 
$$
\SW_{H}(\gamma_{S_{i}},\mathfrak{s}_{i})=\SW_{-}(M,\mathfrak{s})\neq 0.
$$
While for any $i<j$, we apply Proposition \ref{pro: FSW linearly indep}, (\ref{eq: infinite summand eq1}) and (\ref{eq: infinite summand eq2}) to conclude that 
$$
\SW_{H}(\gamma_{S_{i}},\mathfrak{s}_{j})=0.
$$
Therefore, the homomorphism (\ref{eq: homomorphism to Z-infinity}) is surjective. 
\end{proof}

\begin{proof}[Proof of Corollary \ref{cor: construction of examples}] (1) We start with the case $M=E(1)$, obtained by blowing up $\mathbb{CP}^{2}$ along the nine intersection points of two generic cubical curves. We set $\mathfrak{s}$ to be the canonical $\spinc$-structure $\mathfrak{s}_{J}$. Then
$$
c_{1}(\mathfrak{s})=-K=\operatorname{PD}(3[H]-\sum_{i=1}^{9}[C_{i}])=\operatorname{PD}[F_{\infty}].
$$
Here $H$ is a hyperplane in $\mathbb{CP}^{2}$, $C_{i}$ is the $i$-th exceptional divisor and $F_{\infty}$ is a regular fiber. We consider the map
$$
\rho: H_{2}(M\setminus F_{\infty};\mathbb{Z})\rightarrow H_{2}(M;\mathbb{Z}).
$$
Then the $\operatorname{Image}(\rho)$ equals 
$$
K^{\perp}:=\{\alpha \in H(M;\mathbb{Z})\mid \langle K, \alpha \rangle =0\}.
$$
Take an infinite collection of classes in $K^{\perp}$ with square $-2$. For example, we can set 
$$
\alpha_{i} = [C_{1}]-[C_{2}]+i[F_{\infty}]\quad\text{for }i\in \mathbb{N}.
$$
Then by \cite[Theorem 5.1]{FriedmanMorgan}, $\alpha_{i}$ can be represented by an embedded sphere in $S_{i}\hookrightarrow M\setminus F_{\infty}$. This finishes the proof when $M=E(1)$.

For a generic elliptic surface $M$ with $b_{1}(M)=0$, it can be obtained from $E(1)$ by taking logarithmic transformations and taking fiber sums with additional copies of $E(1)$. We apply these constructions along regular fibers near $F_{\infty}$. Then $S_{i}$ is still embedded in $M$. Note that $$[S_{i}]-[S_{j}]=(i-j)[F_{\infty}]\neq 0\in H_{2}(M;\mathbb{Z}).$$
So $\{S_{i}\}_{i\in \mathbb{N}}$ still represent different homology classes in $M$. 

(2) Now we suppose $M$ is a hypersurface or complete intersection with $b^{+}(M)>1$. By the Lefschetz hyperplane theorem, $M$ is simply-connected. We again set $\mathfrak{s}$ to be the canonical $\spinc$ structure. Let $V=K^{\perp}\otimes \mathbb{R}$. In \cite[Page 406]{FriedmanMorgan}, the authors found a subset $\Delta\subset V$ such that any element in $\Delta$ can be represented by an embedded $(-2)$-sphere. Moreover, let $\Gamma_{\Delta}\subset \operatorname{Aut}(V)$ be the group generated by reflections along elements in $\Delta$. Then it is proved that $\Gamma_{\Delta}\cdot \Delta$ is Zariski dense in $V$. In particular, $\Gamma_{\Delta}\cdot \Delta$ is an infinite set. Since any element in  $\Gamma_{\Delta}$ can be realized by a diffeomorphism on $M$ (i.e. a composition of reflections along those $-2$-spheres), we see that any class in $\Gamma_{\Delta}\cdot \Delta$ can be presented by an embedded $-2$-sphere. 

(3) Let $M=M'\#_{T^{2}}E$. Using the Mayer-Vietoris sequence, it can be shown that $b_{1}(M)=0$ and $b^{+}(M)>1$. Since $SW(M',\mathfrak{s})\neq 0$, by the adjuntion inequality, we have $\langle c_{1}(\mathfrak{s}),[F]\rangle =0$. By the gluing formula for the Seiberg-Witten invariants along $T^{3}$ (see \cite[formula  (38.4)]{KronheimerMrowkaMonopole}), we can suitable glue $\mathfrak{s}$ and the canonical $\spinc$ structure $\mathfrak{s}_{J}$ on $E$ together and form a $\spinc$ structure $\mathfrak{s}$ such that $SW(M,\mathfrak{s})\neq 0$, $d(\mathfrak{s})=0$. Let $\{S_{i}\}$ be the $-2$-spheres embedded in $E\setminus F_{\infty}\subset M$ established in (1). Similar to (1), the homology classes of $S_{i}$ differ from each other by nonzero multiples of $F$. So $\{S_{i}\}_{i\in \mathbb{N}}$
represent distinct homology classes (up to sign) in $M$. This finishes the proof.

(4) For $i=0,1$, consider the reflection $r_{i}:M\rightarrow M$ along $S_{i}$ and their induced action $r_{i,*}$ on the subgroup $V\subset H_{2}(M;\mathbb{Z})$ generated by $[S_0]$ and $[S_{1}]$. Suppose $S_0\cdot S_{1}=k$. Then $r_{0,*}, r_{1,*}:V\rightarrow V$  are given by the matrices 
$$
\begin{pmatrix}
  1 & -k\\
  0 & -1
\end{pmatrix},\quad  
\begin{pmatrix}
  -1 & 0\\
  -k & 1
\end{pmatrix}
$$ respectively. Let $r=r_{1}\circ r_{2}$. Then $r_{*}:V\rightarrow V$
is given by the matrix 
$$
\begin{pmatrix}
  k^2-1 & -k\\
  k & -1
\end{pmatrix}.
$$
When $|k|\geq 3$, this matrix does not have an eigenvalue with norm $1$. As a result, the $(-2)$-sphere $\{r^{i}(S_{0})\}_{i\in \mathbb{N}}$ all represent different homology classes. 

(5) This follows directly from the blow up formula for the Seiberg-Witten invariants.
\end{proof}

\begin{proof}[Proof of Corollary \ref{cor: block}] (1) By \cite[Proposition 5.20]{kupers2019some}, $\pi_{1}(\operatorname{haut^{id}}(M))$ is a finitely generated abelian group. Therefore, kernel of the map $
    \pi_{1}(\Diff(M))\rightarrow \pi_{1}(\operatorname{haut^{id}}(M))$ contains a $\mathbb{Z}^{\infty}$ subgroup.

(2) We first consider the simplical group $\widetilde{\operatorname{Homeo}}(M)$ of block homeomorphisms. It is proved in \cite[Proposition 5.22]{kupers2019some} that $\pi_{1}(\widetilde{\operatorname{Homeo}}(M))$ is a finitely generated abelian group. So  kernel of the map 
\begin{equation}\label{eq: diff to block homeo}
\pi_{1}(\Diff(M))\rightarrow \pi_{1}(\widetilde{\operatorname{Homeo}}(M))
\end{equation}
contains a $\mathbb{Z}^{\infty}$ subgroup. Any element in this kernel can be represented by a loop of diffeomorphisms 
$\gamma:S^{1}\rightarrow \Diff(M)$ such that the diffeomorphism
    $$
    \widetilde{\gamma}:S^{1}\times M\rightarrow S^{1}\times M\text{ defined by } (t,x)\mapsto (t,\gamma(t)\cdot x)
    $$
    can be extended to a homeomorphism on $\widehat{\gamma}: D^{2}\times M\rightarrow D^{2}\times M$. By classical smoothing theory (see \cite{KirbySiebenmann}), the unique obstruction class of isotoping  $\widehat{\gamma}$ (relative to $S^{1} \times X$) to a diffeomorphism on $D^{2}\times M$ belongs to 
    $$
    H^{3}(D^{2}\times M, S^{1}\times M; \pi_{3}(\operatorname{Top}/\operatorname{O}))\cong H^{3}(D^{2}\times M, S^{1}\times M; \mathbb{Z}/2)\cong H^{1}(M;\mathbb{Z}/2).
    $$
Here $\operatorname{Top}$ denotes the direct limit of homeomorphism on $\mathbb{R}^{n}$ that preservers the origin and $\operatorname{O}$ denotes the direct limit of $O(n)$. Since $M$ is simply-connected, $H^{1}(M;\mathbb{Z}/2)=0$. So $\widehat{\gamma}$ is isotopic (relative to  boundary) to a diffeomorphism on $D^{2}\times M$. We have proved that the kernel of the map
\begin{equation}\label{eq: diff to block-diff}
    \pi_{1}(\Diff(M))\rightarrow \pi_{1}(\widetilde{\Diff}(M)) 
    \end{equation}
    actually coincides with the kernel of (\ref{eq: diff to block homeo}) so it contains a $\mathbb{Z}^{\infty}$ subgroup too.

(3) For any $\gamma: S^{1}\rightarrow \operatorname{Diff}(M)$ that gives an element in the kernel of (\ref{eq: diff to block-diff}), we use $\gamma$ as the clutching function and form the bundle  $\widetilde{M}_{1}/S^{2}$. Then we have a diffeomorphism from $\widetilde{M}_{1}$ to $S^{2}\times M$ obtained by gluing together $\widehat{\gamma}$ and $\operatorname{Id}_{D^{2}\times M^{2}}$. 

We assume $b^{+}(M)>1$ for now. By the proof of Theorem \ref{thm: linearly independence restated}, we may pick $\gamma$ such that  
\begin{equation}\label{eq: nonzero FSW}
a:=\FSW(\widetilde{M}_{1}/S^{2},\widetilde{\mathfrak{s}}_{1})\neq 0
\end{equation}
for some family $\spinc$ structure $\widetilde{\mathfrak{s}}_{1}$. For any $n\in \mathbb{N}$, we use a degree-$n$ map $S^{2}\rightarrow S^{2}$ to pull back $(\widetilde{M}_{1}/S^{2},\widetilde{\mathfrak{s}}_{1})$. We denote the result by $(\widetilde{M}_{n}/S^{2},\widetilde{\mathfrak{s}}_{n})$. Then by 
By Lemma \ref{lem: FSW pull-back}, we have 
\begin{equation}\label{eq: FSW for pullback}
\FSW(\widetilde{M}_{n}/S^{2},\widetilde{\mathfrak{s}}_{n})=na.
\end{equation}
Now consider another family $\spinc$ structure $\widetilde{\mathfrak{s}}$ on $\widetilde{M}_{n}/\Q$ with $d(\widetilde{\mathfrak{s}},\widetilde{M})=0$. By Lemma \ref{lem: extension of spinc b1=0}, there exists a family $\spinc$ structure on $\widetilde{M}_{1}/\Q$ such that $\widetilde{\mathfrak{s}}'|_{M}=\widetilde{\mathfrak{s}}|_{M}$. Then by Lemma \ref{lem: FSW pull-back} and Lemma \ref{lem: twisting by a line bundle doesn't change FSW}, we have 
$$
\FSW(\widetilde{M}_{n}/S^{2},\widetilde{\mathfrak{s}})=n \FSW(\widetilde{M}_{1}/S^{2},\widetilde{\mathfrak{s}}').
$$
In particular, this implies 
\begin{equation}\label{eq: FSW divisible}
n\mid \FSW(\widetilde{M}_{n}/S^{2},\widetilde{\mathfrak{s}}) \text{ for any }\widetilde{\mathfrak{s}}.    
\end{equation}
Equation (\ref{eq: FSW for pullback}) and (\ref{eq: FSW divisible}) together imply that $\widetilde{M}_{n}/S^{2}$ is not isomorphic to  $\widetilde{M}_{m}/S^{2}$ whenever $m>n|a|$. Therefore, the set $\{\widetilde{M}_{n}/S^{2}\}_{n\in \mathbb{N}}$ contains infinitely many isomorphism classes of smooth bundles over $S^{2}$. This finishes the proof when $b^{+}(M)>1$. When $b^{+}(M)=1$, we fix a choice of $H\in H^{2}(M;\mathbb{Z})$ with $H\cdot H>0$ and repeat the same argument to $\FSW_{H}$. 
\end{proof}


\bibliographystyle{amsplain}
\bibliography{Bbib0320}

\end{document}